\documentclass[12pt]{article}
\usepackage{amssymb,latexsym}
\usepackage{xcolor}
\usepackage{tikz}
\usepackage{enumitem}
\usetikzlibrary {calc}
\usetikzlibrary {decorations.pathreplacing}
\usetikzlibrary {arrows.meta} 

\usepackage{hyperref}
\def\Z{{\cal Z}}
\def\ZZ{{\mathbb Z}}

\def\RR{{\mathbb R}}

\newtheorem{formula}{}[section]
\newtheorem{proposition}[formula]{Proposition}
\newtheorem{definition}[formula]{\indent Definition}
\newtheorem{corollary}[formula]{\indent Corollary}
\newtheorem{remark}[formula]{\indent Remark}
\newtheorem{lemma}[formula]{\indent Lemma}
\newtheorem{theorem}[formula]{\indent Theorem}

\newtheorem{example}[formula]{Example}

\def\thrm{\begin{theorem}}
\def\thrml#1{\begin{theorem}\label{#1}}
\def\ethrm{\end{theorem}}
\def\rmrk{\begin{remark}}
\def\rmrkl#1{\begin{remark}\label{#1}}
\def\ermrk{\end{remark}}
\def\dfntn{\begin{definition}}
\def\dfntnl#1{\begin{definition}\label{#1}}
\def\edfntn{\end{definition}}
\def\nmrt{\begin{enumerate}}
\def\enmrt{\end{enumerate}}

\def\qtn{\begin{equation}}
\def\qtnl#1{\begin{equation}\label{#1}}
\def\eqtn{\end{equation}}
\def\lmm{\begin{lemma}}
\def\lmml#1{\begin{lemma}\label{#1}}
\def\elmm{\end{lemma}}
\def\crllr{\begin{corollary}}
\def\crllrl#1{\begin{corollary}\label{#1}}
\def\ecrllr{\end{corollary}}

\date{}
\begin{document}
\sloppy

\title{A tropical version of Hilbert polynomial (in dimension one)}

\author{Nikita Elizarov$^1$, Dima Grigoriev$^{2}$\\[3pt]
$^1$\small Fakult\"at fuer Mathematik, Universit\"at Bielefeld, 33501, Germany\\
\small \href{mailto:nelizarov@math.uni-bielefeld.de}{nelizarov@math.uni-bielefeld.de}\\
$^2$\small CNRS, Math\'ematiques, Universit\'e de Lille, Villeneuve d'Ascq, 59655, France\\
 \small   \href{mailto:Dmitry.Grigoryev@math.univ-lille1.fr}{Dmitry.Grigoryev@univ-lille.fr}\\
 \small \href{http://en.wikipedia.org/wiki/Dima\_Grigoriev}{http://en.wikipedia.org/wiki/Dima\_Grigoriev}
}
\maketitle

\begin{abstract}
For a tropical univariate polynomial $f$ we define its tropical Hilbert function as the dimension of a tropical linear prevariety of solutions of the tropical Macaulay matrix of the polynomial up to a (growing) degree. We show that the tropical Hilbert function equals (for sufficiently large degrees) a sum of a linear function and a periodic function with an integer period. The leading coefficient of the linear function coincides with the tropical entropy of $f$. Also we establish sharp bounds on the tropical entropy.
\end{abstract}

{\bf keywords:} tropical Hilbert function, tropical Macaulay matrix, tropical entropy \vspace{2mm}


\section*{Introduction}

One can find the basic concepts of tropical algebra in \cite{MS}.

Consider a tropical univariate polynomial $f:=\min_{0\le i\le n} \{iX+a_i\}$ where $a_i\in \ZZ,\, 0\le i\le n$. We call $z:=(z_1,\, z_2,\dots),\, z_j\in \RR,\, j\ge 1$ a {\it tropical recurrent sequence satisfying the vector $\vec a:=(a_0,\dots,a_n)$} \cite{G20} if for any $j\ge 1$ the following tropical (linear) polynomial is satisfied:
\begin{eqnarray}\label{-3}
\min_{0\le i\le n} \{z_{j+i}+a_i\},
\end{eqnarray}
i.~e. the minimum in (\ref{-3}) is attained at least for two different values among $0\le i\le n$.

When one considers classical recurrent sequences $(x_1,\, x_2,\dots)$ satisfying relations $\sum_{0\le i\le n} a_ix_{i+j}=0$ similar to (\ref{-3}), the first $n$ values $x_1,\dots,x_n$ determine the rest of the sequences uniquely. This is not the case for tropical recurrent sequences.

Denote by $D(k)\subset \RR^s$ a tropical linear prevariety \cite{MS} of all the sequences $(z_1,\dots,z_k)$ satisfying  (\ref{-3}) for $1\le j\le k-n$. Therefore, $D(k)$ is a polyhedral complex \cite{MS}. The function $d(k):=d_{\vec a}(k):=\dim (D(k))$ we call {\it the tropical Hilbert function} of the tropical polynomial $f$ (or equivalently, of the vector $\vec a$ of its coefficients). Obviously, $d(k)\le d(k+1)\le d(k)+1$. It is observed in \cite{G20} that $d(k+t)\le d(k)+d(t)$. Therefore, due to Fekete's subadditivity lemma \cite{S} there exists the limit
\begin{eqnarray}\label{-1}
H:=H_{\vec a}=\lim_{k\to \infty} d(k)/k
\end{eqnarray}
which is called \cite{G20} {\it the tropical entropy of the tropical polynomial $f$ or of the vector $\vec a$}. Evidently, $0\le H\le 1$.

In classical commutative algebra the Hilbert function of a polynomial $g=\sum_I g_IX^I\in F[X_1,\dots,X_m]$ is defined as the growth function of the quotient ring $F[X_1,\dots,X_m]/(g)$ in the filtration with respect to degree. For a given degree $e$ this function coincides with the dimension of the space of solutions of a linear system
\begin{eqnarray}\label{-2}
\sum_I g_I Y_{I+J}=0
\end{eqnarray}
for all  vectors $J:=(j_1,\dots,j_m)\in \ZZ^m,\, 0\le j_1,\dots,j_m$ such that for every vector $I=(i_1,\dots,i_m)$ from the support of $g$ we have $i_1+j_1+\cdots+i_m+j_m\le e$. Note that a linear system (\ref{-2}) forms the rows of Macaulay matrix.

Note that we  define  the tropical Hilbert function as the dimension of the space $D(k)$ of 
tropical recurrent sequences.
We mention that multidimensional  tropical recurrent sequences appear also as the solutions of the tropical Macaulay matrix  \cite{G20} (generalizing tropical equations (\ref{-3})). Macaulay matrix emerges in a tropical version of the weak Hilbert Nullstellensatz (see \cite{BE}, \cite{G12}, \cite{GP}, \cite{JM}, \cite{MR}, \cite{ABG}).

The main result of the paper (see Theorem~\ref{period} and Corollary~\ref{quasi-linear}) states that the tropical Hilbert function $d(k)$ is {\it quasi-linear}, i.~e. coincides (for sufficiently big $k$) with a sum $Hk+r(k)$  of a linear function $Hk$ (see (\ref{-1})) and a periodic function $r(k)$ with an integer period.

Recall that in  classical commutative algebra the Hilbert function of an ideal in $F[X_1,\dots,X_m]$ is a polynomial (for sufficiently large degrees). In its turn, the degree of this polynomial is less than $m$ (in particular, in case $m=1$ Hilbert polynomial is a constant). One can directly generalize the definition of a tropical Hilbert function to $m$-variate tropical polynomials based on the tropical Macaulay matrix for $m$-variate polynomials. This function grows asymptotically as $H\cdot k^m$ where $H$ is defined similary to (\ref{-2}) (cf. \cite{G20}).
In case of dimension $m=1$ which we study in the present paper, the tropical Hilbert function $d(k)$ coincides with the linear function $Hk$ up to a periodic function (for sufficiently large $k$).

In \cite{MR} tropical ideals in the semiring of tropical polynomials are introduced, they have some features similar to classical ideals: in particular, it is proved that Hilbert function of a tropical homogeneous ideal is eventually a  polynomial. 

We mention that in \cite{G12}, \cite{G20} it is proved that $H=0$ iff each point $(i,\, a_i)\in \RR^2,\, 0\le i\le n$ is a vertex of the Newton polygon. (The Newton polygon of $f$ is the convex hull of the rays $\{(i,\, x\ge a_i)\}, \, 0\le i\le n$.)

It would be interesting to clarify, whether one can extend the results of the paper to vectors $(a_0,\dots,a_n)$ with $a_i\in \RR \cup \{\infty\}$. Another problem is to improve the bound on the period in the function $r(k)$ and the bound on the minimal $k$ starting with which the tropical Hilbert function coincides with $Hk+r(k)$ (cf. Corollary~\ref{quasi-linear}). 

In section~\ref{one} we prove some auxiliary bounds on tropical recurrent sequences. In section~\ref{two} we describe a directed graph $G:=G_{\vec a}$ and provide a recursive construction, which allows one to produce a tropical recursive sequence corresponding to a path in $G$. Vice versa, to each tropical recurrent sequence corresponds a path in $G$. Thus, all tropical recurrent sequences corresponding to a path $T$ of a length $k$ in $G$ form a polyhedron $Q_T\subset \RR^{k+n}$. We distinguish some edges of $G$ and call them {\it augmenting}. Denote by $d(T)$ the number of augmenting edges in $T$, and by $n(T)\le n$ a certain integer depending just on the first vertex of $T$. The following theorem  (see Theorem~\ref{graph} in section~\ref{three} for more details) relates the tropical recurrent sequences with paths in $G$.

\begin{theorem}\label{sequences_paths}
i) $\dim (Q_T)=d(T)+n(T)$.

ii) The set $D(k+n)$ of all tropical recurrent sequences satisfying $\vec a$ of a length $k+n$ coincides with $\cup_T Q_T$, where $T$ ranges over all paths in $G$ of the length $k\ge 0$.

iii) The tropical Hilbert function $d_{\vec a}(k+n)=\max_T \{d(T)+n(T)\}$.
\end{theorem}

Theorem~\ref{sequences_paths} allows one to express the tropical entropy explicitly in terms of $G$. For a path $T$ in $G$ by $l(T)$ denote its length.

\begin{corollary}\label{entropy_cycle}
(see Corollary~\ref{entropy} in section~\ref{four}). The tropical entropy $H_{\vec a}$ equals the maximum of $d(T)/l(T)$ over all simple cycles $T$ in $G$. In particular, $H_{\vec a}$ is a rational number.
\end{corollary}

The main result of section~\ref{five} describes the behaviour of the tropical Hilbert function (one can find more details in Corollary~\ref{quasi-linear}). Its proof involves Theorem~\ref{sequences_paths}. 

\begin{theorem}\label{tropical_Hilbert}
The tropical Hilbert function $d_{\vec a}(k)=H_{\vec a}\cdot k +r(k)$ for sufficiently big $k$, where $r(k)$ is a periodic function with an integer period. One can compute the function $d_{\vec a}$.
\end{theorem}

In section~\ref{six} we study tropical recurrent sequences satisfying a {\it tropical boolean vector} $\vec a = (a_0,\dots,a_n)$ where $a_0=a_n=0$ and $a_i\in \{0,\infty\}, 0\le i\le n$. We establish similar results to Theorem~\ref{sequences_paths}, Corollary~\ref{entropy_cycle}, Theorem~\ref{tropical_Hilbert} for a tropical boolean vector $\vec a$.

Finally, in section~\ref{seven} we prove the following result separating a positive tropical entropy from zero (see Theorem~\ref{sharp estimate}).

\begin{theorem}
If $H_{\vec a}>0$ then $H_{\vec a}\ge 1/4$ (and the bound is sharp).
\end{theorem}
   
Also we show the sharp upper bound $H_{\vec a}\le 1-2/(n+1)$ in case when Newton polygon of $\vec a$ has a single bounded edge. We conjecture that the latter bound holds for an arbitrary vector $\vec a$.

\section{Bounds on connected coordinates}\label{one}

Let ${\vec a}:=(a_0,\dots,a_n)\in \ZZ^{n+1}$ be a vector, define its {\it amplitude} as
\begin{eqnarray}\label{0}
M:= \max_{0\le i\le n} \{a_i\}-\min_{0\le i\le n} \{a_i\}.
\end{eqnarray}

\begin{definition}
Consider a tropical recurrent sequence $z:=(z_0,\, z_1,\dots),\, z_j\in \RR$ satisfying vector $\vec a$. We call a coordinate $z_{j_0}$ (or, more precisely, $j_0$) {\it connected} if there exists $0\le t_0\le \min\{n,j_0\}$ such that $z_{j_0}+a_{t_0}=\min_{0\le t\le n}\{z_{t+j_0-t_0}+a_t\}$. In other words, one can't diminish the value of $z_{j_0}$ without changing all  other $z_j,\, j\neq j_0$ and keeping the property of being a tropical recurrent sequence satisfying $\vec a$.  Otherwise, we call $z_{j_0}$ {\it disconnected}. We say that connected coordinates $j_0<j_1$ are {\it neighbouring} if any intermediate coordinate $j_0<j<j_1$ is disconnected.
\end{definition}

\begin{lemma}\label{connected}
Let $\vec a\in \ZZ^{n+1}$ be a vector with amplitude $M$, and let $z$ be
a tropical recurrent sequence  satisfying $\vec a$.  Let $j_0<j_1$ be a pair of neighbouring connected coordinates. Then

i) $j_1-j_0\le n$;

ii) $|z_{j_0}-z_{j_1}|\le 2M$.
\end{lemma}

{\bf Proof}. To prove i) suppose the contrary. Then the minimum $\min_{0\le t\le n} \{z_{j_0+t}+a_t\}$ is attained only for $t=0$ which contradicts the assimption that $z$ satisfies $\vec a$.

To prove ii) suppose the contrary. First,  assume that $z_{j_1}\ge z_{j_0}$, hence $z_{j_1}-z_{j_0}>2M$. There exists $0\le t_1\le n$ such that 
\begin{eqnarray}\label{1}
z_{j_1}+a_{t_1}=\min_{0\le t\le n} \{z_{j_1+t-t_1}+a_t\}.
\end{eqnarray}
If $j_1-t_1\le j_0$ then $z_{j_0}+a_{j_0-j_1+t_1}<z_{j_1}-2M+a_{j_0-j_1+t_1}\le z_{j_1}-M+a_{t_1}$. and we get a contradiction with (\ref{1}), thus $j_1-t_1>j_0$.

We claim that the minimum $\min_{0\le t\le n} \{z_{j_0+t}+a_t\}$ is attained only for $t=0$. Indeed, for any connected $j_2\le j_0+n$ we have 
$$z_{j_2}+a_{j_2-j_0}\ge z_{j_1}+a_{k_1}-a_{j_2-j_1+t_1}+a_{j_2-j_0}>z_{j_0}+a_0,$$
\noindent where the first inequality is due to (\ref{1}), while the second inequality follows from $z_{j_1}-z_{j_0}>2M$ and from (\ref{0}). This proves the claim. We come to a contradiction with that $z$ satisfies $\vec a$, which completes the proof of ii) in case $z_{j_1}\ge z_{j_0}$.

The case $z_{j_1}\le z_{j_0}$ is handled in a similar way. The lemma is proved. $\Box$

\begin{corollary}\label{growth}
Let $z$ be a tropical recurrent sequence satisfying a vector $\vec a \in \ZZ^{n+1}$ of amplitude $M$, and let $j$ be a connected coordinate. Then

i)  $z_s\le z_j+2M|s-j|$ for any connected coordinate $s$:

ii)  $z_s\ge z_j-2M\cdot \max\{|s-j|,\, n\}$ for any coordinate $s$:

iii) if $z_{s+n}>\min_{s_0\le t<s+n} \{z_t\}+2Mn$ for some $s_0\ge s,\, s_0\ge 0$ then the coordinate $s+n$ is disconnected.
\end{corollary}

{\bf Proof}. i) follows immediately from Lemma~\ref{connected}~ii).

ii) follows from i)  when a coordinate $s$ is connected, moreover, in this case
\begin{eqnarray}\label{2}
z_s\ge z_j-2M|s-j|. 
\end{eqnarray}

For a disconnected coordinate $s$ one can assume w.l.o.g. that $s>j$. Take the maximal connected coordinate $s_0<s$. Lemma~\ref{connected}~i) implies that $s-s_0<n$. The minimum $\min_{0\le t\le n} \{z_{t+s-n}+a_t\}$ is attained for some connected coordinate $t_0+s-n\le s_0$. Therefore, when $t_0+s-n \ge j$, we obtain 
$$z_s+a_n\ge z_{t_0+s-n}+a_{t_0}\ge z_j-2M(t_0+s-n-j) +a_{t_0}$$
\noindent due to (\ref{2}) which proves ii) in this case.

When $ k_0+s-n < j$, we obtain
$$z_s+a_n\ge z_{t_0+s-n}+a_{t_0}\ge z_j-2M(j-t_0-s+n)+a_{t_0}\ge z_j-2M(n-1)+a_{t_0}$$
\noindent again due to (\ref{2}). This completes the proof of ii).

iii) follows from ii). $\Box$

\begin{lemma}\label{minimum_difference}
Let a vector $\vec a \in \ZZ^n$ fulfill (\ref{0}) and a sequence $(z_0,\dots,z_n)\in \ZZ^{n+1}$ satisfy $\vec a$. Denote $m_L:=\min_{0\le i<n} \{z_i\},\ s_L:=\min \{0\le i<n\ :\ i=m_L\},\ m_R:=\min_{1\le i\le n} \{z_i\},\ s_R:=\min \{1\le i\le n\ :\ i=m_R\}$. Then $|m_L-m_R|\le M,\ s_L\le s_R$.
\end{lemma}

{\bf Proof}. It suffices to consider a case when either $m_L=z_0$ or $m_R=z_n$, (otherwise, $m_L=m_R,\ s_R=s_L$). If $m_L=z_0$ then $z_0+a_0\ge z_t+a_t=\min_{0\le i\le n} \{z_i+a_i\}$ holds for a suitable $1\le t\le n$. Hence  $m_L=z_0\ge z_t+a_t-a_0\ge m_R-M$. If $m_R=z_n$, we obtain the inequality $m_R\ge m_L-M$ is a similar way. Otherwise, if $m_R=z_j,\ 1\le j<n$ then $m_L=z_0\le z_j=m_R$. $\Box$

\section{Construction of a graph of tropical recurrent sequences}\label{two}

We are producing by recursion a tropical recurrent sequence satisfying vector $\vec {a}=(a_0,\dots,a_n)$, and assume that a finite fragment (a prefix) of the sequence is already produced. Possibilities of choices of continuations of the fragment depend only on the last $n$ entries (a suffix) of the fragment. That is why we consider only the last $n$ entries and 
denote them by $(y_1,\dots, y_n) \in \RR^n$.
 We will view $(y_1,\dots, y_n)$ also as the coordinates in $\RR^n$. 
More precisely, if the minimum in $\min_{1\le i\le n} \{y_i+a_{i-1}\}$ is attained 

$\bullet$ once, then a contituation $y_{n+1}\in \RR$ is determined uniquely;

$\bullet$ at least twice, then $y_{n+1}$ ranges over an infinite interval bounded from below.

In this section we construct a directed finite graph $G:=G_{\vec a}$ and for each vertex $v$ of $G$ a polyhedron $P_v\subset \RR^n$. If $(y_1,\dots,y_n)\in P_v$, and $y_{n+1}$ is a possible continuation (i.e. the sequence $(y_1,\dots,y_n,y_{n+1})$ satisfies vector $\vec a$) then $(y_2,\dots,y_n,y_{n+1})\in P_w$ for a suitable vertex $w$ of $G$ such that $(v,w)$ is an edge of $G$. The converse is also valid: if  $(y_1,\dots,y_n)\in P_v$ and $(v,w)$ is an edge of $G$ then there exists a continuation $y_{n+1}$ such that $(y_1,\dots,y_n,y_{n+1})$ satisfies vector $\vec a$ and 
$(y_2,\dots,y_n,y_{n+1})\in P_w$.
Observe that the construction of vertices of $G$ (Definition~\ref{vertex}) depends only on $n, M$, while the construction of its edges (Definitions~\ref{single},~\ref{unbounded},~\ref{bounded}) depends also on $\vec a$.


In section~\ref{three} we show that the tropical recurrent sequences satisfying $\vec a$ are encoded by paths in $G$ (and vice versa).

\subsection{Vertices of graph $G_{\vec a}$}\label{two.1}

\begin{definition}\label{vertex}
We define a vertex $v$ of the graph $G$ and a corresponding polyhedron $P:=P_v \subset \RR^n$.
Each polyhedron $P$ is determined by the following data. We fix a subset $\emptyset \neq B:=B_v \subset \{1,\dots,n\}$, an element $s:= s_v \in B$, for each pair $1\le r<l\le n,\ r,l\in B$ integers $m(r,l):= m_v(r,l)$ and indicators $e(r,l):=e_v(r,l)\in \{0,1\}$ such that $P$ is described by the following system of linear inequalities:


\begin{eqnarray}\label{29}
y_s\le y_j,\, 1\le j\le n,
\end{eqnarray}

\begin{eqnarray}\label{30}
m(s,q)\ge -nM,\ m(p,s)\le (n+s-p)M,\ 1\le p<s<q\le n, 
\end{eqnarray} 

\begin{eqnarray}\label{30a}
|m(r,l)|\le 2nM,\ 1\le r< l\le n,\ p,q,r,l \in B,
\end{eqnarray} 

\begin{eqnarray}\label{31}
y_r=y_l+m(r,l),\ 1\le r<l\le n,\ r,l\in B\ if\ e(r,l)=0, 
\end{eqnarray}
\begin{eqnarray}\label{32}
y_r-1<y_l+m(r,l)<y_r,\ 1\le r<l\le n,\ r,l\in B\ if\ e(r,l)=1,
\end{eqnarray}
\begin{eqnarray}\label{34}
y_j-y_s>jM,\ 1\le j\le n,\ j\notin B.
\end{eqnarray}

For definiteness, we take $s\in B$ to be the minimal possible satisfying (\ref{29}).

Varying $B, s, m(r,l), e(r,l)$ we obtain all vertices of the graph $G$.
 
Coordinates $y_r$ of $\RR^n$ for $r\in B$ we call bounded on $P$, while the coordinates $y_j$ for $j\notin B$ we call unbounded.
\end{definition}




\begin{remark}\label{after_vertex}
i) The graph $G$ consists of just the vertices $v$ for which the
polyhedron $P_v$ is nonempty.  \vspace{1mm}


ii) We define an equivalence relation on $B$ setting $r,l \in B$ to belong to the same equivalence class iff $e(r,l)=0$ (i.e. the difference $y_r-y_l$ is an integer). When it is not an equivalence relation, the polyhedron $P$ is empty. 
\vspace{1mm}

iii) Informally, the difference (repectively, the ceiling function of the difference) of each pair of bounded coordinates
is given in  (\ref{31}) (respectively, in (\ref{32})), while for unbounded coordinates just lower bounds (\ref{34}) via the minimal coordinate,
which is always a bounded one, are given. 
\end{remark}

\begin{remark}\label{global_inequalities}
Note that the inequlities (\ref{30a}) follow from the inequalities (\ref{29}), (\ref{30}), (\ref{31}), (\ref{32}). We keep (\ref{30a}) to have some apriori bound on $m(k,l)$. Indeed, assume that $r<s<l,\ e(r,s)=e(s,l)=1$ (all other cases of orders between $r,s,l$ and the values of $e$ one can  study in a similar manner). Then (\ref{32}) imply that $|y_r-y_l|< \max \{m(r,s)+1, -m(s,l)\}$ taking into account that $y_s\le y_r,y_l$. Therefore, $m(r,l)\le  \max \{m(r,s)+1, -m(s,l)\}$.
\end{remark}


Below we repeatedly make use of the following statements describing the set of solutions of a system of inequalities of the form (\ref{31}), (\ref{32}).

\begin{lemma}\label{determine_minimum}
Let a system of inequalities
\begin{eqnarray}\label{100}
x_r-x_l=m(r,l),\ 1\le r<l\le n,\ e(r,l)=0,
\end{eqnarray}
\begin{eqnarray}\label{101}
m(r,l)-1<x_r-x_l<m(r,l),\ 1\le r<l\le n,\ e(r,l)=1,
\end{eqnarray}
in variables $x_1,\dots,x_n$ be consistent, where $m(r,l)\in \ZZ,\ e(r,l)\in \{0,1\}, 1\le r<l\le n$. W.l.o.g. one can assume that $x_1=0$ replacing $x_l$ by $x_l-x_1,\ 1\le l\le n$. We agree that $m(l,l)=e(l,l)=0$. Denote
\begin{eqnarray}\label{102}
\widehat{x_l}:=x_l+m(1,l),\ 1\le l\le n. 
\end{eqnarray}
Clearly, $0\le \widehat{x_l} <1, 1\le l\le n$. We claim that the system (\ref{100}), (\ref{101}) determines uniquely an order
\begin{eqnarray}\label{103}
0=\widehat{x_1}\le \widehat{x_{\pi(2)}}\le \cdots \le \widehat{x_{\pi(n)}}<1
\end{eqnarray}
for a suitable permutation $\pi \in Sym(n-1)$. The set of solutions of the system (\ref{100}), (\ref{101}) is a polytope (open in its linear hull) isomorphic to the polytope (in $\RR^n$ endowed with the coordinates $\widehat{x_1},\dots,\widehat{x_n}$) given by the system (\ref{103}). The isomorphism is assured by (\ref{102}). 

Vice versa, given $x_1,\dots, x_n$  and integers $m(1,l), 1\le l\le n$ such that (\ref{103}) is satisfied for (\ref{102}), one can find uniquely integers $m(r,l) \in \ZZ, e(r,l)\in \{0,1\}, 1\le r<l\le n$ for which (\ref{100}), (\ref{101}) hold. The integers $m(r,l), e(r,l), 1\le r<l\le n$ are determined just by the integers $m(1,l), 1\le l\le n$ and by the order  (\ref{103}).

In particular, (\ref{103}) together with (\ref{102}) allows one to find $1\le s\le n$ such that $x_s\le x_i,\ 1\le i\le n$ holds for any solution of (\ref{100}), (\ref{101}). Moreover, the set of such $s$ is determined just by the integers $m(1,l), 1\le l\le n$ and by the order  (\ref{103}).
\end{lemma}

{\bf Proof}. First note that the partition of the set $\{1,\dots,n\}$ into classes such that $r<l$ belong to the same class iff $e(r,l)=0$, provides an equivalence relation. 

It holds $-1<\widehat{x_l}-\widehat{x_r}<1,\ 1\le r\le l\le n$. If $e(r,l)=0$ then (because of (\ref{102})) we have
\begin{eqnarray}\label{105}
\widehat{x_l}-\widehat{x_r}=-m(r,l)-m(1,l)+m(1,r)\in \ZZ,
\end{eqnarray}
hence $\widehat{x_l}=\widehat{x_r}$.

For $e(r,l)=1$ we get $\widehat{x_l} \neq \widehat{x_r}$. On the other hand, due to (\ref{101}), (\ref{102}) it holds
\begin{eqnarray}\label{104}
\widehat{m(r,l)}:=m(r,l)-m(1,l)+m(1,r)-1<\widehat{x_l}-\widehat{x_r}<\widehat{m(r,l)}+1.
\end{eqnarray}
We deduce that either $\widehat{m(r,l)}=-1$ or $\widehat{m(r,l)}=0$. In the former case it holds $\widehat{x_l}<\widehat{x_r}$, in the latter case it holds $\widehat{x_l}>\widehat{x_r}$. This assures a required permutation $\pi$ satisfying (\ref{103}). Obviously, $\pi$ is independent from the particular values of $\widehat{x_l},\ 2\le l\le n$, but rather depends only on the order between them.

For any tuple $\widehat{x_l},\ 2\le l\le n$ satisfying (\ref{103}), it holds $\widehat{x_l}=\widehat{x_r}$ when $e(r,l)=0$ (see (\ref{105})), while it holds (\ref{104}) when $e(r,l)=1$. Therefore, $x_l,\ 1\le l\le n$ defined by (\ref{102}), satisfy (\ref{100}), (\ref{101}) (with $x_1=0$).   $\Box$

\begin{corollary}\label{determine_extension}
For any point $(x_1,\dots,x_{n-1})\in \RR^{n-1}$ satisfying the system (\ref{100}), (\ref{101}) (more precisely, its subsystem for $1\le r<l<n$), the set of points $x_n\in \RR$ such that the point $(x_1,\dots,x_n)$ satisfies (\ref{100}), (\ref{101}), consists of either

i) a single point when $e(r,n)=0$ for some $1\le r<n$ or

ii) an open nonempty finite interval when $e(r,n)=1$ for all $1\le r<n$.
\end{corollary}

Note that the system (\ref{31}), (\ref{32}) is similar to the system (\ref{100}), (\ref{101}).

\subsection{An edge of the graph $G$ in case of a unique continuation of a prefix of a tropical recurrent sequence}\label{two.2}

Now we describe when $G$ has an edge from a vertex $v$ to a vertex $w$. Denote the coordinates of $\RR^n$ such that the polyhedron under construction $P_w \subset \RR^n$ by $x_1,\dots,x_n$.
 The polyhedron $P_w$ relates to  $P_v$ informally as follows. For any point $(y_1,\dots,y_n)\in P_v$ there exists a point $(y_2,\dots,y_n,x_n)\in P_w$, and $x_n$ fulfills the conditions described below in Definitions~\ref{single},~\ref{unbounded},~\ref{bounded} (see Theorem~\ref{correctness} below). A value of $x_n$ is either unique or varies in an open interval. Formally, in Definitions~\ref{continuation},~\ref{single},~\ref{unbounded},~\ref{bounded} 
we describe linear inequalities determining $P_w$.

In the following definition we provide a part of equalities and inequalities of the forms (\ref{31}), (\ref{32}) describing  $P_w$ and complete the description in Definitions~\ref{single},~\ref{unbounded},~\ref{bounded}.

\begin{definition}\label{continuation}
Denote the coordinates of $\RR^n$ for which $P_w\subset \RR^n$ by $x_1,\dots,x_n$. First, we impose that  a coordinate $x_{r-1}, 2\le r\le n$ is bounded on $P_w$ iff the coordinate $y_r$ is bounded on $P_v$, in other words $B_w\supset (B\setminus \{1\})-1$ The status of boundness of the coordinate $x_n$ is specified in Definitions~\ref{single},~\ref{unbounded},~\ref{bounded}, i.e. whether $n\in B_w$. 

The description of $P_w$ contains inequalities (cf. (\ref{31}), (\ref{32}))
\begin{eqnarray}\label{93}
x_{r-1}=x_{l-1}+m(r,l),\ 2\le r<l\le n,\ r,l\in B\ if\ e(r,l)=0, 
\end{eqnarray}
\begin{eqnarray}\label{92}
x_{r-1}-1<x_{l-1}+m(r,l)<x_{r-1},\ 2\le r<l\le n,\ r,l\in B\ if\ e(r,l)=1.
\end{eqnarray}
Thus, we put $m_w(r-1,l-1):=m(r,l),\ e_w(r-1,l-1):=e(r,l),\ 2\le r<l\le n$. In addition, the description of $P_w$ contains inequalities (cf. (\ref{30a}))
\begin{eqnarray}\label{91}
|m_w(r-1,l-1)|\le 2nM.
\end{eqnarray}

\end{definition}

In Definitions~\ref{single},~\ref{unbounded},~\ref{bounded} we impose equalities and inequalities of the forms (\ref{31}), (\ref{32})  which involve the coordinate $x_n$, and also inequalities of the forms (\ref{29}), (\ref{34}) for $P_w$.


\begin{lemma}\label{invariance}
For  points $(y_1,\dots,y_n) \in P_v$ the minimum 
\begin{eqnarray}\label{3}
\min_{1\le r\le n} \{y_r+a_{r-1}\}
\end{eqnarray}
is attained on a suitable subset of the bounded coordinates $r\in B$ that do not depend on a choice of 
point 
$(y_1,\dots,y_n) \in P_v$.
\end{lemma}

{\bf Proof}. Due to (\ref{34}) the minimum in (\ref{3}) is attained only on bounded coordinates $y_r$. 
Let two points $(y_1^{(1)},\dots,y_n^{(1)}),\, (y_1^{(2)},\dots,y_n^{(2)})\in P_v$. Assume that for a pair of bounded coordinates $y_r,\, y_t$ an inequality holds $y_r^{(1)}+a_{r-1}\le y_t^{(1)}+a_{t-1}$. Then $y_r^{(2)}+a_{r-1}\le y_t^{(2)}+a_{t-1}$ because of inequalities (\ref{31}), (\ref{32}), taking into the account that $a_{r-1},\, a_{t-1}$ are integers (cf. also Lemma~\ref{determine_minimum}). $\Box$ \vspace{2mm}

\begin{definition}\label{after_invariance}
Denote by $S:=S_v(\subset B)$ the set of 
$r, 1\le r\le n$ on which the minimum in (\ref{3}) is attained.
\end{definition}
 
In particular, all the elements from $S$ belong to the same class (see Remark~\ref{after_vertex}). 
First consider the case when $S$ consists  of a single element
$t$.

\begin{definition}\label{single}
Let the set $S=\{t\}$ be a singleton. We define a unique
edge in $G$ outgoing from the vertex $v$ (to a suitable vertex $w$) and
describe a system of equations and inequalities defining a polyhedron $P_w$. 
Recall that the description of $P_w$ contains inequalities (\ref{93}), (\ref{92}), (\ref{91}).
Declare the coordinate $x_n$ to be bounded, this determines $B_w:= ((B\setminus \{1\})-1) \cup \{n\}$. 

Apply Lemma~\ref{determine_minimum} to the system (\ref{31}), (\ref{32}). We obtain a sequence of the form (\ref{103}) between $\widehat{y_i},\ 1\le i\le n,\ i\in B$. Denote $x_n:=y_t+a_{t-1}-a_n$, then $\widehat{x_n}=x_n-a_{t-1}+a_n+m(1,t)=\widehat{y_t}$ (see Lemma~\ref{determine_minimum}). Extend the obtained sequence by $\widehat{x_n}$. Again applying  Lemma~\ref{determine_minimum} to the extended sequence, we get a system of the form (\ref{100}), (\ref{101}) in the variables $y_i,\ 1\le i\le n,\ i\in B,\ x_n$. Remove from this system equations and inequalities containing $y_1$ (provided that $1\in B$) and replace $y_i,\ 2\le i\le n,\ i\in B$ by $x_{i-1}$, respectively. Thus, we obtain a system which extends the system (\ref{93}), (\ref{92}) (and playing the role of (\ref{31}), (\ref{32}) for $P_w$)
\begin{eqnarray}\label{116}
x_r=x_l+m_w(k,l),\ 1\le r<l\le n,\ r,l\in B_w\ if\ e_w(r,l)=0, 
\end{eqnarray}
\begin{eqnarray}\label{115}
x_r-1<x_l+m_w(r,l)<x_r,\ 1\le r<l\le n,\ r,l\in B_w\ if\ e_w(r,l)=1.
\end{eqnarray}  
Due to Lemma~\ref{determine_minimum} $m_w(r,l), e_w(r,l)$ for $1\le r<l<n$ coincide with the corresponding integers already constructed in Definition~\ref{continuation}.

Applying Lemma~\ref{determine_minimum} to the system (\ref{116}), (\ref{115}) we find the minimal possible $1\le  s_w\le n,\ s_w \in B_w$ such that $x_{s_w}\le x_i,\ i\in B_w$. The system  (\ref{116}), (\ref{115}) together with the following inequalities (playing the role of (\ref{29}), (\ref{30}), (\ref{30a}), (\ref{34}), respectively, for $P_w$):
\begin{eqnarray}\label{106}
x_{s_w}\le x_i,\ 1\le i\le n,
\end{eqnarray}
\begin{eqnarray}\label{107}
m_w(s_w,q)\ge -nM,\ m_w(p,s_w)\le (n+s_w-p)M, 
\end{eqnarray} 

\begin{eqnarray}\label{107a}
|m_w(r,l)|\le 2nM,\ 1\le p<s_w<q\le n,\ 1\le r< l\le n,\ p,q,r,l \in B_w,
\end{eqnarray} 
\begin{eqnarray}\label{108}
x_j-x_{s_w}>jM,\ 1\le j\le n, j\notin B_w
\end{eqnarray}
describe $P_w$.

\end{definition}

\begin{remark}\label{after_single}
Only in case $B_v=\{1\}$ the system (\ref{116}), (\ref{115}) is void, in this case $B_w=\{n\},\ s_w=n$, and $P_w$ is described by inequalities (\ref{108}) with $s_w=n$ and (\ref{106}) (being a consequence of (\ref{108})).
\end{remark}



\subsection{Edges of $G$ in case of non-uniqueness of continuations of a prefix of a tropical recurrent sequence}\label{two.3}

Now we study the case when the set $S$ (see Definition~\ref{after_invariance}) consists of more than one element. Take a minimal $t>1$ such that
$t\in S$. There can be several edges in the graph $G$ outgoing from the vertex $v$. 

\begin{definition}\label{unbounded}
First define a single edge from the vertex $v$ to a vertex $w$ such that the coordinate $x_n$ is unbounded in $P_w$, in other words we put $B_w:=(B\setminus \{1\})-1$. Again recall that the description of $P_w$ contains the inequalities (\ref{93}), (\ref{92}), (\ref{91}). 

Applying Lemma~\ref{determine_minimum} to the system (\ref{93}), (\ref{92}) one can find the minimal possible $1\le s'<n,\ s'\in B_w$ such that $x_{s'}\le x_i,\ 1\le i<n,\ i\in B_w$. We put $s_w:=s'$ (cf. (\ref{29})).
The description of $P_w$ consists of the inequalities (\ref{93}), (\ref{92}), (\ref{91}) together with the following inequalities (playing the role of (\ref{29}), (\ref{30}), (\ref{34}), respectively, for $P_w$):
\begin{eqnarray}\label{90}
x_{s_w}\le x_i,\ 1\le i\le n,
\end{eqnarray}
\begin{eqnarray}\label{89}
m(s_w,q)\ge -nM,\ m(p,s_w)\le (n+s_w-p)M 
\end{eqnarray}
for $1\le p<s_w<q\le n,\ p,q \in B_w$,
\begin{eqnarray}\label{88}
x_n-x_{s'}>nM,
\end{eqnarray}
\begin{eqnarray}\label{88a}
x_j-x_{s_w}>jM,\ 1\le j\le n,\ j\notin B_w.
\end{eqnarray}
\end{definition}

\begin{remark}\label{after_unbounded}
We distinguish  (\ref{88}) among the latter inequalities of the form (\ref{88a}) (when $j=n$) for the sake of easier references below.
\end{remark}

The constructed vertex $w$ is the unique one to which there is an edge in the graph $G$ from the vertex $v$  such that the coordinate $x_n$ is unbounded. Still we assume that $|S|\ge 2, t\in S$ with a minimal possible $t>1$. Now we construct vertices $w$ with a bounded coordinate $x_n$ to which there are edges from $v$.

\begin{definition}\label{bounded}
We declare the coordinate $x_n$ to be bounded, i.e. $B_w:=((B\setminus \{1\})-1) \cup \{n\}$.
 Recall that the description of $P_w$ already contains the equalities and inequalities (\ref{93}), (\ref{92}), (\ref{91}). As in Definition~\ref{unbounded} applying Lemma~\ref{determine_minimum} to the system (\ref{93}), (\ref{92}) one can find the minimal possible $1\le s'<n$ such that $x_{s'}\le x_l,\ 1\le l<n,\ l\in B_w$. 
 
We choose all possible integers $m_w(l,n),\ 0\le e_w(l,n)\le 1,\ 1\le l<n,\ l\in B_w$ for which it holds
\begin{eqnarray}\label{109}
m_w(t-1,n)\le a_n- a_{t-1} - e_w(t-1,n),
\end{eqnarray}
\begin{eqnarray}\label{109a}
m_w(s',n)\ge -nM,
\end{eqnarray} 
\begin{eqnarray}\label{109b}
|m_w(l,n)|\le 2nM,\ 1\le l<n,\ l\in B_w.
\end{eqnarray}   
      
The description of $P_w$ contains inequalities  (\ref{109b}) playing the role of (\ref{30a}) (a part of them are inequalities (\ref{91})) and the following inequalities playing the role of (\ref{31}), (\ref{32}) (a part of them are inequalities (\ref{93}), (\ref{92})):
\begin{eqnarray}\label{110}
x_r-x_l=m_w(r,l),\ 1\le r<l\le n,\ r,l \in B_w,\ e_w(r,l)=0,
\end{eqnarray}
\begin{eqnarray}\label{111}
m_w(r,l)-1<x_r-x_l<m_w(r,l),\ 1\le r<l\le n,\ r,l \in B_w,\ e_w(r,l)=1.
\end{eqnarray}

Applying Lemma~\ref{determine_minimum} to the system (\ref{110}), (\ref{111}) one can find the minimal possible $1\le s_w \le n,\ s_w\in B_w$ such that $x_{s_w}\le x_l,\ 1\le l\le n,\ l\in B_w$. The description of $P_w$ contains the following inequalities (playing the role of (\ref{29}), (\ref{30}), (\ref{34}), respectively):
\begin{eqnarray}\label{112}
x_{s_w}\le x_i,\ 1\le i\le n,
\end{eqnarray}
\begin{eqnarray}\label{113}
m_w(s_w,q)\ge -nM,\ m_w(p,s_w)\le (n+s_w-p)M, 
\end{eqnarray}
for $1\le p<s_w<q\le n,\ p,q\in S_w$,
\begin{eqnarray}\label{114}
x_j-x_{s_w}>jM,\ 1\le j\le n,\ j\notin B_w.
\end{eqnarray}
Thus, the description of $P_w$ consists of the inequalities (\ref{110}) - (\ref{114}) for all possible choices of integers $m_w(l,n),\ 0\le e_w(l,n)\le 1,\ 1\le l<n,\ l\in B_w$ satisfying (\ref{109}), (\ref{109a}), (\ref{109b}), provided that $P_w$ is not empty.

\end{definition}

\begin{remark}\label{after_bounded}
In Definition~\ref{bounded} it holds either $s_w=s'$ or $s_w=n$ (the latter holds iff $x_n<x_{s'}$).
\end{remark}

This completes the description of all the edges outgoing  from the vertex $v$ in the graph $G$.

\section{Description of tropical recurrent sequences via paths in the graph}\label{three}

\subsection{Producing a short tropical recurrent sequence along an edge of the graph}\label{three.5}

In this subsection  for any point $(y_1,\dots,y_n)\in P_v$
we prove the following claim. If a sequence $(y_1,\dots,y_n,x)\in \RR^{n+1}$ satisfies the vector $\vec a$ then for exactly one of the edges $(v,w)$ of the graph $G$ it holds that $(y_2,\dots,y_n,x)\in P_w$. Conversely, for every edge $(v,w)$
of $G$ constructed according to one of Definitions~\ref{single}, \ref{unbounded}, \ref{bounded} there exists a point $(y_2,\dots,y_n,x_n)\in P_w$ such that the point $(y_1,\dots,y_n,x_n)\in \RR^{n+1}$ satisfies the vector $\vec a$ (for more precise statements see Theorem~\ref{correctness}).


\vspace{1mm}
We assume that a point $(y_1,\dots,y_n,x_n)\in \RR^{n+1}$ satisfies the vector $\vec a$.
Denote $x_i:= y_{i+1}, 1\le i\le n-1$ (see Definition~\ref{continuation}). According to Definition~\ref{continuation} a coordinate $y_{i+1}, 1\le i<n$ is bounded on $P_v$ (i.e. $i+1\in B$) iff the coordinate $x_i$ is bounded on $P_w$ (i.e. $i\in B_w$).
Then the bounded coordinates among $x_1,\dots, x_{n-1}$ fulfill the inequalities  (\ref{93}), (\ref{92}), (\ref{91}) introduced in Definition~\ref{continuation}.
\vspace{1mm}

Consider the case of a singleton $S=\{t\}$. Then $x_n=y_t+a_{t-1}-a_n$. We claim that $(y_2,\dots,y_n, x_n)\in P_w$ where the edge $(v,w)$ of $G$ is constructed according to Definition~\ref{single}. Recall that $B_w=((B\setminus \{1\})-1) \cup \{n\}$ in this case. The inequalities (\ref{106}), (\ref{116}), (\ref{115}) are fulfilled by the construction in Definition~\ref{single}.

Now we verify (\ref{107}). First assume that $1\le s_w<n$. Since $s\le s_w+1$ due to Lemma~\ref{minimum_difference}, the inequalities (\ref{107}) for $1\le p<q<n$ follow from (\ref{30}) taking into account that $y_s\le y_{s_w+1}=x_{s_w}$. It holds $x_n=y_t+a_{t-1}-a_n\le x_{s_w}+a_{s_w}-a_n\le x_{s_w}+M$, hence $m_w(s_w,n)\ge -M$ (taking into account the inequalities (\ref{116}), (\ref{115})), which justifies (\ref{107}) in case $1\le s_w<n$. Now assume that $s_w=n$. Lemma~\ref{minimum_difference} implies that $x_n\ge y_s-M$, and (\ref{107}) follows from (\ref{30}).

Now we verify the inequalities (\ref{108}). Recall that $j+1 \notin B, 1\le j<n$ iff $j\notin B_w$. The inequality $y_{j+1}-y_s>(j+1)M, j+1\notin B$ (see (\ref{34})) implies that $x_j-x_{s_w}>jM$ since $x_{s_w}\le y_s+M$ (due to Lemma~\ref{minimum_difference}). This justifies the inequalities (\ref{108}). Thus, the  point $(x_1,\dots,x_n)$ belongs to the polyhedron $P_w$ for an edge $(v,w)$ of $G$ constructed according to Definition~\ref{single}. This proves the claim in case $S=\{t\}$. \vspace{1mm}
 


Now we study the case when $|S|\ge 2$ and the inequality (\ref{88}) is true. We claim that in this case $(x_1,\dots,x_n)\in P_w$ where the edge  $(v,w)$ of $G$ is  constructed according to Definition~\ref{unbounded}. Recall that in this case $B_w=(B\setminus \{1\})-1$, and $2\le t\le n$ is the minimal element of $S\setminus \{1\}$. The inequalities (\ref{90}) are fulfilled according to the construction in Definition~\ref{unbounded}. The inequalities (\ref{89}) follow from (\ref{30}) taking into account that $s\le s_w+1$ due to Lemma~\ref{minimum_difference} and that $y_s\le y_{s_w+1}=x_{s_w}$ (cf. the similar argument in the case $|S|=1$ above). The inequalities (\ref{88a}) for $1\le j<n$ are justified as above in case $|S|=1$. The inequality (\ref{88a}) for $j=n$ coincides with (\ref{88}). Thus, the point $(x_1,\dots,x_n)$ belongs to the polyhedron $P_w$ for an edge $(v,w)$ construced in Definition~\ref{unbounded}. This proves the claim when $|S|\ge 2$ and (\ref{88}) holds.  \vspace{1mm}


Now we assume that $|S|\ge 2$ and (\ref{109a}) hold (in other words, (\ref{88}) is not true). We claim that $(x_1,\dots,x_n)$ belongs to the polyhedron $P_w$ where the edge  $(v,w)$ of $G$ is  constructed according to Definition~\ref{bounded}. Recall that in this case we have $B_w=((B\setminus \{1\})-1) \cup \{n\}$. For the minimal element $2\le t\le n$ of $S\setminus \{1\}$ it holds $x_n\ge x_{t-1}+a_{t-1}-a_n$. The inequalities (\ref{110}), (\ref{111}) for suitable $m_w(r,n), e_w(r,n), 1\le r<n, r\in B_w$ are fulfilled according to the construction in Lemma~\ref{determine_minimum} which we apply to $x_1,\dots,x_n$. Then the inequality $x_n\ge x_{t-1}+a_{t-1}-a_n$ implies (\ref{109}).

We verify the inequalities (\ref{113}) (for $m_w(r,l)$ constructed in Definition~\ref{bounded} invoking Lemma~\ref{determine_minimum}). Observe that it holds either $s_w=s'$ or $s_w=n$ (see Remark~\ref{after_bounded}). First assume that $s_w=s'$. The equalities (\ref{113}) for $q<n$ follow from the inequalities (\ref{30}) taking into account that $s\le s'+1=s_w+1$ and that $y_s\le y_{s'+1}=x_{s'}=x_{s_w}$ (cf. the similar argument in the consideration of the case $|S|=1$ above). The inequality (\ref{113}) for $q=n$ follows from (\ref{109a}) (taking into account the inequalities (\ref{110}), (\ref{111})). Now assume that $s_w=n$. Lemma~\ref{minimum_difference} implies that $x_n\ge y_s-M$, and therefore the inequalities (\ref{113}) follow from (\ref{30}) (cf. the similar argument in the consideration of the case $|S|=1$ above). So, the inequalities   (\ref{113}) are justified.

The inequalities (\ref{114}) we verify as in the consideration of the case $|S|=1$ above. The inequalities (\ref{109b}) follow from the inequalities (\ref{113}) (see the Remark~\ref{global_inequalities}). Thus, the point $(x_1,\dots,x_n)$ belongs to a polyhedron $P_w$ for an appropriate edge $(v,w)$ of $G$ constructed in Definition~\ref{bounded}. This completes the proof of the claim. \vspace{2mm}
 


Conversely, assume that for a point $(x_1,\dots,x_n)\in P_v$
it holds $(x_1,\dots,x_n):=(y_2,\dots,y_n, x_n)\in P_w$ for an edge  $(v,w)$ of $G$ constructed according to one of Definitions~\ref{single}, \ref{unbounded}, \ref{bounded}. First, we study the case

i) there exists $t\in S,\, 2\le t\le n$. If $S=\{t\}$ then $x_n= x_{t-1}+a_{t-1}-a_n$ (see Definition~\ref{single}). Otherwise, if $|S|\ge 2$ then $x_n\ge x_{t-1}+a_{t-1}-a_n$ (see (\ref{88}) in case of Definition~\ref{unbounded} and (\ref{109}) in case of Definition~\ref{bounded}). Therefore, the point $(y_1,\dots,y_n, x_n)\in \RR^{n+1}$ satisfies the vector $\vec a$. 

Observe that if the edge $(v,w)$ is constructed according to Definition~\ref{unbounded} then the values of the coordinate $x_n$ vary in an open infinite interval bounded from below (see (\ref{88})). If the edge $(v,w)$  is constructed according to Definition~\ref{bounded}, and the description of $P_w$ contains an equality $x_l-x_n=m_w(l,n)$ of the form (\ref{110}) for some $1\le l\le n-1$ then the value of the coordinate $x_n$ is unique. Otherwise, if $e_w(l,n)=1, 1\le l<n, l\in B_w$, the values of the coordinate $x_n$ vary in an open finite interval due to Corollary~\ref{determine_extension}. \vspace{1mm}

ii) Now assume that $S=\{1\}$, then the point $(x_1,\dots,x_n)\in P_w$ for an edge $(v,w)$ constructed according to Definition~\ref{single}. If $e(1,l_0)=0$ for some $2\le l_0\le n, l_0\in B$, i.e. the description of $P_v$ contains an equality $y_1=y_{l_0}+m(1,l_0)$ of the form (\ref{31}), then the description of $P_w$ contains the equality $x_n+a_n-a_0=x_{l_0-1}+m(1,l_0)$ (see (\ref{116})). Hence in this case for any point $(y_2,\dots,y_n, x_n)\in P_w$ the point $(y_1,\dots,y_n, x_n) \in \RR^{n+1}$ satisfies the vector $\vec a$ (in fact, $(y_2,\dots,y_n, x_n)\in P_w$ implies that $x_n=y_{l_0}+a_0-a_n+m(1,l_0)$).  

Otherwise, if $e(1,l)=1, 2\le l\le n, l\in B$ then the values  of the coordinate $x_n$ such that $(y_2,\dots,y_n, x_n)\in P_w$ vary in an open  interval (in this case $e_w(r,n)=1, 1\le r<n, r\in B_w$, see (\ref{115})). Observe that the latter interval is finite iff (\ref{115}) is not void, i.e. $e(1,l)=1$ for some $2\le l\le n, l\in B$, in other words $B\setminus \{1\} \neq \emptyset$. If $B=\{1\}$ then this interval is infinite bounded from above (see (\ref{108}) and Remark~\ref{after_single}). In this case only the point $(y_1,\dots,y_n,x_n=y_n+a_0-a_n)$ satisfies the vector $\vec a$. 

Summarizing, we have proved the following theorem.

\begin{theorem}\label{correctness}
Let a point $(y_1,\dots,y_n)\in P_v$. 

If a point $(y_1,\dots,y_n, x_n)\in \RR^{n+1}$ satisfies the vector $\vec a$ then $(y_2,\dots,y_n, x_n)\in P_w$ holds for exactly one edge $(v,w)$ of the graph $G$ constructed according to Definitions~\ref{single}, \ref{unbounded}, \ref{bounded}.
\vspace{1mm}

Conversely, let $(y_2,\dots,y_n, x_n)\in P_w$ for an edge $(v,w)$ of  $G$ constructed according to one of Definitions~\ref{single}, \ref{unbounded}, \ref{bounded}.

i) In case when there exists $t\in S,\, 2\le t\le n$ (see subsection~\ref{two.2}) the point $(y_1,\dots,y_n, x_n)\in \RR^{n+1}$ satisfies the vector $\vec a$. In case of an edge constructed according to

$\bullet$ Definition~\ref{single}, the value of $x_n$ is unique;

$\bullet$ Definition~\ref{unbounded}, the values of $x_n$ vary in an open infinite interval bounded from below;

$\bullet$ Definition~\ref{bounded}, the values of $x_n$ depending on the edge $(v,w)$, can be either unique or vary in an open finite interval. \vspace{1mm}

ii) If $S=\{1\}$ then only for the value $x_n=y_1+a_0-a_n$ the point $(y_1,\dots,y_n, x_n)\in \RR^{n+1}$ satisfies the vector $\vec a$.    

\end{theorem}

\subsection{The polyhedron of tropical recurrent sequences produced along a path of the graph}\label{three.2}

We consider paths in the graph $G$ and describe how they correspond to the tropical recurrent sequences satisfying the vector $\vec a$.  Take an arbitrary vertex $v_0$ as the first vertex in a path and any sequence $y^{(0)}:=(y_1^{(0)},\dots,y_n^{(0)})\in P_{v_0}$. As in subsection~\ref{two.2} consider a subset $S$. If $|S|=1$ then there is a unique edge $(v_0,w_0)$ in $G$ outgoing from $v_0$. In this case one applies Definition~\ref{single}  and obtains a unique $y_{n+1}^{(0)}:=x_n^{(0)}\in \RR$ such that $(y_2^{(0)},\dots,y_{n+1}^{(0)})\in P_{w_0}$ and $(y_1^{(0)},\dots,y_{n+1}^{(0)})$ satisfies vector $\vec a$ (see Theorem~\ref{correctness}). 

Otherwise, if $|S|>1$ then there are several edges outgoing from $v_0$. For each edge  $(v_0,v_1)$ one applies either Definition~\ref{unbounded} or Definition~\ref{bounded}, respectively, and produces $y_{n+1}^{(0)}:=x_n^{(0)}\in \RR$ such that $(y_2^{(0)},\dots,y_{n+1}^{(0)})\in P_{v_1}$ and $(y_1^{(0)},\dots,y_{n+1}^{(0)})$ satisfies vector $\vec a$ (see Theorem~\ref{correctness}). Recall (see Theorem~\ref{correctness}~i)) that for certain edges $(v_0,v_1)$ the value $y_{n+1}^{(0)}$ is unique, while for other edges $y_{n+1}^{(0)}$ runs over an open interval. 

An edge  $(v_0,v_1)$ for which the value $y_{n+1}^{(0)}$ is unique we call {\it rigid}, otherwise if the values run over an open interval we call an edge {\it augmenting}. Due to Theorem~\ref{correctness}~i) the property of an edge to be rigid or augmenting does not depend on a point $y^{(0)}$. Note that in case  of $S$ being a singleton, the edge is rigid, while an edge constructed according to Definition~\ref{unbounded}, is augmenting (cf. also Theorem~\ref{correctness}~ii)).

So far, we have produced a short tropical recurrent sequence $(y_1^{(0)},\dots,y_{n+1}^{(0)})$ corresponding to an edge of $G$. We treat this as a base of recursion. Suppose that we have produced by recursion a tropical recurrent sequence $(y_1^{(0)},\dots, y_{n+k}^{(0)})$ satisfying the vector $\vec a$ corresponding to a path $T$ of the length $k$ in $G$ (the length of a path is defined as the number of its edges). Let $v$ be the last vertex of $T$. Then we apply to $v$ and to the suffix $(y_{k+1}^{(0)},\dots, y_{n+k}^{(0)})$ of the produced sequence one of Definitions~\ref{single}, \ref{unbounded}, \ref{bounded} as above in the base of recursion, choosing an edge $(v,w)$ of $G$ and producing $y_{n+k+1}^{(0)}$. Thereby, we get a tropical recurrent sequence $(y_1^{(0)},\dots, y_{n+k+1}^{(0)})$ satisfying the vector $\vec a$ and corresponding to the path $T_w$ obtained by extending $T$ by an edge $(v,w)$. This completes the recursive step.
Summarizing, we have established in this subsection the following proposition.

\begin{proposition}\label{path}
For any path in the graph $G$
any produced (by the described recursive process) sequence along this path is a tropical recurrent sequence satisfying vector $\vec a$.
\end{proposition}

Denote by 
$Q_T\subset \RR^{k+n}$
a polyhedron of all the tropical recurrent sequences
which are produced along the path $T$ as described above (see Proposition~\ref{path}). Thus, any produced tropical recurrent sequence satisfies the vector $\vec a$.
The polyhedron $Q_T$
is presented by the systems of linear equations and linear  inequalities produced in Definitions~\ref{single}, \ref{unbounded}, \ref{bounded}, respectively, applied to the edges of the path $T$ (see Theorem~\ref{correctness}). Observe that when $S\neq \{1\}$ Theorem~\ref{correctness}~i) implies that for the inequalities describing $Q_T$ just the inequalities describing $P_v$ and $P_w$ suffice, while when $S= \{1\}$ one has to add to the latter inequalities also the equality $x_n=y_1+a_0-a_n$ (see Theorem~\ref{correctness}~ii)). 






Observe that for a rigid edge $(v,w)$ the polyhedron 
$Q_{T_w}\subset \RR^{k+n+1}$ 
is homeomorphic to  $Q_T$, and the homeomorphism is provided by the projection along the last coordinate. For an augmenting edge  $(v,w)$ the polyhedron $Q_{T_w}$ is homeomorphic to the cylinder $Q_T\times \RR$. In particular, in the latter case $\dim (Q_{T_w})= \dim (Q_T) + 1$. Summarizing, we have established the following proposition.

\begin{proposition}\label{polyhedron}
Let $T$ be a finite path of the graph $G$ with an ending vertex $v$,
and $T_w$ be an extension of $T$ by an edge  $(v,w)$. If the edge $(v,w)$ is rigid then the polyhedron $Q_{T_w}$ of all the finite tropical recurrent sequences 
produced along $T_w$ (see Proposition~\ref{path})
 is homeomorphic to the polyhedron $Q_T$, while if $(v,w)$ is augmenting then $Q_{T_w}$ is homeomorphic to the cylinder $Q_T\times \RR$.
\end{proposition} 

\subsection{Completeness of the construction of tropical recurrent sequences}

Now, conversely to Proposition~\ref{path}, we claim that every tropical recurrent sequence $y:=(y_1,y_2,\dots)$ satisfying the vector $\vec a$ emerges along an appropriate path of the graph $G$ (see subsection~\ref{three.2}). 

\begin{proposition}\label{sequence}
i) The union of the polyhedra $P_v$ over all the vertices $v$ of the graph $G$ coincides with $\RR^n$. \vspace{1mm}

ii) For any tropical recurrent sequence $y:=(y_1,y_2,\dots)$ satisfying the vector $\vec a$ and a vertex $v$ of $G$ such that $(y_1,\dots,y_n)\in P_v$ there exists a unique path $T$ of $G$ starting with $v$ such that $y$ is produced along $T$ as described in subsection~\ref{three.2} (see Proposition~\ref{path}).

\end{proposition}

{\bf Proof}. We prove ii) by recursion.
For the base of recursion assume that $y_{s_0}:=\min_{1\le j\le n} \{y_j\}$ holds for the minimal possible $1\le s_0\le n$. We construct a vertex $v_0$ of $G$ such that $(y_1,\dots,y_n)\in P_{v_0}$ as follows. Put $B_{v_0}:= \{1\le j\le n\ :\ y_j-y_{s_0}\le jM\}$ (cf. (\ref{34})). Applying Lemma~\ref{determine_minimum} to the set $\{y_l\ :\ l\in B_{v_0}\}$, one produces integers $m_{v_0}(r,l), e_{v_0}(r,l);\ 1\le r<l\le n;\ r,l\in B_{v_0}$. Then the inequalities similar to (\ref{29}) - (\ref{34}) describe the polyhedron $P_{v_0}$ such that $(y_1,\dots,y_n)\in P_{v_0}$.
 In particular, this proves i).

For the recursive step suppose that a path $T$ of $G$ of a length $k$ is already constructed such that the sequence $(y_1,\dots,y_{n+k})$ is produced along $T$ as in subsection~\ref{three.2} (see Proposition~\ref{path}). Let $v$ be the last vertex of $T$, then  $y^{(k)}:=(y_{k+1},\dots,y_{k+n})\in P_v$. Apply Theorem~\ref{correctness} to $y^{(k)}$, this provides a unique edge $(v,w)$ of $G$ such that $(y_{k+2},\dots,y_{k+n+1})\in P_w$, thus 
the sequence $(y_1,\dots,y_{n+k+1})$ is produced along the extended path $T_w$.
This completes the proof (by recursion) of ii). $\Box$ \vspace{2mm}


Observe that one could choose, perhaps, another initial vertex $v'$ of $G$ such that  $(y_1,\dots,y_n)\in P_{v'}$ (the latter inclusion is the only property of $v'$ we require). In fact, one could declare (in an arbitrary way) any coordinate $y_j, 1\le j\le n$  either bounded on $P_{v'}$ (i.e. $j\in B_{v'}$) or unbounded (i.e. $j\notin B_{v'}$) if it fulfills the inequalities either $jM<y_j-y_{s_0},\ \lfloor y_j-y_{s_0}\rfloor \le (n+s_0-j)M$ when $1\le j<s_0$ or $jM<y_j-y_{s_0}\le nM$ when $s_0<j\le n$ (cf. (\ref{30}), (\ref{34})). Observe that if $y_j-y_{s_0}\le jM$ then $y_j$ should be bounded on $P_{v'}$, i.e. $j\in B_{v'}$ (cf. (\ref{34})), while if either $y_j-y_{s_0}\ge (n+s_0-j)M$ when $1\le j<s_0$ or $y_j-y_{s_0}> nM$ when $s_0<j\le n$, then $y_j$ should be unbounded on $P_{v'}$, i.e. $j\notin B_{v'}$ (cf. (\ref{30})).

 After choosing an initial vertex $v_0$, the rest of a path $T$ in $G$ is constructed uniquely (see Theorem~\ref{correctness} and subsection~\ref{three.2}). Therefore, each tropical recurrent sequence satisfying the vector $\vec a$ corresponds to just a finite number of paths in $G$ as in subsection~\ref{three.2} (see Proposition~\ref{path}). Moreover, this number does not exceed the number of vertices in $G$. Thus, the tropical prevariety of all the tropical recurrent sequences of a length $n+k$ satisfying the vector $\vec a$ has the same dimension as the union of polyhedra $Q_T$ over all the paths $T$ of the length $k$ in $G$.    

For a path $T$ in the graph $G$ denote by $d(T)$ the number of augmenting edges
 in $T$. By $n(T)\le n$ denote the number of (equivalence) classes of the coordinates in the first vertex of $T$ (see Remark~\ref{after_vertex}). 
Thus, we have established
the following theorem taking into account Propositions~\ref{path}, \ref{polyhedron}, \ref{sequence}.

\begin{theorem}\label{graph}
For any vector $\vec a:=(a_0,\dots,a_n)\in \ZZ^{n+1}$ 
the constructed finite directed graph $G:=G_a$ satisfies the following properties. 
For a path $T$ 
of a length $k$ 
in $G$ 
denote by $Q_T\subset \RR^{k+n}$ the polyhedron of all the tropical recurrent sequences satisfying the vector $\vec a$ and being produced along the path $T$ in $G$. Then $\dim (Q_T)=d(T)+n(T)$.
Moreover, the union of polyhedra $Q_T$ over all the paths $T$ of the length $k$ 
coincides with the tropical prevariety of all the tropical recurrent sequences of the length $k+n$ satisfying the vector $\vec a$.
\end{theorem}

\section{Calculating the entropy via the graph of tropical recurrent sequences}\label{four}

In this section we study the tropical Hilbert function $d(k):=d_{\vec a}(k)$ (see the Introduction). Due to Theorem~\ref{graph} $d(k)$ equals the maximum of $n(T) +d(T)$ over all the paths $T$ of the length $k-n$ in the graph $G$.

We call a simple cycle in $G$ {\it optimal} if the quotient of the number of augmenting edges in the cycle to the length of the cycle is the maximal among the simple cycles.  This maximal quotient we denote by ${\cal H}:={\cal H}_{\vec a}$. Later we show that $\cal H$ equals the entropy $H:=H_{\vec a}$ (Corollary~\ref{entropy}). Clearly, $\cal H$ equals the maximum of the same quotient over all the cycles in $G$ (not necessary, simple).

First, we prove a lower bound on the tropical Hilbert function $d(k)$.

\begin{lemma}\label{lower}
$d(k)\ge {\cal H} (k-n)$.
\end{lemma}

{\bf Proof}. Take an optimal simple cycle $U$ in $G$. Denote the length of $U$ by $L$ and the number of augmenting edges in $U$ by $m$, then ${\cal H}=m/L$. Assign to each augmenting edge of $U$ the number $1-{\cal H}$ and to each rigid edge the number $-\cal H$. Then the sum of all these numbers equals 0. Due to the lemma about  leaders \cite{R} there exists a vertex $u$ of $U$ such that the sum of the assigned numbers along any subpath of $U$ starting with $u$, is non-negative.

Consider a path $T$ of a length $k-n$ starting with the vertex previous to $u$ in $U$ and following the cycle $U$ (i.~e. $T$ can wind the cycle $U$ several times). According to Theorem~\ref{graph} $\dim (Q_T) \ge {\cal H}(k-n)$. $\Box$ \vspace{2mm}

Denote by $V$ the number of  vertices in $G$. Now we proceed to an upper bound on the tropical Hilbert function.  

\begin{lemma}\label{upper}
$d(k)\le {\cal H}k+(1-{\cal H})(V+n)$.
\end{lemma}

{\bf Proof}. Consider a path $T$ of a length $L$ in $G$. Take a vertex $v_1$ of $G$ which occurs in $T$ at least twice (provided that it does exist). Then the subpath of $T$ between these two occurrings constitues a cycle of a length $L_1$. Remove this cycle from $T$, and continue removing cycles from the resulting paths, while it is possible. Let $L_2,L_3,\dots,L_q$ be the lengths of the consecutively removed cycles. Then 
$$d(T)\le {\cal H}(L_1+\cdots+L_q)+(L-L_1-\cdots-L_q)\le {\cal H}(L_1+\cdots+L_q)+V$$
\noindent (cf. Theorem~\ref{graph}). Therefore, $d(k)\le {\cal H}(k-n)+(1-{\cal H})V+n$ taking into the account that $L-L_1-\cdots-L_q\le V$. $\Box$ \vspace{2mm}

Lemmata~\ref{lower}, \ref{upper} imply the following corollary (see (\ref{-1})).

\begin{corollary}\label{entropy}
${\cal H}=H$.
\end{corollary}

\begin{remark}\label{rational}
The entropy $H$ is a rational number.
\end{remark}

\section{Quasi-linearity of the tropical Hilbert function}\label{five}

\begin{lemma}\label{optimal}
Any path $T$ of a length $k-n$ greater than $V^2(V+n)+V$ in the graph $G$ such that $n(T)+d(T)=d(k)$, contains a vertex from an optimal simple cycle.
\end{lemma}

{\bf Proof}. First consider the case when $H=0$. Then any simple cycle in $G$ is optimal, and the statement of the lemma is true even with a better bound $k-n>V$. Thus, from now on in the proof of the lemma  we assume that $H>0$.

Recall that according to Theorem~\ref{graph} it holds that $\dim(Q_T)=n(T)+d(T)$. Slightly modifying the construction from the proof of Lemma~\ref{upper}, take the first repetition of some vertex $v$ in $T$ (provided that it is possible). Then the subpath of $T$ between these two occurrences of $v$ constitutes a simple cycle of a length $L_1$ in $T$. Remove this simple cycle from $T$ and continue removing simple cycles from the resulting paths in a similar way, while it is possible. Denote by $L_2,L_3,\dots,L_q$ the lengths of the consecutively removed simple cycles. Denote by $B$ the denominator of $H$ (cf. Remark~\ref{rational}), obviously $B\le V$ (see section~\ref{four}).

Assume the contrary to the claim of the lemma. Then 
$$d(T)\le H(L_1+\cdots+L_q)-q/B+(k-n-L_1-\cdots-L_q)\le H(k-n)-q/B+V.$$
\noindent 
The first inequality follows from the statement that the amount of augmenting edges in the cycle with the length $L_i,\, 1\le i\le q$ is not greater than $H\cdot L_i - \frac{1}{B}.$ 
Making use of Lemma~\ref{lower} we obtain an inequality $q/B\le V+n$, hence $q\le V(V+n)$. 
The path $T$ consists of $q$ simple cycles and a path without cycles. Each simple cycle has length not more than $V$ as well as the path without cycles. 
Therefore, $k-n\le V^2(V+n)+V$ since $L_1,\dots,L_q\le V$. $\Box$ \vspace{2mm}

Denote by $R$ the least common multiple of the lengths of all the optimal simple cycles.

\begin{lemma}\label{shift}
For any $k>(V^2+1)(V+n)$ we have $d(k+R)\ge d(k)+HR$.
\end{lemma}

{\bf Proof}. Take a path $T$ of the length $k-n$ in $G$ such that $n(T)+d(T)=d(k)$ (cf. Theorem~\ref{graph}). Due to Lemma~\ref{optimal} $T$ contains a vertex $v$ which belongs to an optimal simple cycle $C$ of a length $c$. Glue in the path $T$ at the vertex $v$ the number $R/c$ of copies of the cycle $C$, the resulting path of the length $k-n+R$ denote by $T_1$. In other words, in $T_1$ one follows first $T$ till the vertex $v$, then there are $R/c$ windings of the cycle $C$ (finishing at $v$), finally after that one again follows path $T$ (starting at $v$). Clearly, $d(T_1)=d(T) +(R/c)Hc$. $\Box$

\begin{lemma}\label{stable}
If for some $k>(V^2+1)(V+n)$ we have 
$$d(k+iR)=d(k)+HiR,\, 0\le i\le V((1-H)V+n+1)$$ 
\noindent then $d(k+jr)=d(k)+HjR$ for any $j\ge 0$.
\end{lemma}

{\bf Proof}. Due to Lemma~\ref{shift}  it holds $d(k+jR)\ge d(k)+HjR$. Suppose that 
\begin{eqnarray}\label{9}
d(k+jR) >d(k)+HjR
\end{eqnarray}
for some $j>V((1-H)V+n+1)$, and take the minimal such $j$. There exists a path $T$ of the length $k+jR-n$ in $G$ for which $n(T)+d(T)=d(k+jR)$. For $0\le i\le V((1-H)V+n+1)$ denote by $T_i$ the beginning of   $T$ of the length $k+iR-n$. One can represent the path $T=T_i\overline{T_i}$ as a concatenation of two paths.

There exists a subsequence $0\le i_0<i_1<\cdots<i_{(1-H)V+n+1}\le V((1-H)V+n+1)$ such that each path $T_{i_l},\, 0\le l\le (1-H)V+n+1$ ends with the same vertex $v$ of $G$. Assume that there exists $0\le l\le (1-H)V+n$ for which it holds that 
\begin{eqnarray}\label{10}
d(T_{i_{l+1}})\le d(T_{i_l})+H(i_{l+1}-i_l)R.
\end{eqnarray}
Then we consider a concatenation $\overline{T}:=T_{i_l}\overline{T_{i_{l+1}}}$ being a path of the length $k+jR-n-(i_{l+1}-i_l)R$ in $G$. We obtain
$$d(\overline{T})=d(T)+d(T_{i_l})-d(T_{i_{l+1}})>d(k)+HjR-H(i_{l+1}-i_l)R$$
\noindent due to (\ref{9}), (\ref{10}), and we get a contradiction with the choice of the minimal $j$ (see (\ref{9})).

Thus, for every $0\le l\le (1-H)V+n$ we have 
$$d(T_{i_{l+1}})\ge d(T_{i_l})+H(i_{l+1}-i_l)R+1.$$
\noindent Summing up these inequalities for $0\le l\le (1-H)V+n$ we conclude that 
$$d(T_{i_{(1-H)V+n+1}})-d(T_{i_0})\ge H(i_{(1-H)V+n+1}-i_0)R+(1-H)V+n+1$$
\noindent which contradicts to Lemmata~\ref{lower}, \ref{upper}. $\Box$ \vspace{2mm}

Note that $V<(O(Mn))^n$ (see Definition~\ref{vertex}) 
and $R<\exp(V)$.
Lemmata~\ref{lower}, \ref{upper}, \ref{shift}, \ref{stable} entail the following theorem.

\begin{theorem}\label{period}
For $k>(V^2+1)(V+n)+V((1-H)V+n+1)^2$ the tropical Hilbert function $d_{\vec a}(k)$ of the integer vector $\vec a=(a_0,\dots,a_n)$ with an amplitude  $M$ (\ref{0}) fulfils the following equality:
$$d_{\vec a}(k+R)=d_{\vec a}(k)+HR.$$
\noindent for some integer $R<\exp((O(Mn))^n)$ where $H:=H_{\vec a}$ is the tropical entropy of the vector $\vec a$. 
\end{theorem}

We call a function (from the natural numbers to themselves) {\it quasi-linear} if it is a sum of a linear function and a periodic function with an integer period.

\begin{corollary}\label{quasi-linear}
The tropical Hilbert function 
$$d(k)=Hk+r(k)$$
\noindent is quasi-linear for $k>(Mn)^{O(n)}$ where $r(k)$ is a periodic function with an integer period less than $\exp((O(Mn))^n)$.
\end{corollary}




Now we illustrate the constructions in sections~\ref{two} - \ref{five} for three vectors $\vec a \in \ZZ^3$, thus $n=2$ (we use notations  from  sections~\ref{two} - \ref{five}). In each example we construct a graph $G:=G_{\vec a}$ whose vertices $v$ are in a bijective correspondence with polygons $P_v\subset  \RR^2$. Denote by $(y,x)$ coordinates in $\RR^2$.

\begin{example}\label{0,0,0}
First, we consider a vector $\vec a:=(0,0,0)$. Therefore, $M=0$.

The graph $G$ contains three vertices $v$, we list their corresponding polygons $P_v\subset \RR^2$:
$$P_{v_{\infty}}=\{x-y>0\},\ P_{v_0}=\{x-y=0\},\ P_{v_{-\infty}}=\{x-y<0\}.$$
\noindent Note that $B_{v_\infty}=S_{v_\infty}=\{1\},\ B_{v_0}=S_{v_0}=\{1,2\},\ B_{v_{-\infty}}=S_{v_{-\infty}}=\{2\}$ (see Definitions~\ref{vertex},~\ref{after_invariance}).

The edges of $G$ are the following:
$$(v_\infty, v_{-\infty}), (v_0, v_0), (v_0, v_\infty), (v_{-\infty}, v_0).$$
\noindent The edges $(v_\infty, v_{-\infty}), (v_{-\infty}, v_0)$ are constructed according to Definition~\ref{single}, the edge $(v_0, v_0)$ is constructed according to Definition~\ref{bounded}, and the edge $(v_0, v_\infty)$ is constructed according to Definition~\ref{unbounded}. The only augmenting edge is $(v_0, v_\infty)$.

An optimal cycle is $v_{-\infty}, v_0, v_\infty$. It contains a single augmenting edge, hence the entropy $H=1/3$ (see Corollary~\ref{entropy}, also \cite{G20}).

Consider a path
$$T=\underbrace{(v_{-\infty}, v_0, v_\infty) \cdots (v_{-\infty}, v_0, v_\infty)}_p$$
\noindent of the length $3p-1$ in $G$. For any reals $z_0<z_1,\dots,z_{p+1}$ take the point
$(u_1,\dots,u_{3p+1})$ such that $u_{3l-1}=u_{3l}=z_0, 1\le l\le p, u_{3j-2}=z_j, 1\le j\le p+1$. Then $(u_1,\dots,u_{3p+1})$ belongs to the polyhedron $Q_T\subset \RR^{3p+1}$ (see Theorem~\ref{graph}). Vice versa, one can check that any point of $Q_T$ has this form. Therefore, $\dim Q_T = p+2$. One can verify that for every $1\le q\le 3p-1$ the maximum $\max_{T'} \{\dim Q_{T'}\}$, where $T'$ ranges over all paths in $G$ of the length $q$, is attained at the prefix of $T$ of the length $q$. Hence, the tropical Hilbert function $d(k)=\lfloor (k-1)/3\rfloor +2, k\ge 2$ (cf. Theorem~\ref{graph}, Corollary~\ref{quasi-linear}).
\end{example}

\begin{example}\label{1,0,1}
Let a vector $\vec a:= (1,0,1)$, therefore $M=1$.
If to follow the bounds (\ref{30}) formally, then the graph $G:=G_{\vec a}$ should have 13 vertices. We simplify the construction of $G$ imposing stronger bounds than (\ref{30}), namely, $|m(1,2)|\le 1$. One can verify that in case of the chosen $\vec a$ this simplification still computes the entropy $H:=H_{\vec a}$ and the tropical Hilbert function $d:=d_{\vec a}$.

The graph $G$ has 7 vertices $v$, we list their corresponding polygons $P_v\subset \RR^2$: 
$$P_{v_{\infty}}= \{x-y>1\}, P_{v_1}= \{x-y=1\}, P_{v_{0.5}}= \{0<x-y<1\}, P_{v_0}= \{x-y=0\},$$ $$P_{v_{-0.5}}= \{-1<x-y<0\}, P_{v_{-1}}=\{x-y=-1\}, P_{v_{-\infty}}= \{x-y<-1\}$$
\noindent (see Definition~\ref{vertex}). 

Note that the set $B_{v_{\infty}}=\{1\},\ B_{v_{-\infty}}=\{2\}$,  for all other  vertices $v$ of $G$ the set $B_v=\{1,2\}$ (see Definition~\ref{vertex}), the set $S_{v_{\infty}}=\{1\},\ S_{v_1}=\{1,2\}$, for all other vertices $v$ of $G$ the set $S_v=\{2\}$ (see Definition~\ref{after_invariance}). 

The graph $G$ has 12 following edges:
$$e_1:=(v_{\infty}, v_{-\infty}),\ e_2:=(v_1, v_{-1}), e_3:=(v_1, v_{-0.5}),\ e_4:=(v_1, v_0),$$ $$
e_5:=(v_1, v_{0.5}),\ e_6:=(v_1, v_1), e_7:=(v_1, v_{\infty}),\ e_8:=(v_{0.5}, v_{-1}),$$ $$
e_9:=(v_0, v_{-1}),\ e_{10}:=(v_{-0.5}, v_{-1}),\ e_{11}:=(v_{-1}, v_{-1}),\ e_{12}:=´(v_{-\infty}, v_{-1}).$$ 
\noindent Observe that the edges $e_1, e_8, e_9, e_{10}, e_{11}, e_{12}$ (so, the edges not outgoing from the vertex $v_1$) are constructed according to Definition~\ref{single}, the edge $e_7$ is constructed according to Definition~\ref{unbounded}, the edges $e_2, e_3, e_4, e_5, e_6$ are constructed according to Definition~\ref{bounded}. Just the edges $e_3, e_5, e_7$ are augmenting (cf. Proposition~\ref{polyhedron}). 

The graph $G$ has the source $v_1$ (i.e. the vertex without incoming edges), and the sink $v_{-1}$ (i.e. the vertex without outgoing edges). Since no cycle in $G$ contains an augmenting edge, the entropy $H=0$ due to Corollary~\ref{entropy}. This also follows from \cite{G20}.

Every path in $G$ starts with several loops in $v_1$ (perhaps, empty),  ends with several loops in $v_{-1}$ (perhaps, empty), and has at most two edges inbetween. Consider a path
$$T_0=\underbrace{v_1\dots v_1}_p v_{\infty} v_{-\infty} \underbrace{v_{-1}\dots v_{-1}}_q.$$
\noindent Then the polyhedron $Q_{T_0}\subset \RR^{p+q+3}$ (see Theorem~\ref{graph}) is described as follows. For  arbitrary $z_0, z_1>z_0+p+1 \in \RR$ the point
$$(z_0,z_0+1,\dots,z_0+p,z_1,z_0+p,z_0+p-1,\dots, z_0+p-q-1)$$
\noindent belongs to $Q_{T_0}$. Vice versa, these points exhaust $Q_{T_0}$. Therefore, $\dim Q_{T_0}=2$, it holds $n(T_0)=d(T_0)=1$ (cf. Theorem~\ref{graph}). One can verify that any path $T$ in $G$ provides a polyhedron $Q_T$ of dimension at most 2. Indeed, if $T$ starts with the vertex $v_1$ then $n(T)=1, d(T)\le 1$ because any path contains at most one augmenting edge. Otherwise, if $T$ starts with with a vertex different from $v_1$ then $n(T)\le 2, d(T)=0$. Thus, the Hilbert function $d(k)=2, k\ge 2$ (cf. Theorem~\ref{graph}, Corollary~\ref{quasi-linear}).     
\end{example}

\begin{example}\label{0,1,0}
Now consider a vector $\vec a:= (0,1,0)$. Thus, $M=1$. Again as in Example~\ref{1,0,1} we simplify the construction of the graph $G$. It contains the same 7 vertices as in Example~\ref{1,0,1}.

It holds $B_{v_\infty}=\{1\}, B_{v_{-\infty}}=\{2\}$, for all other vertices $v$ of $G$ it holds $B_v=\{1,2\}$. It holds $S_{v_{-1}}=\{1,2\}, S_{v_{-\infty}}=\{2\}$, for all other vertices $v$ of $G$ it holds $S_v=\{1\}$.

The graph $G$ has the following edges:
$$e_1:=(v_\infty, v_{-\infty}),\ e_2:=(v_1, v_{-1}),\ e_3:=(v_{0.5}, v_{-0.5}),\ e_4:=(v_0,v_0),$$ $$e_5:=(v_{-0.5}, v_{0.5}),\ e_6:=(v_{-1}, v_1),\ e_7:= (v_{-1}, v_\infty),\ e_8:=(v_{-\infty}, v_1).$$
\noindent The edge $e_6$ is constructed according to Definition~\ref{bounded}, the edge $e_7$ is constructed according to Definition~\ref{unbounded}, all other edges are constructed according to Definition~\ref{single}. The only augmenting edge is $e_7$.

The unique simple cycle of $G$ which contains the augmenting edge is $(v_{-\infty}, v_1, v_{-1}, v_\infty)$. Therefore, the entropy $H=1/4$ (see also \cite{G20}).

Consider a path
$$T:= \underbrace{(v_{-\infty}, v_1, v_{-1}, v_\infty),\dots, (v_{-\infty}, v_1, v_{-1}, v_\infty)}_p$$
\noindent of the length $4p-1$. Then for any reals $z_0+1<z_1,\dots,z_{p+1}$ take the point $(u_1,\dots,u_{4p+1})$ such that $u_{4l+1}=z_{l+1}, 0\le l\le p,\ u_{4j+2}=u_{4j+4}=u_{4j+3}-1=z_0, 0\le j<p$. The point $(u_1,\dots,u_{4p+1})$ belongs to the polyhedron $Q_T$ (see Theorem~\ref{graph}). Vice versa, one can check that any point of $Q_T$ has the described form, hence $\dim Q_T= p+2$.

One can verify that for every $1\le q\le 4p-1$ the maximum in $\max_{T'} \{\dim Q_{T'}\}$, where $T'$ ranges over all paths of the length $q$ of $G$,  is attained at the prefix of the length $q$ of $T$. Thus, the tropical Hilbert function $d(k)=\lfloor (k-1)/4\rfloor +2$ (cf. Theorem~\ref{graph}, Corollary~\ref{quasi-linear}).  
\end{example}

\begin{remark}\label{regular}
In case when the tropical entropy $H=H_{\vec a}=0$ Lemma~\ref{upper} implies that $d(k)= const$ for sufficiently large $k$, taking into account that $d(k)$ is a non-decreasing function. Recall (see \cite{G20}) that Newton polygon ${\cal N}(\vec a)\subset \RR^2$ for a vector $\vec a=(a_0,\dots,a_n)$ is defined as the convex hull of the rays $\{(i,y)\, : \, y\ge a_i\}$ for $0\le i\le n$. We say that the vector $\vec a$ is {\it regular} \cite{G20} if each point $(i,\, a_i)$ with $a_i<\infty$ is a vertex of ${\cal N}(\vec a)$, and the indices $i$ for which $a_i<\infty$ constitute an arithmetic progression. It was proved in \cite[Corollary 5.7]{G20} that $H_{\vec a}=0$ iff $\vec a$ is regular. For regular $\vec a$ in case when each $(i,\, a_i),\, 0\le i\le n$ is a vertex of ${\cal N}(\vec a)$ one can deduce from \cite[Corollary 4.9]{G12} that $d(k)=k$ for $k\le n$ and $d(k)=n$ for $k\ge n$.
\end{remark}

\section{Tropical boolean vectors}\label{six}
As we already mentioned it would be interesting to extend the results of the paper to arbitrary vectors $\vec a$ involving infinite coordinates. The first step to implementing this idea can be considered as the construction of an appropriate graph $G_{\vec a}$ (cf. section~\ref{two}) for the case when $\vec a$ is a tropical boolean vector (see the Introduction). In this case, the construction looks simpler and contains less technical details comparing to the case considered in the previous sections~\ref{two}, \ref{three}.

\subsection{Construction of a graph for tropical boolean vectors}\label{tropical_construction}
We call a vector $\vec a=(a_0,\dots,a_n)$ \textit{tropical boolean vector} if for all $0\le i\le n$ it holds either $a_i = 0$ or $a_i = \infty$, and $a_0 = a_n = 0.$ 

Below we construct a directed graph $G := G_{\vec a}.$ First we define the vertices of $G$.

\begin{definition}\label{vertices_boolean}
Every  vertex $v$ of $G$ corresponds to an (open in its linear  hull) nonempty polyhedron $P := P_v\subset \mathbb{R}^{n}$ with the condition that for each pair of coordinates $y_r, y_t, 1\le r,t\le n$ a system of equations and strict inequalities defining $P$ contains either $y_r = y_t$ or $y_r < y_t.$ 
\end{definition}

 These linear restrictions set the order on the coordinates $y_1, \dots, y_n.$ The polyhedra $\{P_v\}_v$ constitute a partition of $\mathbb{R}^{n}$. 
 Now we define the edges of $G$.

\begin{definition}\label{edge}
There is an edge $(v,w)$ in $G$ iff there exist vectors $(y_0,\dots,y_{n-1})\in P_v, (y_1,\dots,y_n)\in P_w$ such that the sequence $(y_0,\dots,y_n)\in \RR^{n+1}$ satisfies the vector $\vec a$.
\end{definition}

Similar to subsection~\ref{two.2} for a vertex $v$ of $G$ define $S:=S_v$ as a set of $0\le t\le n-1$ such that $y_t=a_t+y_t=\min_{0\le j\le n-1} \{a_j+y_j\}$. In other words, $t\in S$ iff $a_t=0$ and $y_t\le y_j$ for each $0\le j\le n-1$ such that $a_j=0$. The definition of $S$ does not depend on a choice of a point $(y_0,\dots,y_{n-1})\in P_v$ (cf. Lemma~\ref{invariance}). The following theorem is similar to Theorem~\ref{correctness}.

\begin{theorem}\label{correctness_boolean}
Let $v$ be a vertex of the graph $G:=G_{\vec a}$ (see Definitions~\ref{vertices_boolean}, ~\ref{edge}) and $(y_0,\dots,y_{n-1})\in P_v$. \vspace{2mm}

If a point $(z_0,\dots,z_{n-1})\in P_v$ and a sequence $(z_0,\dots,z_n)\in \RR^{n+1}$ satisfies the vector $\vec a$ then $(z_1,\dots,z_n)\in P_w$ for some edge $(v,w)$ of $G$. \vspace{2mm}

Conversely, let $(y_1,\dots,y_n)\in P_w$ for an edge $(v,w)$ of $G$, and the sequence $(y_0,\dots,y_n)\in \RR^{n+1}$ satisfy the vector $\vec a$. If $t\in S$ for some $0\le t\le n-1$ then $y_n\ge y_t$. \vspace{1mm}

i) Let $t\in S$ for some $1\le t\le n-1$ and $y_r=y_n$ for some $1\le r\le n-1$. Assume that a point  $(z_0,\dots,z_{n-1})\in P_v$. If a point $(z_1,\dots,z_{n-1}, z)\in P_w$ then $z=z_r$. The point $(z_1,\dots,z_{n-1}, z_r)\in P_w$, and the sequence $(z_0,\dots,z_{n-1}, z_r) \in \RR^{n+1}$ satisfies the vector $\vec a$. \vspace{1mm}

ii) Let $t\in S$ for some $1\le t\le n-1$. Assume that $y_{r_1}<y_n$ for some $1\le r_1\le n-1$ and for every $1\le r\le n-1$ neither $y_{r_1}<y_r\le y_n$ nor  $y_n\le y_r$ holds. Then for any point  $(z_0,\dots,z_{n-1})\in P_v$ if a point $(z_1,\dots,z_{n-1}, z)\in P_w$ then $z_{r_1}<z$ and for every $1\le r\le n-1$ neither $z_{r_1}<z_r\le z_n$ nor $z_n\le z_r$ holds. For any $z_{r_1}<z_n\in \RR$ the point $(z_1,\dots,z_n)\in P_w$ and the sequence $(z_0,\dots,z_n) \in \RR^{n+1}$ satisfies the vector $\vec a$. \vspace{1mm}

iii) Let $t\in S$ for some $1\le t\le n-1$. Assume that $y_{r_1}<y_n <y_{r_2}$ for some $1\le r_1, r_2\le n-1$, and for every $1\le r\le n-1$ neither $y_{r_1}<y_r\le y_n$ nor $y_n\le y_r< y_{r_2}$ holds. Then for any point $(z_0,\dots,z_{n-1})\in P_v$ if a point $(z_1,\dots,z_{n-1}, z)\in P_w$ then $z_{r_1}<z<z_{r_2}$ and for every $1\le r\le n-1$ neither $z_{r_1}<z_r\le z_n$ nor $z_n\le z_r< z_{r_2}$ holds. For any $z_n\in \RR, z_{r_1}<z_n<z_{r_2}$ the point $(z_1,\dots,z_n)\in P_w$ and the sequence $(z_0,\dots,z_n) \in \RR^{n+1}$ satisfies the vector $\vec a$. \vspace{1mm}

iv) Let $S=\{0\}$. Then $y_n=y_0$. For any point $(z_0,\dots,z_{n-1})\in P_v$ the point $(z_1,\dots,z_{n-1}, z_0)\in P_w$ and the sequence $(z_0,\dots,z_{n-1}, z_0) \in \RR^{n+1}$ satisfies the vector $\vec a$.
\end{theorem} 

{\bf Proof}. An informal idea of the proof is to transfer inequalities on the differences between the coordinates $y$ to the corresponding inequalities on the coordinates $z$, and back.

Let $(z_0,\dots,z_{n-1})\in P_v$ and a sequence $(z_0,\dots,z_{n-1}, z_n)\in \RR^{n+1}$ satisfy the vector $\vec a$. Assume that $t\in S$ for some $1\le t \le n-1$, then
\begin{eqnarray}\label{40}
z_t=a_t+z_t=\min_{0\le j\le n} \{a_j+z_j\}.
\end{eqnarray}

First, consider the case when $t\in S$ for some $1\le t \le n-1$ and $z_n=z_r$ for some $1\le r\le n$ (cf. i)). Then the sequence $(y_0,\dots,y_{n-1}, y_n=y_r)\in \RR^{n+1}$ also satisfies the vector $\vec a$. Indeed, (\ref{40}) implies that $y_t=a_t+y_t=\min_{0\le j\le n} \{a_j+y_j\}.$ Therefore, due to Definition~\ref{edge} there exists an edge $(v,w)$ of $G$ such that $(y_1,\dots,y_{n-1}, y_n=y_r)\in P_w$. Hence $(z_1,\dots,z_{n-1}, z_n=z_r)\in P_w$ as well. This proves the first statement of the theorem in the case under consideration.

Now consider the case when $t\in S$ for some $1\le t \le n-1$ and $z_n>z_r$ for each $1\le r\le n-1$ (cf. ii)). Then for any $y>\max_{1\le j\le n-1} \{y_j\}$ the sequence $(y_0,\dots,y_{n-1},y)\in \RR^{n+1}$  satisfies the vector $\vec a$. Indeed, (\ref{40}) implies that $y_t=a_t+y_t=\min\{\min_{0\le j\le n-1} \{a_j+y_j\}, y\}$. Due to Definition~\ref{edge} there exists an edge $(v,w)$ (independent of $y$) of $G$ such that $(y_1,\dots,y_{n-1},y)\in P_w$. Hence  $(z_1,\dots,z_{n-1},z_n)\in P_w$ as well. This proves the first statement of the theorem in the case under consideration. 

The case when $t\in S$ for some $1\le t \le n-1$ and $z_{r_1}<z_n<z_{r_2}$ for some $1\le r_1,r_2\le n-1$ such that for each $1\le r\le n-1$ neither $z_{r_1}<z_r\le z_n$ nor $z_n\le z_r<z_{r_2}$ holds (cf. iii)) can be studied in a similar manner as the previous case.

Finally, consider the case $S=\{0\}$ (cf. iv)). Then $z_0<z_l$ for each $1\le l\le n-1$ for which $a_l=0$. Therefore, $z_n=z_0$. Hence the sequence $(y_0,\dots,y_{n-1},y_0)\in \RR^{n+1}$  satisfies the vector $\vec a$. Due to Definition~\ref{edge} there exists an edge $(v,w)$ of $G$ such that $(y_1,\dots,y_{n-1},y_0)\in P_w$. Therefore $(z_1,\dots,z_{n-1},z_0)\in P_w$ as well. This proves the first statement of the theorem.

One can directly verify the second statement of the theorem. $\Box$ 

\begin{corollary}
The edges of the graph $G$ (see Definition~\ref{edge}) do not depend on choices of points $(y_0,\dots,y_{n-1})\in P_v$. 
\end{corollary}

\begin{remark}\label{after_correctness}
Let an edge $(v,w)$ fulfill the assumptions of one of the items Theorem~\ref{correctness_boolean}~i), ii), iii) and a point $(z_0,\dots,z_{n-1})\in P_v$.  Then for any $z\in \RR$ such that $(z_1,\dots,z_{n-1},z)\in P_w$ the sequence  $(z_0,\dots,z_{n-1},z)\in \RR^{n+1}$ satisfies the vector $\vec a$. In contrast, in case of Theorem~\ref{correctness_boolean}~iv) only for the value $z=z_0$ it holds that the sequence $(z_0,\dots,z_{n-1},z)$ satisfies the vector $\vec a$ (cf. Theorem~\ref{correctness}).
\end{remark}

\subsection{The polyhedron of tropical recurrent sequences produced along a path of the graph}\label{trop Q}
Consider an arbitrary path $T$ of a length $k$ with vertices $v_0,\dots,v_k$ in the graph $G_{\vec a}.$ Similar to subsection~\ref{two.2} we describe a recursive process producing along $T$ tropical recurrent sequences satisfying the vector $\vec a$. For the first vertex $v_0$ take any vector $(y_1,\dots,y_n)\in P_{v_0}$. Assume by recursion that a tropical recurrent sequence $(y_1,\dots,y_{k+n})$ is already produced along $T$. Then $(y_{k+1},\dots,y_{k+n})\in P_{v_k}$. Take an edge $(v_k,w)$ of $G$ and denote by $T_w$ the extension of $T$ by  $(v_k,w)$. We choose $y_{k+n+1}\in \RR$ such that $(y_{k+2},\dots,y_{k+n+1})\in P_w$ and the sequence $(y_{k+1},\dots,y_{k+n+1})\in \RR^{n+1}$ satisfies the vector $\vec a$. Thus, the tropical recurrent sequence $(y_1,\dots,y_{k+n+1})$ is produced along $T_w$. Theorem~\ref{correctness_boolean} justifies that a required $y_{k+n+1}$ exists and moreover, Theorem~\ref{correctness_boolean} describes all possible $y_{k+n+1}$. This completes the description of the recursive process.

Denote by $Q_T\subset \RR^{k+n}$ the set of all the tropical recurrent sequences  produced along $T$  by the described recursive process. One can define $Q_T$ by imposing linear inequalities for each edge of $T$. Say, for an edge $(v_i,v_{i+1}), 0\le i\le k-1$ we impose that the point $(y_{i+1},\dots,y_{i+n+1})$ belongs to $P_{v_i}$, the point $(y_{i+2},\dots,y_{i+n+2})$ belongs to $P_{v_{i+1}}$. This suffices for edges $(v_i,v_{i+1})$ fulfilling the items Theorem~\ref{correctness_boolean}~i), ii), iii). In case of Theorem~\ref{correctness_boolean}~iv) one has to impose an extra condition that $y_{i+1}=y_{i+n+2}$, i.e. the sequence $(y_{i+1},\dots,y_{i+n+2})\in \RR^{n+1}$ satisfies the vector $\vec a$. Thus, $Q_T$ is (an open in its linear hull) polyhedron.

If an edge  $(v_i,v_{i+1})$ fulfills one of the items  Theorem~\ref{correctness_boolean}~i), iv) we call the edge {\it rigid}, otherwise, if the edge fulfills one of the items  Theorem~\ref{correctness_boolean}~ii), iii) we call the edge {\it augmenting}. Similar to subsection~\ref{two.2} when the  edge $(v_k,w)$ is rigid the value of $y_{k+n+1}$ is unique, while when the edge is augmenting  the values of $y_{k+n+1}$ vary in an open interval. Therefore, when the edge $(v_k,w)$ is rigid the polyhedron $Q_{T_w}$ is homeomophic to $Q_T$, while when the edge is augmenting the polyhedron $Q_{T_w}$ is homeomophic to $Q_T\times \RR$.

Conversely,  Theorem~\ref{correctness_boolean} implies that any tropical recurrent sequence satisfying the vector $\vec a$ emerges along a suitable path of $G$ in the described above recursive process. Thus, the tropical prevariety of  all tropical recurrent sequences of a length $k+n$ satisfying the vector $\vec a$ coincides with the union of polyhedra $Q_T$ over all the paths of the length $k$ in $G$.

For a path $T$ in the graph $G$ denote by $d(T)$ the number of augmenting edges 
in $T$. By $n(T)\le n$ denote the number of the pairwise distinct coordinates in
$(y_1,\dots,y_n)\in P_{v_0}$ for the first vertex $v_0$ of $T$. 
We summarize the proved above in the following theorem which is analogous to the Theorem~\ref{graph} for the case when $\vec a$ is a tropical boolean vector.

\begin{theorem}\label{boolean_graph}
For any tropical boolean vector $\vec a:=(a_0,\dots,a_n)$
(i.~e. $a_0=a_n=0$ and each $a_i,\, 0\le i\le n$ equals either $0$ or $\infty$)
a finite directed graph $G:=G_{\vec a}$ is constructed  with the following properties. 
For an arbitrary  path $T$ 
of a length $k$ 
in $G$ 
denote by $Q_T\subset \RR^{k+n}$ the polyhedron of all the tropical recurrent sequences satisfying the vector $\vec a$ and corresponding (as described above in this subsection) to the path $T$ in $G$. Then $\dim (Q_T)=d(T)+n(T)$.
Moreover, the union of polyhedra $Q_T$ over all the paths $T$ of the length $k$ 
coincides with the tropical prevariety of all the tropical recurrent sequences of the length $k+n$ satisfying the vector $\vec a$.
\end{theorem}

Now let us notice that all the arguments presented in sections~\ref{four} and \ref{five} for the graph constructed in section~\ref{two} are also true in the case of tropical boolean vectors. Indeed, both definitions of $n(T)$ and $d(T)$ and thereby, Hilbert function $d_{\vec a}(s)$ coincide, respectively, with the definitions for the case when $\vec a$ has a finite amplitude. Moreover,  an analogue of Theorem~\ref{graph} holds in the tropical boolean case (Theorem~\ref{boolean_graph}). As all the statements from sections~\ref{four} and \ref{five} (except of Theorem~\ref{period} and Corollary~\ref{quasi-linear}) depend only on $d_{\vec a}(s)$ and on Theorem~\ref{graph}, we can formulate the following corollaries.

\begin{corollary}\label{boolean cor}
Lemmata~\ref{lower}, \ref{upper}, \ref{optimal}, \ref{shift}, \ref{stable} and Corollary~\ref{entropy} hold when $\vec a$ is a tropical boolean vector.
\end{corollary}
{\bf Proof}. Follows from the proofs of the mentioned statements. $\Box$
\begin{corollary}\label{boolean theorems}
Theorem~\ref{period} and Corollary~\ref{quasi-linear} hold when $\vec a$ is a tropical boolean vector putting in the bounds $M = 1.$
\end{corollary}
{\bf Proof}. From subsection~\ref{tropical_construction} it follows that the number of vertices $V$ in $G$ is less than the amount of orders on an $n$-element set, hence it is less than $n^n.$ Thus, we can put $M = 1$ in the bounds. The remaining part of the proof is literally as in the proofs of the mentioned statements. $\Box$

$\quad$

Let us note that for any tropical boolean vector $\vec a$ every vertex of $G_{\vec a}$ could be presented as a sequence of  numbers from $0$ to $n$ which reflects the order between the coordinates of the corresponding polyhedron. For example polyhedron $P_v\subset \RR^3$ defined by a system of equalities and inequalities $\{y_1 < y_2,\:\: y_1 < y_3,\:\: y_2 = y_3\}$ could be presented as sequence $\{0, 1, 1\}.$
\begin{example}\label{(0, 0, infty, 0)}
Now we illustrate constructions in this section for vector $\vec a = (0, 0, \infty, 0).$

It is not hard to see that there are 13 different orderings on 3 coordinates. We list corresponding polyhedra (see Definition~\ref{vertices_boolean}) $P_v\subset \RR^3:$
\begin{center}
$P_{v_{0, 0, 0}} = \{y_1 = y_2 = y_3\},\:\: P_{v_{1, 0, 0}} = \{y_1 > y_2,\:\: y_2 = y_3\}$,

$P_{v_{0, 0, 1}} = \{y_1 = y_2,\:\: y_2 < y_3\},\:\:P_{v_{0, 1, 0}} = \{y_1 < y_2,\:\: y_1 = y_3\}$,

$P_{v_{1, 0, 1}} = \{y_1 > y_2,\:\: y_1 = y_2\},\:\:P_{v_{1, 1, 0}} = \{y_1 = y_2,\:\: y_2 > y_3\}$,

$P_{v_{0, 1, 1}} = \{y_1 < y_2,\:\: y_2 = y_3\},\:\:P_{v_{0, 2, 1}} = \{y_1 < y_2,\:\: y_1 < y_3\},
$

$P_{v_{1, 2, 0}} = \{y_1 > y_3,\:\: y_2 > y_1\},\:\:P_{v_{2, 1, 0}} = \{y_1 > y_2,\:\: y_2 > y_3\},
$

$P_{v_{2, 0, 1}} = \{y_1 > y_3,\:\: y_2 < y_3\},\:\:P_{v_{0, 1, 2}} = \{y_1 < y_2,\:\: y_2 < y_3\},
$

$P_{v_{1, 0, 2}} = \{y_1 > y_2,\:\: y_1 < y_3\}.
$
\end{center}

The graph $G_{\vec a}$ has the following 13 edges (see Definition~\ref{edge}):
\begin{center}
    $e_1 := (v_{0, 0, 0}, v_{0, 0, 1}),\:\: e_2 = (v_{1, 0, 0}, v_{0, 0, 0}), \:\: e_3 = (v_{0, 0, 1}, v_{0, 1, 0}),$
    
    $e_4 := (v_{0, 0, 1}, v_{0, 1, 1}),\:\: e_5 = (v_{0, 0, 1}, v_{0, 2, 1}), \:\: e_6 = (v_{0, 0, 1}, v_{0, 1, 2}),$
    
    $e_7 := (v_{0, 1, 0}, v_{1, 0, 0}),\:\: e_8 = (v_{1, 0, 1}, v_{0, 1, 0}), \:\: e_9 = (v_{0, 1, 1}, v_{1, 1, 0}),$

    $e_{10} := (v_{0, 2, 1}, v_{2, 1, 0}),\:\: e_{11} = (v_{2, 0, 1}, v_{0, 1, 0}), \:\: e_{12} = (v_{0, 1, 2}, v_{1, 2, 0}),$

    $e_{13} := (v_{1, 0, 2}, v_{0, 1, 0}).$
\end{center}
Moreover, just the edges $e_1, e_5$ and $e_6$ are augmenting (see Theorem~\ref{correctness_boolean}). 

There is only one cycle in this graph: 
$$(v_{0, 0, 0}, v_{0, 0, 1}, v_{0, 1, 0}, v_{1, 0, 0}).$$

As only one of these edges is augmenting we obtain that $H_{\vec a} = \frac{1}{4}$ (cf. Corollary~\ref{boolean cor}). 

Now we give the detailed description of the tropical Hilbert function $d(k)$ in this case (cf. Corollary~\ref{boolean theorems}):
\begin{itemize}
    \item if $k = 1, 2, 3$ then $d(k) = 3$. The maximum of the dimension $\dim (Q_T)$ over paths of a given length in $G_{\vec a}$ (see Theorem~\ref{boolean_graph}) is attained at the path
    $$T := e_{11}e_7e_2e_1;$$
    \item if $k = 4p,$ where $p\ge 1$ then $d(k) = 3 + p.$ The maximum is attained at the path
    $$T := e_{11}\underbrace{(v_{0, 0, 0}, v_{0, 0, 1}, v_{0, 1, 0}, v_{1, 0, 0}) \cdots (v_{0, 0, 0}, v_{0, 0, 1}, v_{0, 1, 0}, v_{1, 0, 0})}_{p - 1}e_7e_2e_1;$$
    \item if $k = 4p + 1,$ where $p\ge 1$ then $d(k) = 3 + (p + 1).$ The maximum is attained at the path
    $$T := e_{11}\underbrace{(v_{0, 0, 0}, v_{0, 0, 1}, v_{0, 1, 0}, v_{1, 0, 0}) \cdots (v_{0, 0, 0}, v_{0, 0, 1}, v_{0, 1, 0}, v_{1, 0, 0})}_{p - 1}e_7e_2e_1e_5;$$
    \item if $k = 4p + 2,$ where $p\ge 1$ then $d(k) = 3 + (p + 1).$ The maximum is attained at the path
    $$T := e_{11}\underbrace{(v_{0, 0, 0}, v_{0, 0, 1}, v_{0, 1, 0}, v_{1, 0, 0}) \cdots (v_{0, 0, 0}, v_{0, 0, 1}, v_{0, 1, 0}, v_{1, 0, 0})}_{p - 1}e_7e_2e_1e_6e_{12};$$
    \item if $k = 4p + 3,$ where $p\ge 1$ then $d(k) = 3 + p.$ The maximum is attained at the path
    $$T := e_{11}\underbrace{(v_{0, 0, 0}, v_{0, 0, 1}, v_{0, 1, 0}, v_{1, 0, 0}) \cdots (v_{0, 0, 0}, v_{0, 0, 1}, v_{0, 1, 0}, v_{1, 0, 0})}_{p}e_7e_2.$$
\end{itemize}
\end{example}

\section{Sharp bounds on the  tropical entropy}\label{seven}

\subsection{Sharp lower bound on the positive entropy}

In this section our main goal is to  prove that if for a vector $\vec a=(a_0,\dots,a_n)\in \ZZ^{n+1}$ its tropical entropy $H(\vec a)>0$  then $H(\vec a)\ge \frac{1}{4}.$ Together with the example \cite[Example 5.5]{G20} demonstrating that $H(0,\, 1,\, 0)=1/4$ (cf. also Example~\ref{0,1,0}) we will conclude that this bound is sharp. This result is the answer to the hypothesis that was formulated in  \cite[Remark 5.6]{G20} (for the criterion of positivity of the tropical entropy see  \cite[Corollary 5.7]{G20}, cf. also Remark~\ref{regular}). 

For convenience, we use the following assumptions in this section:

\begin{itemize}
	\item we consider tropical sequences $\{z_{\mathcal{I}}\}_{\mathcal{I}\in \mathbb{Z}}$ infinite in both directions. To obtain finite tropical sequence it is enough to consider only $\mathcal{I} \in \mathbb{N};$
	\item for most of the considered cases we attach diagrams with simple examples. To save place on these diagrams the $\ge$ sign is replaced by sign $\thicksim$ above the symbol.

\end{itemize}
\begin{theorem}\label{sharp estimate}
If a vector $\vec a$ is not regular then $H(\vec a)\ge \frac{1}{4}.$
\end{theorem}
{\bf Proof}. 
Consider Newton polygon ${\cal N}(\vec a)$ of the vector $\vec a$ (see Remark~\ref{regular}). It has several bounded edges and two unbounded edges. First, assume that there is a bounded edge of ${\cal N}(\vec a)$ such that there are at least three points of $\vec a$, i.~e. of the form $(i, a_i)$ (in this case we follow the proof of \cite[Theorem 5.5]{G20}). Making a suitable affine transformation of the plane one can suppose w.l.o.g. that this edge lies on the abscissas axis and $(0, 0)$ is its left end-point (the transformation, perhaps, converts the tropical polynomial $f$, see (\ref{-3}), into a tropical Laurent polynomial, the proof still goes through for the latter). Consider the points of $\vec a$ located on this edge: $E_0:=\{(e, 0) : a_e = 0\}$, then $|E_0|\ge 3$ by our assumption.  One can assume w.l.o.g. that the greatest common divisor $GCD(E_0)$ of the differences $e_1 - e_2$ of all the pairs of the elements $e_1,\:e_2\in E_0$ of $E_0$ equals $1$. Otherwise, one can consider separately all $GCD(E_0)$ arithmetic progressions with the difference $GCD(E_0)$.

Pick any three elements of $E_0$ not all with the same parity, say $0, 2v, u$ w.l.o.g. where $v \ge 1$ and $u$ being odd. Consider the following tropical recurrent sequence $\{z_{\mathcal{I}}\}_{\mathcal{I}\in\mathbb{Z}}$ satisfying $a:$
 \begin{itemize}
        \item $z_{2l + 1} = 0,$ for $0\le l\in\mathbb{Z};$
        \item $z_{2(2qv + r)} = 0,$ for $q\in \mathbb{Z}$ and $0\le r < v;$
        \item $z_{2((2q+1)v + r)} \ge 0,$ for $q\in\mathbb{Z}$ and $0\le r < v.$
\end{itemize}

\begin{center}
\begin{tikzpicture}
    \node[orange] at (-1.2, 1.1){$\{z_{\mathcal{I}}\}_{\mathcal{I}\in\mathbb{Z}}$};
    \node[orange] at (-1.2, 0.6){sequence};
    \node[orange] at (-1.2, 0.1){values:};
    \node[blue] at (-1.2, -0.5){Index:};

	\foreach \i in {0, ...,  16}
		\fill [orange] (\i * 0.7, 0) circle (3pt);
	\node[blue] at (0 * 0.7, -0.5){0};
	\node[blue] at (4 * 0.7, -0.5){$2v$};
	\node[blue] at (8 * 0.7, -0.5){$4v$};
	\node[blue] at (12 * 0.7, -0.5){$6v$};
 	\node[blue] at (16 * 0.7, -0.5){$8v$};
	\node[blue] at (3 * 0.7, -0.5){$u$};

	\node[red] at (0 * 0.7, 0.4){0};
	\node[orange] at (1 * 0.7, 0.4){0};
	\node[orange] at (2 * 0.7, 0.4){0};
	\node[orange] at (3 * 0.7, 0.4){0};
	\node[red] at (4 * 0.7, 0.5){$\widetilde{0}$};
	\node[orange] at (5 * 0.7, 0.4){0};
	\node[orange] at (6 * 0.7, 0.5){$\widetilde{0}$};
	\node[orange] at (7 * 0.7, 0.4){0};
	\node[red] at (8 * 0.7, 0.4){0};
	\node[orange] at (9 * 0.7, 0.4){0};
	\node[orange] at (10 * 0.7, 0.4){0};
	\node[orange] at (11 * 0.7, 0.4){0};
	\node[red] at (12 * 0.7, 0.5){$\widetilde{0}$};
	\node[orange] at (13 * 0.7, 0.4){0};
	\node[orange] at (14 * 0.7, 0.5){$\widetilde{0}$};
	\node[orange] at (15 * 0.7, 0.4){0};
 	\node[red] at (16 * 0.7, 0.4){0};

    \draw [<->][red] (0.7 * 12 + 0.06, 0.85) -- (0.7 * 16 - 0.06, 0.85);
    \draw [<->][red] (0.7 * 8 + 0.06, 0.85) -- (0.7 * 12 - 0.06, 0.85);
    \draw [<->][red] (0.7 * 4 + 0.06, 0.85) -- (0.7 * 8 - 0.06, 0.85);
    \draw [<->][red] (0.7 * 0 + 0.06, 0.85) -- (0.7 * 4 - 0.06, 0.85);

    \node[red] at (0.7 * 2, 1.1){$2v$};
    \node[red] at (0.7 * 6, 1.1){$2v$};
    \node[red] at (0.7 * 10, 1.1){$2v$};
    \node[red] at (0.7 * 14, 1.1){$2v$};

\end{tikzpicture}
\end{center}

Taking finite fragments $(z_1, \dots, z_N)$ with growing $N$ we conclude that $H(\vec a)\ge \frac{1}{2 + 1 + 1} = \frac{1}{4}.$ \vspace{2mm}

In the other case we have $0,$ $2v + 1$ and $2u + 1$ $\in E_0$ and thus we have the following sequence $\{z_{\mathcal{I}}\}_{\mathcal{I}\in\mathbb{Z}}:$

\begin{itemize}
        \item $z_{2l} = 0,$ for $0\le l\in\mathbb{Z};$
	\item $z_{2(2q(u - v) + r) + 1} = 0,$ for $q\in \mathbb{Z}$ and $0\le r < u - v;$
	\item $z_{2((2q+1)(u - v) + r) + 1} \ge 0,$ for $q\in\mathbb{Z}$ and $0\le r < u - v.$
\end{itemize}

\begin{center}
\begin{tikzpicture}
    \node[orange] at (-1.2, 1.1){$\{z_{\mathcal{I}}\}_{\mathcal{I}\in\mathbb{Z}}$};
    \node[orange] at (-1.2, 0.6){sequence};
    \node[orange] at (-1.2, 0.1){values:};
    \node[blue] at (-1.2, -0.5){Index:};

	\foreach \i in {0, ...,  9}
		\fill [orange] (\i * 1.0, 0) circle (3pt);
	\node[blue] at (0 , -0.5){0};
	\node[blue] at (1.0 * 3, -0.5){$2v + 1$};
	\node[blue] at (1.0 * 5, -0.5){$2u + 1$};
 
	\node[orange] at (0, 0.5){0};
	\node[red] at (1.0, 0.5){0};
	\node[orange] at (1.0 * 2, 0.5){0};
	\node[orange] at (1.0 * 3, 0.5){0};
	\node[orange] at (1.0 * 4, 0.5){0};
	\node[red] at (1.0 * 5, 0.5){$\widetilde{0}$};
	\node[orange] at (1.0 * 6, 0.5){0};
	\node[orange] at (1.0 * 7, 0.5){$\widetilde{0}$};
	\node[orange] at (1.0 * 8, 0.5){0};
	\node[red] at (1.0 * 9, 0.5){0};

    \draw [<->][red] (1.0 + 0.06, 0.85) -- (1.0 * 5 - 0.06, 0.85);
    \draw [<->][red] (1.0 * 5 + 0.06, 0.85) -- (1.0 * 9 - 0.06, 0.85);
    \node[red] at (1.0 * 3, 1.1){$2|u - v|$};
    \node[red] at (1.0 * 7, 1.1){$2|u - v|$};

\end{tikzpicture}
\end{center}

Now we assume that no edge of ${\cal N}(\vec a)$ contains a point of $\vec a$ other than two vertices of this edge. We take an edge of ${\cal N}(\vec a)$ with the biggest difference of indices of its vertices. Due to a suitable affine transformation we suppose w.l.o.g. that these vertices are $(0, 0)$ and $(n_0, 0).$ There exists $i\in J$ such that $n_0$ does not divide $i,$ since $\vec a$ is not regular. Among such $i$ we pick $i_0$ for which $c := a_{i_0}$ is minimal. Then $c > 0.$ Denote $k = GCD(n_0, i_0).$ When $\frac{n_0}{k}$ is even we consider the sequence $\{z_{\mathcal{I}}\}_{\mathcal{I}\in\mathbb{Z}}$:
\begin{itemize}
    \item $z_{qn_0 - 2ji_0 + i} = 0, $ when $0\le 2j \le \frac{n_0}{k};$
    \item $z_{2qn_0 - (2j + 1)i_0 + i} = c,$ when $0 < 2j + 1 < \frac{n_0}{k};$
    \item $z_{(2q + 1)n_0 - (2j + 1)i_0 + i}\ge c,$ when $0 < 2j + 1 < \frac{n_0}{k},$
\end{itemize}
for $q\in\ZZ, 0\le i < k.$

\begin{center}
\begin{tikzpicture}
    \node[orange] at (-1.2, 1.1){$\{z_{\mathcal{I}}\}_{\mathcal{I}\in\mathbb{Z}}$};
    \node[orange] at (-1.2, 0.6){Sequence};
    \node[orange] at (-1.2, 0.1){values:};
    \node[blue] at (-1.2, -0.5){Index:};

	\foreach \i in {0, ...,  25}
		\fill [orange] (\i * 0.45, 0) circle (3pt);
	\node[blue] at (0 * 0.45 , -0.5){0};
	\node[blue] at (6 * 0.45, -0.5){$i_0$};
	\node[blue] at (8 * 0.45, -0.5){$n_0$};
	\node[blue] at (16 * 0.45, -0.5){$2n_0$};
	\node[blue] at (24 * 0.45, -0.5){$3n_0$};

	\node[orange] at (0 * 0.45, 0.5){0};
	\node[orange] at (1 * 0.45, 0.5){0};
	\node[orange] at (2 * 0.45, 0.5){$\widetilde c$};
	\node[orange] at (3 * 0.45, 0.5){$\widetilde c$};
	\node[orange] at (4 * 0.45, 0.5){0};
	\node[orange] at (5 * 0.45, 0.5){0};
	\node[red] at (6 * 0.45, 0.5){$\widetilde c$};
	\node[orange] at (7 * 0.45, 0.5){$\widetilde c$};
	\node[orange] at (8 * 0.45, 0.5){0};
	\node[orange] at (9 * 0.45, 0.5){0};
	\node[orange] at (10 * 0.45, 0.5){c};
	\node[orange] at (11 * 0.45, 0.5){c};
	\node[red] at (12 * 0.45, 0.5){0};
	\node[orange] at (13 * 0.45, 0.5){0};
	\node[orange] at (14 * 0.45, 0.5){c};
	\node[orange] at (15 * 0.45, 0.5){c};
	\node[orange] at (16 * 0.45, 0.5){0};
	\node[orange] at (17 * 0.45, 0.5){0};
	\node[red] at (18 * 0.45, 0.5){$\widetilde c$};
	\node[orange] at (19 * 0.45, 0.5){$\widetilde c$};
	\node[orange] at (20 * 0.45, 0.5){0};
	\node[orange] at (21 * 0.45, 0.5){0};
	\node[orange] at (22 * 0.45, 0.5){$\widetilde c$};
	\node[orange] at (23 * 0.45, 0.5){$\widetilde c$};
	\node[red] at (24 * 0.45, 0.5){0};
	\node[orange] at (25 * 0.45, 0.5){0};
 
    \draw [<->][red] (18 * 0.45 + 0.06, 0.75) -- (24 * 0.45 - 0.06, 0.75);
    \draw [<->][red] (12 * 0.45 + 0.06, 0.75) -- (18 * 0.45 - 0.06, 0.75);
    \draw [<->][red] (6 * 0.45 + 0.06, 0.75) -- (12 * 0.45 - 0.06, 0.75);
    \node[red] at (21 * 0.45, 1.0){$i_0$};
    \node[red] at (15 * 0.45, 1.0){$i_0$};
    \node[red] at (9 * 0.45, 1.0){$i_0$};
    
    \draw [decorate, decoration={brace, amplitude = 3}][violet] (0, 0.75) -- (0.45, 0.75);
    \draw [decorate, decoration={brace, amplitude = 3}][violet] (2 * 0.45, 0.75) -- (3 * 0.45, 0.75);
    \draw [decorate, decoration={brace, amplitude = 3}][violet] (4 * 0.45, 0.75) -- (5 * 0.45, 0.75);
    \node[violet] at (0.45 / 2, 1){$k$};
    \node[violet] at (0.45 * 2.5, 1){$k$};
    \node[violet] at (0.45 * 4.5, 1){$k$};
\end{tikzpicture}
\end{center}

This sequence satisfies $\vec a$ and taking finite fragments $(z_1, \dots, z_N)$ with growing $N$ we conclude that $H(\vec a)\ge \frac{1}{2 + 1 + 1} = \frac{1}{4}.$ Thus further we suppose that $\frac{n_0}{k}$ is odd.

We denote the first (respectively, the last) index of $a$ by $\mathcal{B}$ (respectively, by $\mathcal{E}$). Thus, the projection of ${\cal N}(\vec a)$ is the interval from $\mathcal{B}$ to $\mathcal{E}$ on the abscissas axis. 
Before we prove the statement of the theorem in general case let us prove the following lemma.

\begin{lemma}\label{twomins}
If there exists $i_1\neq i_0$ such that $n_0\nmid i_1$ and $a_{i_1}=a_{i_0}$ 
then $H(\vec a)\ge\frac{1}{4}.$
\end{lemma}
{\bf Proof}.

	\large\textbf{L.I}
	\normalsize
	$\:\:$Let $n_0|(i_1 - i_0).$
	
	Then we consider a sequence $\{z_\mathcal{I}\}_{\mathcal{I}\in\ZZ}$ such that:
    \begin{itemize}
        \item $z_{qn_0 - 2ji_0 + i} = 0$ when $0\le 2j < \frac{n_0}{k};$
        \item $z_{qn_0 - (2j + 1)i_0 + i}\ge c$ when $0 < 2j + 1 < \frac{n_0}{k}$
    \end{itemize}
    for $q\in \ZZ$, $0\le i < k.$ 

\begin{center}
\begin{tikzpicture}
    \node[orange] at (-1.2, 0.6){Sequence};
    \node[orange] at (-1.2, 1.1){$\{z_{\mathcal{I}}\}_{\mathcal{I}\in\mathbb{Z}}$};
    \node[orange] at (-1.2, 0.1){values:};
    \node[blue] at (-1.2, -0.5){Index:};

	\foreach \i in {0, ...,  15}
		\fill [orange] (\i *0.6, 0) circle (3pt);
	\node[blue] at (0 * 0.6 , -0.5){0};
	\node[blue] at (2 * 0.6, -0.5){$i_0$};
	\node[blue] at (5 * 0.6, -0.5){$n_0$};
	\node[blue] at (7 * 0.6, -0.5){$i_1$};
	\node[blue] at (10 * 0.6, -0.5){$2n_0$};
	\node[blue] at (15 * 0.6, -0.5){$3n_0$};

	\node[orange] at (0 * 0.6, 0.5){0};
	\node[orange] at (1 * 0.6, 0.5){0};
	\node[red] at (2 * 0.6, 0.5){0};
	\node[orange] at (3 * 0.6, 0.5){$\widetilde c$};
	\node[red] at (4 * 0.6, 0.5){$\widetilde c$};
	\node[orange] at (5 * 0.6, 0.5){0};
	\node[red] at (6 * 0.6, 0.5){0};
	\node[orange] at (7 * 0.6, 0.5){0};
	\node[red] at (8 * 0.6, 0.5){$\widetilde c$};
	\node[orange] at (9 * 0.6, 0.5){$\widetilde c$};
	\node[red] at (10 * 0.6, 0.5){0};
	\node[orange] at (11 * 0.6, 0.5){0};
	\node[orange] at (12 * 0.6, 0.5){0};
	\node[orange] at (13 * 0.6, 0.5){$\widetilde c$};
	\node[orange] at (14 * 0.6, 0.5){$\widetilde c$};
	\node[orange] at (15 * 0.6, 0.5){0};

    \draw [<->][red] (0.6 * 2 + 0.06, 0.85) -- (0.6 * 4 - 0.06, 0.85);
    \draw [<->][red] (0.6 * 4 + 0.06, 0.85) -- (0.6 * 6 - 0.06, 0.85);
    \draw [<->][red] (0.6 * 6 + 0.06, 0.85) -- (0.6 * 8 - 0.06, 0.85);
    \draw [<->][red] (0.6 * 8 + 0.06, 0.85) -- (0.6 * 10 - 0.06, 0.85);
    \node[red] at (0.6 * 9, 1.1){$i_0$};
    \node[red] at (0.6 * 7, 1.1){$i_0$};
    \node[red] at (0.6 * 5, 1.1){$i_0$};
    \node[red] at (0.6 * 3, 1.1){$i_0$};
\end{tikzpicture}
\end{center}

This sequence satisfies $a.$
    
    Indeed,
    \begin{itemize}
        \item For $m = qn_0 - 2ji_0 + i$ we have $\min_{\mathcal{B}\le v\le \mathcal{E}}\{a_v + z_{v + m}\} = 0,$ the minimum is attained at indices $m$ and $m + n_0.$ 
        \item For $m = qn_0 - (2j + 1)i_0 + i$ we have $\min_{\mathcal{B}\le v\le \mathcal{E}}\{a_v + z_{v + m}\} = c,$ the minimum is attained at indices $m + i_0$ and $m + i_1.$
    \end{itemize} Taking finite fragments $(z_1, \dots, z_N)$ with growing $N$ we conclude that $H(\vec a)\ge \frac{1}{2}$ for even $\frac{n_0}{k}$ and $H(\vec a)\ge \frac{\frac{\frac{n_0}{k} - 1}{2}\cdot k}{n_0}\ge\frac{1}{3}.$
    
	$\quad$

     \large\textbf{L.II}
     \normalsize
     $\:\:$Let $n_0 \nmid (i_1 - i_0).$

     $\quad$

	\large\textbf{L.II.1}
	\normalsize
	$\:\:$Assume that $\mathbf{k = 1}$ (since we consider the case where $\frac{n_0}{k}$ is odd, thus $n_0$ is odd).
        
		First, consider sequence $\{w_\mathcal{I}\}_{\mathcal{I}\in\ZZ}$ such that:
        \begin{itemize}
            \item $w_{qn_0 - 2ji_0} = 0$ when $0\le 2j\le n_0;$
            \item $w_{2qn_0 - (2j + 1)i_0} = c$ when $0 < 2j + 1 < n_0;$
            \item $w_{(2q+1)n_0 - (2j + 1)i_0} \ge c$ when $0 < 2j + 1 < n_0$
        \end{itemize}
    for $q \in \ZZ.$

This sequence satisfies $a.$

    Indeed, 
    \begin{itemize}
    \item For $m = qn_0 - 2ji_0$ we have $\min_{\mathcal{B}\le v\le \mathcal{E}}\{a_v + w_{v + m}\} = 0,$ the minimum is attained at indices $m$ and $m + n_0.$ 
    \item For $m = 2qn_0 - (2j + 1)i_0$ we have $\min_{\mathcal{B}\le v\le \mathcal{E}}\{a_v + w_{v + m}\} = c,$ the minimum is attained at indices $m + i_0$ and $m.$ 
    \item For $m = (2q+1)n_0 - (2j + 1)i_0$ we have $\min_{\mathcal{B}\le v\le \mathcal{E}}\{a_v + w_{v + m}\} = c,$ the minimum is attained at indices $m + i_0$ and $m + n_0.$
    \end{itemize}
    
    $\quad$
    
    Now, we claim that there exists $0 < 2l + 1 < n_0$ such that $w_{qn_0 - (2l + 1)i_0 + i_1} = 0$ for all $q\in\ZZ.$ Note, that if we found $2l' + 1$ such that $w_{qn_0 - (2l' + 1)i_0 + i_1'} = 0$ for all $q\in\ZZ$ for some $i_1'\equiv i_1,$ then $w_{qn_0 - (2l' + 1)i_0 + i_1} = 0$ for all $q\in\ZZ.$ 
    
    Recall that $GCD(i_0, n_0) = 1$ and $n_0 \nmid i_1$, therefore there exists $m_{i_1}$ such that $m_{i_1}i_0\equiv i_1\:(mod\:\:n_0)$ and $0 < m_{i_1} < n_0.$
    \begin{itemize}
        \item If $m_{i_1}$ is odd then the required $2l + 1$ equals $n_0 - 2.$ Indeed, $qn_0 - (n_0 - 2)i_0 + m_{i_1}i_0 = qn_0 - (n_0 - 2 - m_{i_1})i_0.$ $0\le n_0 - 2 - m_{i_1} < n_0 - 2$ and $(n_0 - 1 - m_{i_1})$ is even, thus $w_{qn_0 - (n_0 - 2 - m_{i_1})i_0} = 0.$
        \item If $m_{i_1}$ is even then the required $2l + 1$ equals $m_{i_1} - 1.$ Indeed, $qn_0 - (m_{i_1} - 1)i_0 + m_{i_1}i_0 = qn_0 + i_0 = (q + i_0)n_0 - (n_0 - 1)i_0.$ $0 < n_0 - 1 < n_0$ and $n_0 - 1$ is even, thus $w_{(q + i_0)n_0 - (n_0 - 1)i_0} = 0.$
    \end{itemize}
			     Now consider a sequence $\{z_\mathcal{I}\}_{\mathcal{I}\in\ZZ}$ such that:
    \begin{itemize}
        \item $z_{qn_0 - 2ji_0} = 0$ when $0\le 2j\le n_0;$
        \item $z_{2qn_0 - (2j + 1)i_0} = c$ when $0 < 2j + 1 < n_0$ and $l\ne j;$
        \item $z_{2qn_0 - (2l + 1)i_0} \ge 0;$
        \item $z_{(2q+1)n_0 - (2j + 1)i_0} \ge c$ when $0 < 2j + 1 < n_0$ 
    \end{itemize}
    for $q \in \ZZ.$ 

\begin{center}
\begin{tikzpicture}
    \node[orange] at (-1.2, 0.6){Sequence};
    \node[orange] at (-1.2, 1.1){$\{w_{\mathcal{I}}\}_{\mathcal{I}\in\mathbb{Z}}$};
    \node[orange] at (-1.2, 0.1){values:};
    \node[blue] at (-1.2, -0.5){Index:};

	\foreach \i in {0, ...,  12}
		\fill [orange] (\i * 0.8, 0) circle (3pt);
	\node[blue] at (0 * 1.0 , -0.5){0};
	\node[blue] at (2 * 0.8, -0.5){$i_0$};
    	\node[blue] at (3 * 0.8, -0.5){$i_1$};
	\node[blue] at (5 * 0.8, -0.5){$n_0$};
	\node[blue] at (10 * 0.8, -0.5){$i_0 n_0$};

	\node[orange] at (0 * 0.8, 0.5){0};
	\node[orange] at (1 * 0.8, 0.5){0};
	\node[red] at (2 * 0.8, 0.5){0};
	\node[orange] at (3 * 0.8, 0.5){$\widetilde c$};
	\node[red] at (4 * 0.8, 0.5){c};
	\node[orange] at (5 * 0.8, 0.5){0};
	\node[red] at (6 * 0.8, 0.5){0};
	\node[orange] at (7 * 0.8, 0.5){0};
	\node[red] at (8 * 0.8, 0.5){c};
	\node[orange] at (9 * 0.8, 0.5){$\widetilde c$};
	\node[red] at (10 * 0.8, 0.5){0};
	\node[orange] at (11 * 0.8, 0.5){0};
 	\node[orange] at (12 * 0.8, 0.5){0};

    \draw [<->][red] (8 * 0.8 + 0.06, 0.75) -- (10 * 0.8 - 0.06, 0.75);
    \draw [<->][red] (6 * 0.8 + 0.06, 0.75) -- (8 * 0.8 - 0.06, 0.75);
    \draw [<->][red] (4 * 0.8 + 0.06, 0.75) -- (6 * 0.8 - 0.06, 0.75);
    \draw [<->][red] (2 * 0.8 + 0.06, 0.75) -- (4 * 0.8 - 0.06, 0.75);
    \node[red] at (9 * 0.8, 1.0){$i_0$};
    \node[red] at (7 * 0.8, 1.0){$i_0$};
    \node[red] at (5 * 0.8, 1.0){$i_0$};
    \node[red] at (3 * 0.8, 1.0){$i_0$};

    \draw [-][green!75!black] (4 * 0.8 + 0.12, 0.25) -- (4 * 0.8 - 0.12, 0.25);
    \draw [-][green!75!black] (4 * 0.8 + 0.12, 0.20) -- (4 * 0.8 - 0.12, 0.20);

    \draw [green!75!black][dashed][arrows = {-Latex[width'=0pt .5, length=15pt]}] (4 * 0.8, -0.75) -- (4 * 0.8,-2);
    
    \node[orange] at (-1.2, 0.6 - 3.0){Sequence};
    \node[orange] at (-1.2, 1.1 - 3.0){$\{z_{\mathcal{I}}\}_{\mathcal{I}\in\mathbb{Z}}$};
    \node[orange] at (-1.2, 0.1 - 3.0){values:};
    \node[blue] at (-1.2, -0.5 - 3.0){Index:};

	\foreach \i in {0, ...,  12}
		\fill [orange] (\i * 0.8, - 3.0) circle (3pt);
	\node[blue] at (0 * 1.0 , -0.5 - 3.0){0};
	\node[blue] at (2 * 0.8, -0.5 - 3.0){$i_0$};
    	\node[blue] at (3 * 0.8, -0.5 - 3.0){$i_1$};
	\node[blue] at (5 * 0.8, -0.5 - 3.0){$n_0$};
	\node[blue] at (10 * 0.8, -0.5 - 3.0){$i_0 n_0$};

	\node[orange] at (0 * 0.8, 0.5 - 3.0){0};
	\node[orange] at (1 * 0.8, 0.5 - 3.0){0};
	\node[orange] at (2 * 0.8, 0.5 - 3.0){0};
	\node[orange] at (3 * 0.8, 0.5 - 3.0){$\widetilde c$};
	\node[green!75!black] at (4 * 0.8, 0.5 - 3.0){$\widetilde c$};
	\node[orange] at (5 * 0.8, 0.5 - 3.0){0};
	\node[orange] at (6 * 0.8, 0.5 - 3.0){0};
	\node[orange] at (7 * 0.8, 0.5 - 3.0){0};
	\node[orange] at (8 * 0.8, 0.5 - 3.0){c};
	\node[orange] at (9 * 0.8, 0.5 - 3.0){$\widetilde c$};
	\node[orange] at (10 * 0.8, 0.5 - 3.0){0};
	\node[orange] at (11 * 0.8, 0.5 - 3.0){0};
 	\node[orange] at (12 * 0.8, 0.5 - 3.0){0};

\end{tikzpicture}
\end{center}

    This sequence satisfies $a.$ Indeed, 
    \begin{itemize}
    \item For $m = qn_0 - 2ji_0$ we have $\min_{\mathcal{B}\le v\le \mathcal{E}}\{a_v + z_{v + m}\} = 0,$ the minimum is attained at indices $m$ and $m + n_0.$ 
    \item For $m = 2qn_0 - (2j + 1)i_0$ and $j\ne l$ we have $\min_{\mathcal{B}\le v\le \mathcal{E}}\{a_v + z_{v + m}\} = c,$ the minimum is attained at indices $m + i_0$ and $m.$
    \item For $m = 2qn_0 - (2l + 1)i_0$ we have $\min_{\mathcal{B}\le v\le \mathcal{E}}\{a_v + z_{v + m}\} = c,$ the minimum is attained at indices $m + i_0$ and $m + i_1.$
    \item For $m = (2q+1)n_0 - (2j + 1)i_0$ we have $\min_{\mathcal{B}\le v\le \mathcal{E}}\{a_v + z_{v + m}\} = c,$ the minimum is attained at indices $m + i_0$ and $m + n_0.$
    \end{itemize}
    
    Taking finite fragments $(z_1, \dots, z_N)$ with growing $N$ we conclude that $H(\vec a)\ge \frac{\frac{n - 1}{2} + 1}{2n} \ge \frac{1}{4}.$

    $\quad$
    
     \large\textbf{L.II.2}
     \normalsize
     $\:\:$Now assume that  $\mathbf{k > 1.}$

   Define $k_1 := GCD(i_1, n_0).$ W.l.o.g. we can consider that $k_1 \ge k$ (otherwise we can swap $i_0$ and $i_1$).
   
   We will find $k$ different indices $l_1,\dots l_k$ such that $w_{qn_0 + l_j + i_1} = 0$ for all $q\in\ZZ$ and for all $1\le j \le k$ and $l_{j_1}\not\equiv l_{j_2}$ for all $j_1\ne j_2.$  Note, that if $i_1'\equiv i_1\:(mod\:\:n_0)$ and $w_{qn_0 + l_j + i_1'} = 0$ for all $q\in\ZZ$ and for all $1\le j \le k$ then it is true for $i_1.$ Thus, we can assume that $0\le i_1 < n_0.$
    We can represent $i_1$ as $s\cdot k + r,$ where $0 \le r < k.$ We study two different cases:

	$\quad$
    
	\large\textbf{L.II.2.1}
	\normalsize
	$\:\: r = 0.$ 
        
		Denote $n' := \frac{n_0}{k},$  $i_0' = \frac{i_0}{k}$ and $i_1' : = \frac{i_1}{k}.$Consider $a' = (a_j)_{j\equiv 0\:(mod\:\:k)}.$ Sequence $\{w'_\mathcal{I}\}_{\mathcal{I}\in\ZZ}$ is defined as follows:
            \begin{itemize}
                \item $w'_{q'\frac{n_0}{k} - 2j\frac{i_0}{k}} = 0,$ where $0\le 2j < \frac{n_0}{k};$
                \item $w'_{2q'\frac{n_0}{k} - (2j + 1)\frac{i_0}{k}} = c,$ where $0 < 2j + 1 < \frac{n_0}{k};$
                \item $w'_{(2q' + 1)\frac{n_0}{k} - (2j + 1)\frac{i_0}{k}} \ge c,$ where $0 < 2j + 1 < \frac{n_0}{k}$
            \end{itemize}
        for $q'\in\ZZ.$    

 Similar to the previous case $(k = 1)$ we can consider sequence $\{z'_\mathcal{I}\}_{\mathcal{I}\in\ZZ}$ that provides the bound $H(a')\ge\frac{1}{4}$ for $n_0'$, $i_0'$, $i_1'$ and $a'$.
        Now take sequence $\{z_\mathcal{I}\}_{\mathcal{I}\in\ZZ}$ as follows:
        \begin{itemize}
            \item $z_{i\cdot k + r} = z'_i,$ for $0\le i\in\ZZ$ and $0\le r < k.$
        \end{itemize}

\begin{center}
\begin{tikzpicture}
    \node[orange] at (-1.2, 1.1){$\{w'_{\mathcal{I}}\}_{\mathcal{I}\in\mathbb{Z}}$};
    \node[orange] at (-1.2, 0.6){Sequence};
    \node[orange] at (-1.2, 0.1){values:};
    \node[blue] at (-1.2, -0.5){Index:};

	\foreach \i in {0, ...,  7}
		\fill [orange] (0.5 + 3 * \i / 2, 0) circle (3pt);
	\node[blue] at (0.5 , -0.5){0};
	\node[blue] at (2, -0.5){$i_0'$};
    \node[blue] at (3.5, -0.5){$i_1'$};
	\node[blue] at (5, -0.5){$n_0'$};
	\node[blue] at (9.5, -0.5){$2n_0'$};

	\node[orange] at (0.5, 0.5){0};
	\node[orange] at (0.5 + 1.5, 0.5){0};
	\node[orange] at (0.5 + 3, 0.5){$\widetilde c$};
	\node[orange] at (0.5 + 4.5, 0.5){0};
	\node[orange] at (0.5 + 6, 0.5){0};
	\node[orange] at (0.5 + 7.5, 0.5){c};
	\node[orange] at (0.5 + 9, 0.5){0};
 	\node[orange] at (0.5 + 10.5, 0.5){0};

    \draw [-][green!75!black] (8 + 0.12, 0.25) -- (8 - 0.12, 0.25);
    \draw [-][green!75!black] (8 + 0.12, 0.20) -- (8 - 0.12, 0.20);

	\draw [green!75!black][dashed][arrows = {-Latex[width'=0pt .5, length=15pt]}] (8,-0.75) -- (8,-2);
    
    \node[orange] at (-1.2, 1.1 - 3.0){$\{z'_{\mathcal{I}}\}_{\mathcal{I}\in\mathbb{Z}}$};
    \node[orange] at (-1.2, 0.6 - 3.0){Sequence};
    \node[orange] at (-1.2, 0.1 - 3.0){values:};
    \node[blue] at (-1.2, -0.5 - 3.0){Index:};

	\foreach \i in {0, ...,  7}
		\fill [orange] (0.5 + 3 * \i / 2, - 3.0) circle (3pt);
	\node[blue] at (0.5 , -0.5 - 3.0){0};
	\node[blue] at (2, -0.5 - 3.0){$i_0'$};
    \node[blue] at (3.5, -0.5 - 3.0){$i_1'$};
	\node[blue] at (5, -0.5 - 3.0){$n_0'$};
	\node[blue] at (9.5, -0.5 - 3.0){$2n_0'$};

	\node[orange] at (0.5, 0.5 - 3.0){0};
	\node[orange] at (0.5 + 1.5, 0.5 - 3.0){0};
	\node[orange] at (0.5 + 3, 0.5 - 3.0){$\widetilde c$};
	\node[orange] at (0.5 + 4.5, 0.5 - 3.0){0};
	\node[orange] at (0.5 + 6, 0.5 - 3.0){0};
	\node[green!75!black] at (0.5 + 7.5, 0.5 - 3.0){$\widetilde c$};
	\node[orange] at (0.5 + 9, 0.5 - 3.0){0};
 	\node[orange] at (0.5 + 10.5, 0.5 - 3.0){0};

  	\foreach \i in {0, ...,  7}
        \draw [green!75!black][dashed][arrows = {-Latex[width'=0pt .2, length=8pt]}] (0.5 + 1.5 * \i, -0.70 - 3.0) -- (1.5 * \i, -3.4 - 3.0);
    \foreach \i in {0, ...,  7}
        \draw [green!75!black][dashed][arrows = {-Latex[width'=0pt .2, length=8pt]}] (0.5 + 1.5 * \i, -0.70 - 3.0) -- (0.5 + 1.5 * \i, -3.4 - 3.0);
    \foreach \i in {0, ...,  7}
        \draw [green!75!black][dashed][arrows = {-Latex[width'=0pt .2, length=8pt]}] (0.5 + 1.5 * \i, -0.70 - 3.0) -- (1.0 + 1.5 * \i, -3.4 - 3.0);
  
    \node[orange] at (-1.2, 1.1 - 4 - 3){$\{z_{\mathcal{I}}\}_{\mathcal{I}\in\mathbb{Z}}$};
    \node[orange] at (-1.2, 0.6 - 4 - 3){Sequence};
    \node[orange] at (-1.2, 0.1 - 4 - 3){values:};
    \node[blue] at (-1.2, -0.5 - 4 - 3){Index:};

	\foreach \i in {0, ...,  23}
		\fill [orange] (\i / 2, -4 - 3) circle (3pt);
	\node[blue] at (0 , -4.5 - 3){0};
	\node[blue] at (1.5, -4.5 - 3){$i_0$};
    \node[blue] at (3.0, -4.5 - 3){$i_1$};
	\node[blue] at (4.5, -4.5 - 3){$n_0$};
	\node[blue] at (9.0, -4.5 - 3){$2n_0$};

  	\node[orange] at (0, 0.5 - 4 - 3){0};
	\node[orange] at (0.5, 0.5 - 4 - 3){0};
	\node[orange] at (1, 0.5 - 4 - 3){0};
	\node[orange] at (1.5, 0.5 - 4 - 3){0};
	\node[orange] at (2.0, 0.5 - 4 - 3){0};
	\node[orange] at (2.5, 0.5 - 4 - 3){0};
	\node[orange] at (3.0, 0.5 - 4 - 3){$\widetilde c$};
 	\node[orange] at (3.5, 0.5 - 4 - 3){$\widetilde c$};
    \node[orange] at (4.0, 0.5 - 4 - 3){$\widetilde c$};
	\node[orange] at (4.5, 0.5 - 4 - 3){0};
	\node[orange] at (5.0, 0.5 - 4 - 3){0};
	\node[orange] at (5.5, 0.5 - 4 - 3){0};
	\node[orange] at (6.0, 0.5 - 4 - 3){0};
	\node[orange] at (6.5, 0.5 - 4 - 3){0};
	\node[orange] at (7.0, 0.5 - 4 - 3){0};
 	\node[orange] at (7.5, 0.5 - 4 - 3){$\widetilde c$};
    \node[orange] at (8.0, 0.5 - 4 - 3){$\widetilde c$};
	\node[orange] at (8.5, 0.5 - 4 - 3){$\widetilde c$};
	\node[orange] at (9.0, 0.5 - 4 - 3){0};
	\node[orange] at (9.5, 0.5 - 4 - 3){0};
	\node[orange] at (10.0, 0.5 - 4 - 3){0};
	\node[orange] at (10.5, 0.5 - 4 - 3){0};
	\node[orange] at (11.0, 0.5 - 4 - 3){0};
 	\node[orange] at (11.5, 0.5 - 4 - 3){0};
    
\end{tikzpicture}
\end{center}

         This provides us the bound $H(\vec a)\ge\frac{1}{4}$.

	 $\quad$

	\large\textbf{L.II.2.2}
	\normalsize
	$\:\: r\ne 0$
	
		Note, that in this case $k_1 > k$ and thus $s\ne 0$ and $s + 1 \ne\frac{n_0}{k}.$         
		Here we have three different cases:

		$\quad$

	    \large\textbf{L.II.2.2a}
	    \normalsize
	    $\:\: s \equiv \frac{i_0}{k}\:(mod\:\:\frac{n_0}{k}).$
            
            From the proof for $k = 1$ we know that there exists $0 < 2l + 1 < \frac{n_0}{k}$ such that $w'_{q'\frac{n_0}{k} - (2l + 1)\frac{i_0}{k} + (s + 1)} = 0$ for all $q'\in\ZZ.$
            
			Consider $\{z_\mathcal{I}\}_{\mathcal{I}\in Z}$ as follows:
            \begin{itemize}
                \item $z_{qn_0 - 2ji_0 + i},$ where $0\le 2j < n_0$ and $0\le i < k;$
                \item $z_{2qn_0 - (2j + 1)i_0 + i} = c,$ where $0 < 2j + 1 < n_0,$ $j\ne l$ and $0\le i < k;$
                \item $z_{2qn_0 - (2l + 1)i_0 + i}\ge c,$ where $0\le i < k;$
                \item $z_{(2q + 1)n_0 - (2l + 1)i_0 + i},$ where $0 < 2j + 1 < n_0$ and $0\le i < k$
            \end{itemize}
            for $q\in\ZZ.$

\begin{center}
\begin{tikzpicture}
    \node[orange] at (-1.2, 1.1){$\{w'_{\mathcal{I}}\}_{\mathcal{I}\in\mathbb{Z}}$};
    \node[orange] at (-1.2, 0.6){Sequence};
    \node[orange] at (-1.2, 0.1){values:};
    \node[blue] at (-1.2, -0.5){Index:};

	\foreach \i in {0, ...,  7}
		\fill [orange] (0.45 + 3 * \i * 0.45 , 0) circle (3pt);
	\node[blue] at (0.45, -0.5){0};
	\node[blue] at (0.45 * 4, -0.5){$i_0'$};
    \node[blue] at (0.45 * 7, -0.5){$s + 1$};
	\node[blue] at (0.45 * 10, -0.5){$n_0'$};
	\node[blue] at (0.45 * 19, -0.5){$2n_0'$};

	\node[orange] at (0.45 + 0.45 * 0, 0.5){0};
	\node[orange] at (0.45 + 0.45 * 3, 0.5){0};
	\node[orange] at (0.45 + 0.45 * 6, 0.5){$\widetilde c$};
	\node[orange] at (0.45 + 0.45 * 9, 0.5){0};
	\node[orange] at (0.45 + 0.45 * 12, 0.5){0};
	\node[orange] at (0.45 + 0.45 * 15, 0.5){c};
	\node[orange] at (0.45 + 0.45 * 18, 0.5){0};
 	\node[orange] at (0.45 + 0.45 * 21, 0.5){0};

    \draw [-][green!75!black] (0.45 * 16 + 0.12, 0.25) -- (0.45 * 16 - 0.12, 0.25);
    \draw [-][green!75!black] (0.45 * 16 + 0.12, 0.20) -- (0.45 * 16 - 0.12, 0.20);

	\draw [green!75!black][dashed][arrows = {-Latex[width'=0pt .5, length=15pt]}] (0.45 * 16, -0.75) -- (0.45 * 16, -2);

    \node[orange] at (-1.2, 1.1 - 3.0){$\{z'_{\mathcal{I}}\}_{\mathcal{I}\in\mathbb{Z}}$};
    \node[orange] at (-1.2, 0.6 - 3.0){Sequence};
    \node[orange] at (-1.2, 0.1 - 3.0){values:};
    \node[blue] at (-1.2, -0.5 - 3.0){Index:};

	\foreach \i in {0, ...,  7}
		\fill [orange] (0.45 + 3 * \i * 0.45, - 3.0) circle (3pt);
	\node[blue] at (0.45 , -0.5 - 3.0){0};
	\node[blue] at (0.45 * 4, -0.5 - 3.0){$i_0'$};
    \node[blue] at (0.45 * 7, -0.5 - 3.0){$s + 1$};
	\node[blue] at (0.45 * 10, -0.5 - 3.0){$n_0'$};
	\node[blue] at (0.45 * 19, -0.5 - 3.0){$2n_0'$};

	\node[orange] at (0.45 + 0.45 * 0, 0.5 - 3.0){0};
	\node[orange] at (0.45 + 0.45 * 3, 0.5 - 3.0){0};
	\node[orange] at (0.45 + 0.45 * 6, 0.5 - 3.0){$\widetilde c$};
	\node[orange] at (0.45 + 0.45 * 9, 0.5 - 3.0){0};
	\node[orange] at (0.45 + 0.45 * 12, 0.5 - 3.0){0};
	\node[green!75!black] at (0.45 + 0.45 * 15, 0.5 - 3.0){$\widetilde c$};
	\node[orange] at (0.45 + 0.45 * 18, 0.5 - 3.0){0};
 	\node[orange] at (0.45 + 0.45 * 21, 0.5 - 3.0){0};

  	\foreach \i in {0, ...,  7}
        \draw [green!75!black][dashed][arrows = {-Latex[width'=0pt .2, length=8pt]}] (0.45 + 0.45 * 3 * \i, -0.70 - 3.0) -- (0.45 * 3 * \i, -3.4 - 3.0);
    \foreach \i in {0, ...,  7}
        \draw [green!75!black][dashed][arrows = {-Latex[width'=0pt .2, length=8pt]}] (0.45 + 0.45 * 3 * \i, -0.70 - 3.0) -- (0.45 + 0.45 * 3 * \i, -3.4 - 3.0);
    \foreach \i in {0, ...,  7}
        \draw [green!75!black][dashed][arrows = {-Latex[width'=0pt .2, length=8pt]}] (0.45 + 0.45 * 3 * \i, -0.70 - 3.0) -- (0.45 * 2 + 0.45 * 3 * \i, -3.4 - 3.0);
  
    \node[orange] at (-1.2, 1.1 - 4 - 3){$\{z_{\mathcal{I}}\}_{\mathcal{I}\in\mathbb{Z}}$};
    \node[orange] at (-1.2, 0.6 - 4 - 3){Sequence};
    \node[orange] at (-1.2, 0.1 - 4 - 3){values:};
    \node[blue] at (-1.2, -0.5 - 4 - 3){Index:};

	\foreach \i in {0, ...,  23}
		\fill [orange] (\i * 0.45, -4 - 3) circle (3pt);
	\node[blue] at (0 , -4.5 - 3){0};
	\node[blue] at (0.45 * 3, -4.5 - 3){$i_0$};
    \node[blue] at (0.45 * 4, -4.5 - 3){$i_1$};
	\node[blue] at (0.45 * 9, -4.5 - 3){$n_0$};
	\node[blue] at (0.45 * 18, -4.5 - 3){$2n_0$};

  	\node[orange] at (0.45 * 0, 0.5 - 4 - 3){0};
	\node[orange] at (0.45 * 1, 0.5 - 4 - 3){0};
	\node[orange] at (0.45 * 2, 0.5 - 4 - 3){0};
	\node[orange] at (0.45 * 3, 0.5 - 4 - 3){0};
	\node[orange] at (0.45 * 4, 0.5 - 4 - 3){0};
	\node[orange] at (0.45 * 5, 0.5 - 4 - 3){0};
	\node[orange] at (0.45 * 6, 0.5 - 4 - 3){$\widetilde c$};
 	\node[orange] at (0.45 * 7, 0.5 - 4 - 3){$\widetilde c$};
    \node[orange] at (0.45 * 8, 0.5 - 4 - 3){$\widetilde c$};
	\node[orange] at (0.45 * 9, 0.5 - 4 - 3){0};
	\node[orange] at (0.45 * 10, 0.5 - 4 - 3){0};
	\node[orange] at (0.45 * 11, 0.5 - 4 - 3){0};
	\node[orange] at (0.45 * 12, 0.5 - 4 - 3){0};
	\node[orange] at (0.45 * 13, 0.5 - 4 - 3){0};
	\node[orange] at (0.45 * 14, 0.5 - 4 - 3){0};
 	\node[orange] at (0.45 * 15, 0.5 - 4 - 3){$\widetilde c$};
    \node[orange] at (0.45 * 16, 0.5 - 4 - 3){$\widetilde c$};
	\node[orange] at (0.45 * 17, 0.5 - 4 - 3){$\widetilde c$};
	\node[orange] at (0.45 * 18, 0.5 - 4 - 3){0};
	\node[orange] at (0.45 * 19, 0.5 - 4 - 3){0};
	\node[orange] at (0.45 * 20, 0.5 - 4 - 3){0};
	\node[orange] at (0.45 * 21, 0.5 - 4 - 3){0};
	\node[orange] at (0.45 * 22, 0.5 - 4 - 3){0};
 	\node[orange] at (0.45 * 23, 0.5 - 4 - 3){0};
    
\end{tikzpicture}
\end{center}

	We claim that $\{z_\mathcal{I}\}_{\mathcal{I}\in\mathbb{Z}}$ satisfies $a.$ Indeed,
            \begin{itemize}
                \item For $m = qn_0 - 2ji_0 + i$ we have $\min_{\mathcal{B}\le v\le \mathcal{E}}\{a_v + z_{v + m}\} = 0,$ the minimum is attained at indices $m$ and $m + n_0.$ 
                \item For $m = 2qn_0 - (2j + 1)i_0 + i,$ $j\ne l$ we have $\min_{\mathcal{B}\le v\le \mathcal{E}}\{a_v + z_{v + m}\} = c,$ the minimum is attained at indices $m$ and $m + i_0.$
                \item For $m = (2q + 1)n_0 - (2j + 1)i_0 + i,$ $j\ne l$ we have $\min_{\mathcal{B}\le v\le \mathcal{E}}\{a_v + z_{v + m}\} = c,$ the minimum is attained at indices $m + n_0$ and $m + i_0.$
                \item For $m = qn_0 - (2l + 1)i_0 + i$ we have $\min_{\mathcal{B}\le v\le \mathcal{E}}\{a_v + z_{v + m}\} = c,$ the minimum is attained at indices $m + i_1$ and $m + i_0.$
            \end{itemize}
            
             Taking finite fragments $(z_1, \dots, z_N)$ with growing $N$ we conclude that $H(\vec a)\ge \frac{\frac{\frac{n}{k} - 1}{2}\cdot k  + k}{2n} \ge \frac{1}{4}$
            
	    $\quad$

	    \large\textbf{L.II.2.2b}
	    \normalsize
	    $\:\: s + 1 \equiv \frac{i_0}{k}\:(mod\:\:\frac{n_0}{k}).$
            
            This case is the same as the previous one except that we need to find $0 < 2l + 1 < \frac{n_0}{k}$ such that $w'_{q'\frac{n_0}{k} - (2l + 1)\frac{i_0}{k} + s} = 0$ for all $q'\in\ZZ.$
            
\begin{center}
\begin{tikzpicture}
    \node[orange] at (-1.2, 1.1){$\{w'_{\mathcal{I}}\}_{\mathcal{I}\in\mathbb{Z}}$};
    \node[orange] at (-1.2, 0.6){Sequence};
    \node[orange] at (-1.2, 0.1){values:};
    \node[blue] at (-1.2, -0.5){Index:};

	\foreach \i in {0, ...,  10}
		\fill [orange] (0.35 + 0.35 * 3 * \i, 0) circle (3pt);
	\node[blue] at (0.35, -0.5){0};
	\node[blue] at (0.35 + 0.35 * 3 * 4, -0.5){$i_0'$};
    \node[blue] at (0.35 + 0.35 * 3 * 3, -0.5){$s$};
	\node[blue] at (0.35 + 0.35 * 3 * 5, -0.5){$n_0'$};
	\node[blue] at (0.35 + 0.35 * 3 * 10, -0.5){$2n_0'$};

	\node[orange] at (0.35, 0.5){0};
	\node[orange] at (0.35 + 0.35 * 3 * 1, 0.5){c};
	\node[orange] at (0.35 + 0.35 * 3 * 2, 0.5){$\widetilde c$};
	\node[orange] at (0.35 + 0.35 * 3 * 3, 0.5){0};
	\node[orange] at (0.35 + 0.35 * 3 * 4, 0.5){0};
 	\node[orange] at (0.35 + 0.35 * 3 * 5, 0.5){0};
  	\node[orange] at (0.35 + 0.35 * 3 * 6, 0.5){$\widetilde c$};
   	\node[orange] at (0.35 + 0.35 * 3 * 7, 0.5){c};
	\node[orange] at (0.35 + 0.35 * 3 * 8, 0.5){0};
 	\node[orange] at (0.35 + 0.35 * 3 * 9, 0.5){0};
  	\node[orange] at (0.35 + 0.35 * 3 * 10, 0.5){0};

    \draw [-][green!75!black] (0.35 + 0.35 * 3 * 1 + 0.12, 0.25) -- (0.35 + 0.35 * 3 * 1 - 0.12, 0.25);
    \draw [-][green!75!black] (0.35 + 0.35 * 3 * 1 + 0.12, 0.20) -- (0.35 + 0.35 * 3 * 1 - 0.12, 0.20);

    \draw [green!75!black][dashed][arrows = {-Latex[width'=0pt .5, length=15pt]}] (0.35 + 0.35 * 3 * 1,-0.75) -- (0.35 + 0.35 * 3 * 1,-2);
    
    \node[orange] at (-1.2, 1.1 - 3.0){$\{z'_{\mathcal{I}}\}_{\mathcal{I}\in\mathbb{Z}}$};
    \node[orange] at (-1.2, 0.6 - 3.0){Sequence};
    \node[orange] at (-1.2, 0.1 - 3.0){values:};
    \node[blue] at (-1.2, -0.5 - 3.0){Index:};
	\foreach \i in {0, ...,  10}
		\fill [orange] (0.35 +  0.35 * 3 * \i, -3.0) circle (3pt);
	\node[blue] at (0.35 , -0.5 - 3.0){0};
	\node[blue] at (0.35 + 0.35 * 3 * 4, -0.5 - 3.0){$i_0'$};
    \node[blue] at (0.35 + 0.35 * 3 * 3, -0.5 - 3.0){$s$};
	\node[blue] at (0.35 + 0.35 * 3 * 5, -0.5 - 3.0){$n_0'$};
	\node[blue] at (0.35 + 0.35 * 3 * 10, -0.5 - 3.0){$2n_0'$};

	\node[orange] at (0.35, 0.5 - 3.0){0};
	\node[green!75!black] at (0.35 + 0.35 * 3 * 1, 0.5 - 3.0){$\widetilde c$};
	\node[orange] at (0.35 + 0.35 * 3 * 2, 0.5 - 3.0){$\widetilde c$};
	\node[orange] at (0.35 + 0.35 * 3 * 3, 0.5 - 3.0){0};
	\node[orange] at (0.35 + 0.35 * 3 * 4, 0.5 - 3.0){0};
 	\node[orange] at (0.35 + 0.35 * 3 * 5, 0.5 - 3.0){0};
  	\node[orange] at (0.35 + 0.35 * 3 * 6, 0.5 - 3.0){$\widetilde c$};
   	\node[orange] at (0.35 + 0.35 * 3 * 7, 0.5 - 3.0){c};
	\node[orange] at (0.35 + 0.35 * 3 * 8, 0.5 - 3.0){0};
 	\node[orange] at (0.35 + 0.35 * 3 * 9, 0.5 - 3.0){0};
  	\node[orange] at (0.35 + 0.35 * 3 * 10, 0.5 - 3.0){0};

  	\foreach \i in {0, ...,  10}
        \draw [green!75!black][dashed][arrows = {-Latex[width'=0pt .2, length=8pt]}] (0.35 + 0.35 * 3 * \i, -0.70 - 3.0) -- (0.35 * 3 * \i, -3.4 - 3.0);
    \foreach \i in {0, ...,  10}
        \draw [green!75!black][dashed][arrows = {-Latex[width'=0pt .2, length=8pt]}] (0.35 + 0.35 * 3 * \i, -0.70 - 3.0) -- (0.35 + 0.35 * 3 * \i, -3.4 - 3.0);
    \foreach \i in {0, ...,  10}
        \draw [green!75!black][dashed][arrows = {-Latex[width'=0pt .2, length=8pt]}] (0.35 + 0.35 * 3 * \i, -0.70 - 3.0) -- (0.8 + 0.35 * 3 * \i, -3.4 - 3.0);
  
    \node[orange] at (-1.2, 1.1 - 4 - 3){$\{z_{\mathcal{I}}\}_{\mathcal{I}\in\mathbb{Z}}$};
    \node[orange] at (-1.2, 0.6 - 4 - 3){Sequence};
    \node[orange] at (-1.2, 0.1 - 4 - 3){values:};
    \node[blue] at (-1.2, -0.5 - 4 - 3){Index:};

	\foreach \i in {0, ...,  32}
		\fill [orange] (\i * 0.35, -4 - 3) circle (3pt);
	\node[blue] at (0 , -0.5 - 4 - 3){0};
	\node[blue] at (0.35 * 3 * 4, -0.5 - 4 - 3){$i_0$};
    \node[blue] at (0.35 * 3 * 3 + 0.35, -0.5 - 4 - 3){$i_1$};
	\node[blue] at (0.35 * 3 * 5, -0.5 - 4 - 3){$n_0$};
	\node[blue] at (0.35 * 3 * 10, -0.5 - 4 - 3){$2n_0$};
    \draw [-][blue] (0.35 * 3 * 4 + 0.12, 0.25 - 4 - 3) -- (0.35 * 3 * 4 - 0.12, 0.25 - 4 - 3);
    \draw [-][blue] (0.35 * 3 * 4 + 0.12, 0.20 - 4 - 3) -- (0.35 * 3 * 4 - 0.12, 0.20 - 4 - 3);
    \draw [-][blue] (0.35 * 3 * 3 + 0.35 + 0.12, 0.25 - 4 - 3) -- (0.35 * 3 * 3 + 0.35 - 0.12, 0.25 - 4 - 3);
    \draw [-][blue] (0.35 * 3 * 3 + 0.35 + 0.12, 0.20 - 4 - 3) -- (0.35 * 3 * 3 + 0.35 - 0.12, 0.20 - 4 - 3);
    \draw [-][blue] (0.35 * 3 * 5 + 0.12, 0.25 - 4 - 3) -- (0.35 * 3 * 5 - 0.12, 0.25 - 4 - 3);
    \draw [-][blue] (0.35 * 3 * 5 + 0.12, 0.20 - 4 - 3) -- (0.35 * 3 * 5 - 0.12, 0.20 - 4 - 3);
    \draw [-][blue] (0.35 * 3 * 10 + 0.12, 0.25 - 4 - 3) -- (0.35 * 3 * 10 - 0.12, 0.25 - 4 - 3);
    \draw [-][blue] (0.35 * 3 * 10 + 0.12, 0.20 - 4 - 3) -- (0.35 * 3 * 10 - 0.12, 0.20 - 4 - 3);

    \node[orange] at (0.35 * 0, 0.5 - 4 - 3){0};
    \node[orange] at (0.35 * 1, 0.5 - 4 - 3){0};
    \node[orange] at (0.35 * 2, 0.5 - 4 - 3){0};
    \node[orange] at (0.35 * 3, 0.5 - 4 - 3){$\widetilde c$};
    \node[orange] at (0.35 * 4, 0.5 - 4 - 3){$\widetilde c$};
    \node[orange] at (0.35 * 5, 0.5 - 4 - 3){$\widetilde c$};
    \node[orange] at (0.35 * 6, 0.5 - 4 - 3){$\widetilde c$};
    \node[orange] at (0.35 * 7, 0.5 - 4 - 3){$\widetilde c$};
    \node[orange] at (0.35 * 8, 0.5 - 4 - 3){$\widetilde c$};
    \node[orange] at (0.35 * 9, 0.5 - 4 - 3){0};
    \node[orange] at (0.35 * 10, 0.5 - 4 - 3){0};
    \node[orange] at (0.35 * 11, 0.5 - 4 - 3){0};
    \node[orange] at (0.35 * 12, 0.5 - 4 - 3){0};
    \node[orange] at (0.35 * 13, 0.5 - 4 - 3){0};
    \node[orange] at (0.35 * 14, 0.5 - 4 - 3){0};
    \node[orange] at (0.35 * 15, 0.5 - 4 - 3){0};
    \node[orange] at (0.35 * 16, 0.5 - 4 - 3){0};
    \node[orange] at (0.35 * 17, 0.5 - 4 - 3){0};
    \node[orange] at (0.35 * 18, 0.5 - 4 - 3){$\widetilde c$};
    \node[orange] at (0.35 * 19, 0.5 - 4 - 3){$\widetilde c$};
    \node[orange] at (0.35 * 20, 0.5 - 4 - 3){$\widetilde c$};
    \node[orange] at (0.35 * 21, 0.5 - 4 - 3){c};
    \node[orange] at (0.35 * 22, 0.5 - 4 - 3){c};
    \node[orange] at (0.35 * 23, 0.5 - 4 - 3){c};
    \node[orange] at (0.35 * 24, 0.5 - 4 - 3){0};
    \node[orange] at (0.35 * 25, 0.5 - 4 - 3){0};
    \node[orange] at (0.35 * 26, 0.5 - 4 - 3){0};
    \node[orange] at (0.35 * 27, 0.5 - 4 - 3){0};
    \node[orange] at (0.35 * 28, 0.5 - 4 - 3){0};
    \node[orange] at (0.35 * 29, 0.5 - 4 - 3){0};
    \node[orange] at (0.35 * 30, 0.5 - 4 - 3){0};
    \node[orange] at (0.35 * 31, 0.5 - 4 - 3){0};
    \node[orange] at (0.35 * 32, 0.5 - 4 - 3){0};

\end{tikzpicture}
\end{center}

            Taking finite fragments $(z_1, \dots, z_N)$ with growing $N$ we conclude that $H(\vec a)\ge \frac{\frac{\frac{n_0}{k} - 1}{2}\cdot k  + k}{2n_0} \ge \frac{1}{4}$
            
	$\quad$

	    \large\textbf{L.II.2.2c}
	    \normalsize
	    $\:\: s,\:\: s + 1 \not\equiv \frac{i_0}{k}\:(mod\:\:\frac{n_0}{k}).$
            
            From the proof for $k = 1$ we know that there exist $0 < 2l + 1,$  $2l' + 1 < \frac{n_0}{k}$ such that $w'_{q'\frac{n_0}{k} - (2l + 1)\frac{i_0}{k} + s} = 0$ and $w'_{q'\frac{n_0}{k} - (2l' + 1)\frac{i_0}{k} + (s + 1)} = 0$ for all $0\le q'\in\ZZ.$
            
	    Consider $\{z_\mathcal{I}\}_{\mathcal{I}\in \ZZ}$ as follows:
            \begin{itemize}
                \item $z_{qn_0 - 2ji_0 + i},$ where $0\le 2j < n_0$ and $0\le i < k;$
                \item $z_{2qn_0 - (2j + 1)i_0 + i} = c,$ where $0 < 2j + 1 < n_0,$ $j\ne l,\:\:l'$ and $0\le i < k;$
                \item $z_{2qn_0 - (2l + 1)i_0 + i} \ge c,$ where $0 \le i < k - r;$
                \item $z_{2qn_0 - (2l + 1)i_0 + i} = c,$ where $k - r \le i < k;$
                \item $z_{2qn_0 - (2l' + 1)i_0 + i} = c,$ where $0 \le i < k - r;$
                \item $z_{2qn_0 - (2l' + 1)i_0 + i} \ge c,$ where $k - r \le i < k;$
                \item $z_{2qn_0 - (2l + 1)i_0 + i}\ge c,$ where $0\le i < k;$
                \item $z_{(2q + 1)n_0 - (2l + 1)i_0 + i},$ where $0 < 2j + 1 < n_0$ and $0\le i < k$
            \end{itemize}
            for $q\in\ZZ.$ We claim that $\{z_\mathcal{I}\}_{\mathcal{I}\in \ZZ}$ satisfies $a.$
            
            Indeed,
            \begin{itemize}
                \item For $m = qn_0 - 2ji_0 + i$ we have $\min_{\mathcal{B}\le v\le \mathcal{E}}\{a_v + z_{v + m}\} = 0,$ the minimum is attained at indices $m$ and $m + n_0.$ 
                \item For $m = 2qn_0 - (2j + 1)i_0 + i,$ $j\ne l,\:\: l'$ we have $\min_{\mathcal{B}\le v\le \mathcal{E}}\{a_v + z_{v + m}\} = c,$ the minimum is attained at indices $m$ and $m + i_0.$
                \item For $m = 2qn_0 - (2l + 1)i_0 + i,$ $k - r\le i < k$ we have $\min_{\mathcal{B}\le v\le \mathcal{E}}\{a_v + z_{v + m}\} = c,$ the minimum is attained at indices $m$ and $m + i_0.$
                \item For $m = 2qn_0 - (2l' + 1)i_0 + i,$ we have $0\le i < k - r$ $\min_{\mathcal{B}\le v\le \mathcal{E}}\{a_v + z_{v + m}\} = c,$ the minimum is attained at indices $m$ and $m + i_0.$
                \item For $m = (2q + 1)n_0 - (2j + 1)i_0 + i,$ $j\ne l,\:\: l'$ we have $\min_{\mathcal{B}\le v\le \mathcal{E}}\{a_v + z_{v + m}\} = c,$ the minimum is attained at indices $m + n_0$ and $m + i_0.$
                \item For $m = (2q + 1)n_0 - (2l + 1)i_0 + i,$ $k - r\le i < k$ we have $\min_{\mathcal{B}\le v\le \mathcal{E}}\{a_v + z_{v + m}\} = c,$ the minimum is attained at indices $m + n_0$ and $m + i_0.$
                \item For $m = (2q + 1)n_0 - (2l' + 1)i_0 + i,$ $0\le i < k - r$ we have $\min_{\mathcal{B}\le v\le \mathcal{E}}\{a_v + z_{v + m}\} = c,$ the minimum is attained at indices $m + n_0$ and $m + i_0.$
                \item For $m = qn_0 - (2l + 1)i_0 + i,$ $0\le i < k - r$ we have $\min_{\mathcal{B}\le v\le \mathcal{E}}\{a_v + z_{v + m}\} = c,$ the minimum is attained at indices $m + i_1$ and $m + i_0.$
                \item For $m = qn_0 - (2l' + 1)i_0 + i,$ $k - r\le i < k$ we have $\min_{\mathcal{B}\le v\le \mathcal{E}}\{a_v + z_{v + m}\} = c,$ the minimum is attained at indices $m + i_1$ and $m + i_0.$
            \end{itemize}
            
            Taking finite fragments $(z_1, \dots, z_N)$ with growing $N$ we conclude that $H(\vec a)\ge \frac{\frac{\frac{n_0}{k} - 1}{2}\cdot k  + k}{2n_0} \ge \frac{1}{4}.$
$\Box$ \vspace{2mm}

Now we are returning to the proof of the theorem. Assume that there is no such $n\nmid i_i,$ $i_0\ne i_1$  that $a_{i_1} = a_{i_0}.$ As in the proof of  Lemma~\ref{twomins} we will consider two different cases.

$\quad$

    \large\textbf{T.1}\normalsize
    $\:\: k = 1.$
    
	    Define the following sequence $\{w_\mathcal{I}\}_{\mathcal{I}\in\ZZ}$:
    \begin{itemize}
            \item $w_{qn_0 - 2ji_0} = 0$ when $0\le 2j < n_0;$
            \item $w_{2qn_0 - (2j + 1)i_0} = c$ when $0 < 2j + 1 < n_0;$
            \item $w_{(2q+1)n_0 - (2j + 1)i_0} \ge c$ when $0 < 2j + 1 < n_0$
    \end{itemize}
    for $q\in\ZZ.$
    
		We define $\mathcal{L}_0 := \{\mathcal{B}\le v \le \mathcal{E},\:\: n_0\nmid v,\:\: v\ne 0,\:n_0,  $ such that $ w_{qn_0 - (n_0 - 1)i_0 + v} = 0 $ for all  $ q\in\ZZ\}.$ Set $x := \min\:\:\{a_v\:\:|\:\:v\in \mathcal{L}_0\}.$ Also define $i_x$ by the equation $a_{i_x} = x.$ If such $i_x$ is not unique then we choose any $i_x$ with such property. 

		$\quad$

		\large\textbf{T.1.1}\normalsize 
		$\:\:$First assume that $x\le 2c.$
        
		In this case we define a sequence $\{z_\mathcal{I}\}_{\mathcal{I}\in \ZZ}$ as follows:
        \begin{itemize}
            \item $z_{qn_0 - 2ji_0} = 0$ when $0\le 2j < n_0,$ $2j\ne n_0 - 1;$
            \item $z_{2qn_0 - (2j + 1)i_0} = c$ when $0 < 2j + 1 < n_0;$
            \item $z_{(2q+1)n_0 - (2j + 1)i_0} \ge x$ when $0 < 2j + 1 < n_0$
            \item $z_{2qn_0 - (n_0 - 1)i_0} = x;$
            \item $z_{(2q+1)n_0 - (n_0 - 1)i_0} \ge x$
        \end{itemize}
        for $q\in\ZZ.$

\begin{center}
\begin{tikzpicture}
    \node[orange] at (-1.2, 1.1){$\{w_{\mathcal{I}}\}_{\mathcal{I}\in\mathbb{Z}}$};
    \node[orange] at (-1.2, 0.6){Sequence};
    \node[orange] at (-1.2, 0.1){values:};
    \node[blue] at (-1.2, -0.5){Index:};

	\foreach \i in {0, ...,  10}
		\fill [orange] (0.35 +  0.35 * 3 * \i, 0) circle (3pt);
	\node[blue] at (0.35 , -0.5){0};
	\node[blue] at (0.35 + 0.35 * 3 * 2, -0.5){$i_0$};
	\node[blue] at (0.35 + 0.35 * 3 * 5, -0.5){$n_0$};
	\node[blue] at (0.35 + 0.35 * 3 * 10, -0.5){$2n_0$};

	\node[orange] at (0.35, 0.5){0};
	\node[orange] at (0.35 + 0.35 * 3 * 1, 0.5){0};
	\node[orange] at (0.35 + 0.35 * 3 * 2, 0.5){0};
	\node[orange] at (0.35 + 0.35 * 3 * 3, 0.5){$\widetilde{c}$};
	\node[orange] at (0.35 + 0.35 * 3 * 4, 0.5){c};
 	\node[orange] at (0.35 + 0.35 * 3 * 5, 0.5){0};
  	\node[orange] at (0.35 + 0.35 * 3 * 6, 0.5){0};
   	\node[orange] at (0.35 + 0.35 * 3 * 7, 0.5){0};
	\node[orange] at (0.35 + 0.35 * 3 * 8, 0.5){c};
 	\node[orange] at (0.35 + 0.35 * 3 * 9, 0.5){$\widetilde{c}$};
  	\node[orange] at (0.35 + 0.35 * 3 * 10, 0.5){0};

    \draw [green!75!black][dashed][arrows = {-Latex[width'=0pt .5, length=15pt]}] (0.35 + 0.35 * 3 * 3,-0.75) -- (0.35 + 0.35 * 3 * 3,-2);
    \draw [green!75!black][dashed][arrows = {-Latex[width'=0pt .5, length=15pt]}] (0.35 + 0.35 * 3 * 2,-0.75) -- (0.35 + 0.35 * 3 * 2,-2);
    \draw [green!75!black][dashed][arrows = {-Latex[width'=0pt .5, length=15pt]}] (0.35 + 0.35 * 3 * 7,-0.75) -- (0.35 + 0.35 * 3 * 7,-2);
    \draw [green!75!black][dashed][arrows = {-Latex[width'=0pt .5, length=15pt]}] (0.35 + 0.35 * 3 * 9,-0.75) -- (0.35 + 0.35 * 3 * 9,-2);
    
    \node[orange] at (-1.2, 1.1 - 3.0){$\{z_{\mathcal{I}}\}_{\mathcal{I}\in\mathbb{Z}}$};
    \node[orange] at (-1.2, 0.6- 3.0){Sequence};
    \node[orange] at (-1.2, 0.1- 3.0){values:};
    \node[blue] at (-1.2, -0.5- 3.0){Index:};

	\foreach \i in {0, ...,  10}
		\fill [orange] (0.35 + 0.35 * 3 * \i, 0- 3.0) circle (3pt);
	\node[blue] at (0.35 , -0.5- 3.0){0};
	\node[blue] at (0.35 + 0.35 * 3 * 2, -0.5- 3.0){$i_0$};
	\node[blue] at (0.35 + 0.35 * 3 * 5, -0.5- 3.0){$n_0$};
	\node[blue] at (0.35 + 0.35 * 3 * 10, -0.5- 3.0){$2n_0$};

	\node[orange] at (0.35, 0.5- 3.0){0};
	\node[orange] at (0.35 + 0.35 * 3 * 1, 0.5- 3.0){0};
	\node[green!75!black] at (0.35 + 0.35 * 3 * 2, 0.5- 3.0){$x$};
	\node[green!75!black] at (0.35 + 0.35 * 3 * 3, 0.5- 3.0){$\widetilde{x}$};
	\node[orange] at (0.35 + 0.35 * 3 * 4, 0.5- 3.0){c};
 	\node[orange] at (0.35 + 0.35 * 3 * 5, 0.5- 3.0){0};
  	\node[orange] at (0.35 + 0.35 * 3 * 6, 0.5- 3.0){0};
   	\node[green!75!black] at (0.35 + 0.35 * 3 * 7, 0.5- 3.0){$\widetilde{x}$};
	\node[orange] at (0.35 + 0.35 * 3 * 8, 0.5- 3.0){c};
 	\node[green!75!black] at (0.35 + 0.35 * 3 * 9, 0.5- 3.0){$\widetilde{x}$};
  	\node[orange] at (0.35 + 0.35 * 3 * 10, 0.5- 3.0){0};

\end{tikzpicture}
\end{center}

	We claim that $\{z_\mathcal{I}\}_{\mathcal{I}\in \ZZ}$ satisfies $a.$
        
        Indeed,
        \begin{itemize}
            \item For $m = qn_0 - 2ji_0,$  $2j\ne n_0 - 1$ we have $\min_{\mathcal{B}\le v\le \mathcal{E}}\{a_v + z_{v + m}\} = 0,$ the minimum is attained at indices $m$ and $m + n_0.$
            \item For $m = 2qn_0 - (2j + 1)i_0,$ $2j\ne n_0 - 1$ we have $\min_{\mathcal{B}\le v\le \mathcal{E}}\{a_v + z_{v + m}\} = c,$ the minimum is attained at indices $m$ and $m + i_0.$
            \item For $m = (2q+1)n_0 - (2j + 1)i_0$ we have $\min_{\mathcal{B}\le v\le \mathcal{E}}\{a_v + z_{v + m}\} = c,$ the minimum is attained at indices $m + n_0$ and $m + i_0.$
            \item For $m = 2qn_0 - (n - 1)i_0,$  we have $\min_{\mathcal{B}\le v\le \mathcal{E}}\{a_v + z_{v + m}\} = x,$ the minimum is attained at indices $m$ and $m + i_x$ (because $a_v + z_{v + m}$ is at least $x$ if $v\in \mathcal{L}_0$ and $a_v + z_{v + m}\ge min\{c + c, c + x\}\ge x$ if $v\not\in \mathcal{L}_0$).
            \item For $m = (2q+1)n_0 - (n_0 - 1)i_0$ we have $\min_{\mathcal{B}\le v\le \mathcal{E}}\{a_v + z_{v + m}\} = x,$ the minimum is attained at indices $m + n_0$ and $m + i_x$ (because $a_v + z_{v + m}$ is at least $x$ if $v\in \mathcal{L}_0$ and $a_v + z_{v + m}\ge min\{c + c, c + x\}\ge x$ if $v\not\in \mathcal{L}_0$).
        \end{itemize}
        
        Taking finite fragments $(z_1, \dots, z_N)$ with growing $N$ we conclude that $H(\vec a)\ge \frac{\frac{n_0 - 1}{2} + 1}{2n_0} \ge \frac{1}{4}.$

	$\quad$
        
	\large\textbf{T.1.2}\normalsize\
	$\:\:$Now we assume that $x > 2c.$
        
        Denote $\min_{v\ne i_0,\:\: v\nmid n_0}\{a_v\}$ by $s.$ Note, that $s > c.$ Indeed, otherwise we can use lemma \ref{twomins} and get the required bound. Denote $\min_{v\ne 0, n_0,\:\: v\mid n_0}\{a_v\}$ by $d.$ Note, that $d > 0.$ Finally, set $y := \min\{s + c, x, 2c + d\}.$
        
	Define a sequence $\{z_\mathcal{I}\}_{\mathcal{I}\in \ZZ}$ as follows:
        \begin{itemize}
             \item $z_{qn_0 - 2ji_0} = 0$ when $0\le 2j < n_0,$ $2j\ne n_0 - 1;$
            \item $z_{2qn_0 - (2j + 1)i_0} = c$ when $0 < 2j + 1 < n_0;$
            \item $z_{(2q+1)n_0 - (2j + 1)i_0} \ge y$ when $0 < 2j + 1 < n_0$
            \item $z_{4qn_0 - (n_0 - 1)i_0} = 2c;$
            \item $z_{(4q + 1)n_0 - (n_0 - 1)i_0} = t_q,$ is a free variable, $t_q\in[2c, y];$
            \item $z_{(4q + 2)n_0 - (n_0 - 1)i_0} = t_q,$ is a free variable, $t_q\in [2c, y];$
            \item $z_{(4q + 3)n_0 - (n_0 - 1)i_0} = 2c$ 
        \end{itemize}
        for $q\in\ZZ.$
	
\begin{center}
\begin{tikzpicture}
    \node[orange] at (-1.2, 1.1){$\{w_{\mathcal{I}}\}_{\mathcal{I}\in\mathbb{Z}}$};
    \node[orange] at (-1.2, 0.55){Sequence};
    \node[orange] at (-1.2, 0.1){values:};
    \node[blue] at (-1.2, -0.5){Index:};
	\foreach \i in {0, ...,  19}
		\fill [orange] (\i * 0.55, 0) circle (3pt);
	\node[blue] at (0 , -0.5){$2n_0$};
	\node[blue] at (0.55 * 5, -0.5){$3n_0$};
    \node[blue] at (0.55 * 10, -0.5){$4n_0$};
	\node[blue] at (0.55 * 15, -0.5){$5n_0$};

    \node[orange] at (0.55 * 0, 0.5){0};
    \node[orange] at (0.55 * 1, 0.5){0};
    \node[orange] at (0.55 * 2, 0.5){0};
    \node[orange] at (0.55 * 3, 0.5){$\widetilde c$};
    \node[orange] at (0.55 * 4, 0.5){c};
    \node[orange] at (0.55 * 5, 0.5){0};
    \node[orange] at (0.55 * 6, 0.5){0};
    \node[orange] at (0.55 * 7, 0.5){0};
    \node[orange] at (0.55 * 8, 0.5){c};
    \node[orange] at (0.55 * 9, 0.5){$\widetilde c$};
    \node[orange] at (0.55 * 10, 0.5){0};
    \node[orange] at (0.55 * 11, 0.5){0};
    \node[orange] at (0.55 * 12, 0.5){0};
    \node[orange] at (0.55 * 13, 0.5){$\widetilde c$};
    \node[orange] at (0.55 * 14, 0.5){c};
    \node[orange] at (0.55 * 15, 0.5){0};
    \node[orange] at (0.55 * 16, 0.5){0};
    \node[orange] at (0.55 * 17, 0.5){0};
    \node[orange] at (0.55 * 18, 0.5){c};
    \node[orange] at (0.55 * 19, 0.5){$\widetilde c$};
    
    \draw [green!75!black][dashed][arrows = {-Latex[width'=0pt .5, length=15pt]}] (0.55 * 2, -0.75) -- (0.55 * 2, -2);
    \draw [green!75!black][dashed][arrows = {-Latex[width'=0pt .5, length=15pt]}] (0.55 * 3, -0.75) -- (0.55 * 3, -2);
    \draw [green!75!black][dashed][arrows = {-Latex[width'=0pt .5, length=15pt]}] (0.55 * 7, -0.75) -- (0.55 * 7, -2);
    \draw [green!75!black][dashed][arrows = {-Latex[width'=0pt .5, length=15pt]}] (0.55 * 9, -0.75) -- (0.55 * 9, -2);
    \draw [green!75!black][dashed][arrows = {-Latex[width'=0pt .5, length=15pt]}] (0.55 * 12, -0.75) -- (0.55 * 12, -2);
    \draw [green!75!black][dashed][arrows = {-Latex[width'=0pt .5, length=15pt]}] (0.55 * 13, -0.75) -- (0.55 * 13, -2);
    \draw [green!75!black][dashed][arrows = {-Latex[width'=0pt .5, length=15pt]}] (0.55 * 17, -0.75) -- (0.55 * 17, -2);
    \draw [green!75!black][dashed][arrows = {-Latex[width'=0pt .5, length=15pt]}] (0.55 * 19, -0.75) -- (0.55 * 19, -2);

    \node[orange] at (-1.2, 1.1 - 3.0){$\{z_{\mathcal{I}}\}_{\mathcal{I}\in\mathbb{Z}}$};
    \node[orange] at (-1.2, 0.55- 3.0){Sequence};
    \node[orange] at (-1.2, 0.1- 3.0){values:};
    \node[blue] at (-1.2, -0.5- 3.0){Index:};
	\foreach \i in {0, ...,  19}
		\fill [orange] (\i * 0.55, 0- 3.0) circle (3pt);
	\node[blue] at (0 , -0.5- 3.0){$2n_0$};
	\node[blue] at (0.55 * 5, -0.5- 3.0){$3n_0$};
    \node[blue] at (0.55 * 10, -0.5- 3.0){$4n_0$};
	\node[blue] at (0.55 * 15, -0.5- 3.0){$5n_0$};

    \node[orange] at (0.55 * 0, 0.5- 3.0){0};
    \node[orange] at (0.55 * 1, 0.5- 3.0){0};
    \node[green!75!black] at (0.55 * 2, 0.5- 3.0){2c};
    \node[green!75!black] at (0.55 * 3, 0.5- 3.0){$\widetilde y$};
    \node[orange] at (0.55 * 4, 0.5- 3.0){c};
    \node[orange] at (0.55 * 5, 0.5- 3.0){0};
    \node[orange] at (0.55 * 6, 0.5- 3.0){0};
    \node[green!75!black] at (0.55 * 7, 0.5- 3.0){$t_1$};
    \node[orange] at (0.55 * 8, 0.5- 3.0){c};
    \node[green!75!black] at (0.55 * 9, 0.5- 3.0){$\widetilde y$};
    \node[orange] at (0.55 * 10, 0.5- 3.0){0};
    \node[orange] at (0.55 * 11, 0.5- 3.0){0};
    \node[green!75!black] at (0.55 * 12, 0.5- 3.0){$t_1$};
    \node[green!75!black] at (0.55 * 13, 0.5- 3.0){$\widetilde y$};
    \node[orange] at (0.55 * 14, 0.5- 3.0){c};
    \node[orange] at (0.55 * 15, 0.5- 3.0){0};
    \node[orange] at (0.55 * 16, 0.5- 3.0){0};
    \node[green!75!black] at (0.55 * 17, 0.5- 3.0){2c};
    \node[orange] at (0.55 * 18, 0.5- 3.0){c};
    \node[green!75!black] at (0.55 * 19, 0.5- 3.0){$\widetilde y$};

\end{tikzpicture}
\end{center}

	We claim that $\{z_\mathcal{I}\}_{\mathcal{I}\in \Z}$ satisfies $a.$

    Indeed,
        \begin{itemize}
            \item For $m = qn_0 - 2ji_0,$ we have $2j\ne n_0 - 1$ $\min_{\mathcal{B}\le v\le \mathcal{E}}\{a_v + z_{v + m}\} = 0,$ the minimum is attained at indices $m$ and $m + n_0.$
            \item For $m = 2qn_0 - (2j + 1)i_0,$ $2j\ne n_0 - 1$ we have $\min_{\mathcal{B}\le v\le \mathcal{E}}\{a_v + z_{v + m}\} = c,$ the minimum is attained at indices $m$ and $m + i_0.$
            \item For $m = (2q+1)n_0 - (2j + 1)i_0$ we have $\min_{\mathcal{B}\le v\le \mathcal{E}}\{a_v + z_{v + m}\} = c,$ the minimum is attained at indices $m + n_0$ and $m + i_0.$
            \item For $m = 4qn_0 - (n_0 - 1)i_0,$ we have $\min_{\mathcal{B}\le v\le \mathcal{E}}\{a_v + z_{v + m}\} = 2c,$ the minimum is attained at indices $m$ and $m + i_0$ (because $a_v + z_{v + m}$ is at least $x > 2c$ if $v\in \mathcal{L}_0$, and $a_v + z_{v + m}\ge min\{s + c, 0 + t_q, d + 2c\}\ge 2c$ if $v\not\in \mathcal{L}_0$).
            \item For $m = (4q+1)n_0 - (n_0 - 1)i_0$ we have $\min_{\mathcal{B}\le v\le \mathcal{E}}\{a_v + z_{v + m}\} = t_q,$ the minimum is attained at indices $m$ and $m + n_0$ (because $a_v + z_{v + m}$ is at least $x > y\ge t_q$ if $v\in \mathcal{L}_0$, and $a_v + z_{v + m}\ge min\{c + y, s + c, d + 2c\}\ge y \ge t_q$ if $v\not\in \mathcal{L}_0$).
            \item For $m = (4q + 2)n_0 - (n_0 - 1)i_0,$ we have $\min_{\mathcal{B}\le v\le \mathcal{E}}\{a_v + z_{v + m}\} = 2c,$ the minimum is attained at indices $m + n_0$ and $m + i_0$ (because $a_v + z_{v + m}$ is at least $x > 2c$ if $v\in \mathcal{L}_0$, and $a_v + z_{v + m}\ge min\{s + c, 0 + t_q, d + 2c\}\ge 2c$ if $v\not\in \mathcal{L}_0$).
            \item For $m = (4q + 3)n_0 - (n_0 - 1)i_0,$ we have $\min_{\mathcal{B}\le v\le \mathcal{E}}\{a_v + z_{v + m}\} = 2c,$ the minimum is attained at indices $m$ and $m + n_0$ (because $a_v + z_{v + m}$ is at least $x > 2c$ if $v\in \mathcal{L}_0$, and $a_v + z_{v + m}\ge min\{c + y, s + c, 0 + t_q, d + 2c\}\ge 2c$ if $v\not\in \mathcal{L}_0$).
        \end{itemize}
        
         Taking finite fragments $(z_1, \dots, z_N)$ with growing $N$ we conclude that $H(\vec a)\ge \frac{n_0 - 1 + 1}{4n_0} = \frac{1}{4}.$

	 $\quad$
        
	\large\textbf{T.2}\normalsize
	$\:\: k > 1.$
        
	Consider the following sequence $\{z'_\mathcal{I}\}_{\mathcal{I}\in\ZZ}$:
        \begin{itemize}
            \item $z'_{qn_0 - 2ji_0 + i} = 0$ when $0\le 2j\le n_0;$
            \item $z'_{2qn_0 - (2j + 1)i_0 + i} = c$ when $0 < 2j + 1 < n_0;$
            \item $z'_{(2q+1)n_0 - (2j + 1)i_0 + i} \ge c$ when $0 < 2j + 1 < n_0$
        \end{itemize}
    for $q \in \ZZ$ and $0\le i < k.$
    
    For $0\le i < k$ define $mathcal{L}_{0,r} := \{\mathcal{B}\le v \le \mathcal{E},\:\: n_0\nmid v, \:\: v\ne 0, n_0 $ such that $ z'_{qn_0 - (n_0 - 1)i_0 + r + v} = 0 $ for all $q\in\ZZ $ and such that $qn_0 - (n_0 - 1)i_0 + r + v\ne q'n_0 - (n_0 - 1)i_0 + r' $ for any $q'$ and for any $0\le r' < k\}.$ Set $x_r := \min\{a_v\:\: |\:\: v\in \mathcal{L}_{0, r}\}.$ Define $i_{x, r}$ by the equation $a_{i_{x, r}} = x_r.$
    
    $\quad$
    
    Denote $\min_{v\ne i_0,\:\: v\nmid n_0}\{a_v\}$ by $s.$ Note, that $s > c.$ Indeed, otherwise we can use lemma \ref{twomins} and get the required bound. Denote $\min_{v\ne 0, n_0,\:\: v\mid n_0}\{a_v\}$ by $d.$ Note, that $d > 0.$ For $0\le r < k$ set $y_r := min\{s + c, x_r, 2c + d\}.$ Finally, define $M := \max_{0\le r < k}\{x_r, y_r\}.$
    
    $\quad$
    
    Define a sequence $\{z_\mathcal{I}\}_{\mathcal{I}\in\Z}$ as follows:
    \begin{itemize}
        \item[$\bullet$] $z_{qn_0 - 2ji_0 + r} = 0$ when $0\le 2j\le n_0,$ $2j\ne (n_0 - 1),$ where $0\le r < k;$
        \item[$\bullet$] $z_{2qn_0 - (2j + 1)i_0 + r} = c$ when $0 < 2j + 1 < n_0,$ where $0\le r < k;$
        \item[$\bullet$] $z_{(2q+1)n_0 - (2j + 1)i_0 + r} \ge M$ when $0 < 2j + 1 < n_0,$ where $0\le r < k;$
        \item[] For $0\le r < k$ set:
        \begin{enumerate}
            \item if $x_r\le 2c$ then:
            \begin{itemize}
                \item[$\bullet$] $z_{2qn_0 - (n_0 - 1)i_0 + r} = x_r;$
                \item[$\bullet$] $z_{(2q + 1)n_0 - (n_0 - 1)i_0 + r} \ge M;$
            \end{itemize}
            \item if $x_r > 2c$ then:
            \begin{itemize}
                \item[$\bullet$] $z_{4qn_0 - (n_0 - 1)i_0 + r} = 2c;$
                \item[$\bullet$] $z_{(4q + 1)n_0 - (n_0 - 1)i_0 + r} = t_{q, r},$ is a free variable, $t_{q, r}\in [2c, y_r];$
                \item[$\bullet$] $z_{(4q + 2)n_0 - (n_0 - 1)i_0 + r} = t_{q, r},$ is a free variable, $t_{q, r}\in [2c, y_r];$
                \item[$\bullet$] $z_{(4q + 3)n_0 - (n_0 - 1)i_0 + r} = 2c$
            \end{itemize}
        \end{enumerate}
    \end{itemize}
    for $q\in\ZZ.$

\begin{center}
\begin{tikzpicture}
    \node[orange] at (-1.2, 1.1){$\{z_{\mathcal{I}}\}_{\mathcal{I}\in\mathbb{Z}}$};
    \node[orange] at (-1.2, 0.6){Sequence};
    \node[orange] at (-1.2, 0.1){values:};
    \node[blue] at (-1.2, -0.5){Index:};
	\foreach \i in {0, ...,  9}
		\fill [orange] (\i * 1.2, 0) circle (3pt);
	\foreach \i in {0, ...,  9}
		\fill [brown] (0.6 + \i * 1.2, 0) circle (3pt);
	\node[blue] at (0 , -0.5){$2n_0$};
	\node[blue] at (0.6 * 10, -0.5){$3n_0$};

    \node[orange] at (1.2 * 0, 0.5){0};
    \node[orange] at (1.2 * 1, 0.5){0};
    \node[orange] at (1.2 * 2, 0.5){$x_0$};
    \node[orange] at (1.2 * 3, 0.5){$\widetilde M$};
    \node[orange] at (1.2 * 4, 0.5){c};
    \node[orange] at (1.2 * 5, 0.5){0};
    \node[orange] at (1.2 * 6, 0.5){0};
    \node[orange] at (1.2 * 7, 0.5){$\widetilde x_0$};
    \node[orange] at (1.2 * 8, 0.5){c};
    \node[orange] at (1.2 * 9, 0.5){$\widetilde M$};

    \node[brown] at (0.6 + 1.2 * 0, 0.5){0};
    \node[brown] at (0.6 + 1.2 * 1, 0.5){0};
    \node[brown] at (0.6 + 1.2 * 2, 0.5){2c};
    \node[brown] at (0.6 + 1.2 * 3, 0.5){$\widetilde M$};
    \node[brown] at (0.6 + 1.2 * 4, 0.5){c};
    \node[brown] at (0.6 + 1.2 * 5, 0.5){0};
    \node[brown] at (0.6 + 1.2 * 6, 0.5){0};
	\node[brown] at (0.6 + 1.2 * 7, 0.5){$t_{1, 1}$};
    \node[brown] at (0.6 + 1.2 * 8, 0.5){c};
    \node[brown] at (0.6 + 1.2 * 9, 0.5){$\widetilde M$};
	
    \node[blue] at (-1.2, -0.5 - 2.0){Index:};
    \node[blue] at (0.6 * 0, -0.5 - 2.0){$4n_0$};
	\node[blue] at (0.6 * 10, -0.5 - 2.0){$5n_0$};
	\foreach \i in {0, ...,  9}
		\fill [orange] (\i * 1.2, 0 - 2.0) circle (3pt);
	\foreach \i in {0, ...,  9}
		\fill [brown] (0.6 + \i * 1.2, 0 - 2.0) circle (3pt);
    \node[orange] at (1.2 * 0, 0.5 - 2.0){0};
    \node[orange] at (1.2 * 1, 0.5 - 2.0){0};
    \node[orange] at (1.2 * 2, 0.5 - 2.0){$x_0$};
    \node[orange] at (1.2 * 3, 0.5 - 2.0){$\widetilde M$};
    \node[orange] at (1.2 * 4, 0.5 - 2.0){c};
    \node[orange] at (1.2 * 5, 0.5 - 2.0){0};
    \node[orange] at (1.2 * 6, 0.5 - 2.0){0};
    \node[orange] at (1.2 * 7, 0.5 - 2.0){$\widetilde x_0$};
    \node[orange] at (1.2 * 8, 0.5 - 2.0){c};
    \node[orange] at (1.2 * 9, 0.5 - 2.0){$\widetilde M$};

    \node[brown] at (0.6 + 1.2 * 0, 0.5 - 2.0){0};
    \node[brown] at (0.6 + 1.2 * 1, 0.5 - 2.0){0};
	\node[brown] at (0.6 + 1.2 * 2, 0.5 - 2.0){$t_{1, 1}$};
    \node[brown] at (0.6 + 1.2 * 3, 0.5 - 2.0){$\widetilde M$};
    \node[brown] at (0.6 + 1.2 * 4, 0.5 - 2.0){c};
    \node[brown] at (0.6 + 1.2 * 5, 0.5 - 2.0){0};
    \node[brown] at (0.6 + 1.2 * 6, 0.5 - 2.0){0};
    \node[brown] at (0.6 + 1.2 * 7, 0.5 - 2.0){2c};
    \node[brown] at (0.6 + 1.2 * 8, 0.5 - 2.0){c};
    \node[brown] at (0.6 + 1.2 * 9, 0.5 - 2.0){$\widetilde M$};

\end{tikzpicture}
\end{center}

    We claim that this sequence satisfies $a.$ It is sufficient to check that a subsequence $\{z_{qn_0 - (n - 1)i_0 + r'}\}_{q\in\ZZ}$ does not change the minima  in the subsequence $\{z_{qn_0 - (n - 1)i_0 + r}\}_{0\le q\in\ZZ}$ with $r\ne r'$ in the definition of satisfiability of the vector $\vec a$ (see (\ref{-3})). The latter is true because $z_{qn_0 - (n_0 - 1)i_0 + r'}\ge c$ and thus $z_{qn_0 - (n_0 - 1)i_0 + r'} + a_v \ge c + s \ge x_r$ (if $x_r\le 2c$) and $z_{qn_0 - (n_0 - 1)i_0 + r'} + a_v \ge c + s \ge y_r$ (if $x_r > 2c$).
    
    Taking finite fragments $(z_1, \dots, z_N)$ with growing $N$ we conclude that in the worst case $H(\vec a)\ge \frac{(\frac{n_0}{k} - 1)k + k}{4n_0} = \frac{1}{4}.$
$\Box$

\subsection{Sharp upper bound on the tropical entropy in case of a single bounded edge of Newton polygon}

The last theorem is an upper bound on $H(\vec a)$ in case of a single bounded edge of Newton polygon ${\cal N} (\vec a)$. We conjecture that this bound holds for an arbitrary vector $\vec a$. We mention that in \cite{G20} a weaker upper bound $1-1/n$ was established for an arbitrary vector $\vec a$. 
Together with the result $H(\vec a)=1-2/(n+1)$ for a vector $a=(a_0,\dots,a_n)$ with $a_0=\cdots=a_n=0$ \cite[Example 5.2]{G20} it demonstrates the sharpness of the obtained upper bound. The full proof will be provided in the future.

\begin{theorem}\label{verh}
If Newton polygon for $\vec a$ has only one bounded edge then $H(\vec a)\le 1 - \frac{2}{n + 1}.$
\end{theorem}
{\bf Proof}.
For convenience we make a suitable affine transformation such that $a_0 = a_n = 0.$

Consider the polyhedral complex $D(s).$ It is a union of a finite number of polyhedra such that each of these polyhedra $Q$ satisfies the following conditions. For every $0\le j\le s - n$ there exists a pair $0\le i_1 < i_2\le n$ such that 
\begin{equation}\label{restriction}
z_{j + i_1} + a_{i_1} = z_{j} + a_{i_2} = \min_{0\le p \le n}\{z_{p + j} + a_p\}
\end{equation}

for any $(z_1, \dots, z_s)\in Q.$ 

$\quad$

For every $Q$ we consider the following \textit{restriction graph} $RG(Q):$
\begin{itemize}
    \item vertices are the indices of coordinates from $1$ to $s;$
    \item there is an edge between vertices $i$ and $j$ if there is a linear condition of the form $y_i + \gamma = y_j$ which is true for all $(y_1, \dots, y_s)\in Q.$
\end{itemize}
Let us notice that $RG(Q)$ is the union of connected components where each component is the complete subgraph. Moreover, the dimension of $Q$ equals  the number of components of $RG(Q)$ (cf. \cite{G22}).

$\quad$

Let us fix some $Q$ from the finite union above. For arbitrary $(t_1, \dots, t_s)\in Q$ we construct the following sequence by recursion:
\begin{itemize}
    \item The first element of the sequence equals  the least index $i_0$ such that $t_{i_0} = \min_{1\le f \le s}{t_f};$
    \item Let $i_v$ be the last current constructed element of the sequence. If $i_v + n > s$ then we terminate the process and declare $i_v$ to be the last constructed element of the sequence.
    \item If $i_v + n \le s$ then we consider $\min_{0\le p \le n}\{t_{i_v + p} + a_p\}.$ According to the definition of a tropical sequence and the definition of $Q$ there exist $0 \le p_1 < p_2$ such that $\min_{0\le p \le n}\{z_{i_v + p} + a_p\} = z_{i_v + p_1} + a_{p_1} = z_{i_v + p_2} + a_{p_2}$ for all $(z_1, \dots, z_s)\in Q.$ If $p_1 > 0$ then we set $i_{v+1} = i_v + p_1$ and $i_{v+2} = i_v + p_2.$ Otherwise, we just set $i_{v + 1} = i_v + p_2.$ 
    
    Note that there can be more than two indices where $\min_{0\le p \le n}\{z_{i_v + p} + a_p\}$ is attained for all $(z_1, \dots, z_s)\in Q.$ We pick
some pair $p_1 < p_2$.
    
\end{itemize}

We will call this sequence an \textit{equality row} for $(t_1, \dots, t_s).$ Now we claim two important statements:
\begin{itemize}
    \item 

\begin{equation}\label{start_eq_row}
    i_0 < n + 1
\end{equation}
Indeed, suppose the contrary. Then consider $\min_{0\le p \le n}\{t_{i_0 - n + p} + a_p\}.$ As $t_{i_0} = \min_{1\le f\le s}\{t_f\}$ and $a_n = 0$ then this minimum equals  $\min_{1\le f\le s}\{t_f\}$ and there exist $p_1 < p_2\le n$ such that $t_{i_0 - n + p_1} + a_{p_1} = t_{i_0 - n + p_2} + a_{p_2} = \min_{1\le f\le s}\{t_f\}.$ As $a_p\ge 0$ then we obtain that $a_{p_2} = a_{p_1} = 0$ and $t_{i_0 - n + p_1} = t_{i_0 - n + p_2} = t_{i_0}.$ However, $i_0 - n + p_1 < i_0$ and we get a contradiction with that $i_0$ is the least index such that $t_{i_0} = \min_{1\le f\le s}\{t_f\}.$

\item 
\begin{equation}\label{eq_mins}
    t_{i_v} = t_{i_0},
\end{equation}
for all $i_v$ in the equality row.

We prove this by recursion. For $i_0$ the statement is already true. Suppose we have proved this statement for $i_v$ and we consider $\min_{0\le p \le n}\{t_{i_v + p} + a_p\}$ then either $t_{i_v} + a_0 = t_{i_{v + 1}} + a_{p_2}$ equals  this minimum or $t_{i_{v + 1}} + a_{p_1} = t_{i_{v + 2}} + a_{p_2}.$ However, this minimum is less or equal to $t_{i_v} + a_0 = t_{i_v} = \min_{1\le f\le s}\{t_f\}.$ Recalling the fact that $a_p \ge 0$ for $0\le p \le n$ we obtain that $a_{p_1} = a_{p_2} = 0$ and either $t_{i_{v + 1}} = t_{i_v} = t_{i_0}$ or $t_{i_{v + 2}} = t_{i_{v + 1}} = t_{i_v} = t_{i_0}.$
\end{itemize}

$\quad$

Let us fix $(t_1, \dots, t_s)\in Q$ and its equality row  $\{i_0, \dots, i_E\}.$ Consider another arbitrary point $(t_1', \dots, t_s')\in Q.$ We prove the following lemma:
\begin{lemma}\label{eq_row_constr}
If for some $v$ it is true that $t_{i_v}' = \min_{0\le f \le s}\{t'_f\}$ then for all $w \ge v$ it is true that $t'_{i_{w}} = t'_{i_v}.$
\end{lemma}
{\bf Proof of lemma}. 
Indeed, during the recursive construction of the equality row $i_v$ for $(t_1, \dots, t_s)$ there could appear one of the following three possibilities:
\begin{itemize}
    \item $v = 0.$ Then the processes of construction of equality row for $(t_1, \dots, t_s)$ and for $(t_1', \dots, t_s')$ completely coincide.

    \item We considered $\min_{0\le p \le n}\{t_{i_{v - 1} + p} + a_p\}$ which is equal to $t_{i_{v - 1} + p_1} + a_{p_1} = t_{i_{v - 1} + p_2} + a_{p_2}$ for some $p_1 < p_2$ and $i_v = i_{v - 1} + p_2.$ Then the processes of construction of equality row for $(t_1, \dots, t_s)$ and for $(t_1', \dots, t_s')$ completely coincide starting from the next step.
    
    \item We considered $\min_{0\le p \le n}\{t_{i_{v - 1} + p} + a_p\}$ which is equal to $t_{i_{v - 1} + p_1} + a_{p_1} = t_{i_{v - 1} + p_2} + a_{p_2}$ for some $p_1 < p_2$ and $i_v = i_{v - 1} + p_1.$ We recall that these equalities are true for arbitrary $(z_1, \dots, z_s)\in Q$ and so they are true for $(t_1',\dots, t_s').$ Thus  $t'_{i_{v - 1} + p_2}$ also equals  $\min_{1\le f\le s}\{t_f'\}$ and we come to the previous case.
\end{itemize}

$\quad$

Now we define $Q_b$ as
\begin{center}
$\{(y_1, \dots, y_s)\in Q$ and $b$ is the least index such that $y_b = \min_{1\le f \le s}\{t_f\}\}$
\end{center}

According to the statement \ref{start_eq_row} $Q = \bigcup\limits_{b = 1}^n Q_b.$ Next we prove the crucial lemma.
\begin{lemma}\label{eq_end}
The number of connected components in the $RG(Q_b)$ is not greater than $s + 4 - \frac{2s}{n + 1}.$
\end{lemma}
{\bf Proof of lemma}.
According to the definition of $Q_b$ and according to lemma \ref{eq_row_constr} for every $i_v > b$ from the equality row, $(b, i_v)$ is an edge in $RG(Q_b).$ We partition $[b, s]$ into disjoint intervals each of length  $n + 1$ starting from $b.$ Now we produce the following sequence $\{G'_r\}_{r = 0}^{[\frac{s - q}{n + 1}]}$ of subgraphs by recursion on an interval number:
\begin{itemize}
    \item $G'_0$ is just $RG(Q_b)$ without edges;
    \item Suppose we have produced $G'_r$ and now we are considering $(r + 1)$-th interval of length  $(n + 1).$ The interval contains at least one element $i_v$ from the equality row. If there are at least two elements from the equality row then  for each $i_v$ from this interval we add an edge $(b, i_v)$ to the graph $G_r'$ and obtain $G_{r + 1}'.$ 
    
    Otherwise, we consider $(r + 1)$-th interval: 
    $$[b + (n + 1)\cdot r;\:\: b + (n + 1)\cdot r + n].$$
    Consider $\min_{0\le p \le n}\{y_{b + (n + 1)\cdot r + p} + a_p\}$. According to the definition of the tropical sequence there exist $p_1 < p_2$ such that this minimum equals $y_{b + (n + 1)\cdot r + p_1} + a_{p_1} = y_{b + (n + 1)\cdot r + p_2} + a_{p_2}$ for all $(y_1, \dots, y_s)\in Q_b.$ Thus there is an edge from $RG(Q_b)$ whose vertices have indices from the $(r + 1)$-th interval and at least one of them does not lie in the equality row. We call this edge a non-equality edge. Then we set $G'_{r + 1}$ as $G_r'$ with one added edge $(b, i_v)$ and one added non-equality edge.
\end{itemize}

We claim that for every $r$ the number of components in $G'_{r}$ is at least by two less than $G'_{r - 1}.$ It follows from the fact that at each step all edges have at least one end-point which does not belong to the transitive closure of previous subgraph.

Thus we obtain that the number of components is less than $s - 2\cdot [\frac{s - b}{n + 1}] \le s  + 2 - 2\frac{s - b}{n + 1} \le s + 4 - 2\frac{s}{n + 1}.$ 

$\quad$

Now we note that $\dim{Q} = \max_{1\le b \le n}\{\dim{Q_b}\}$ and therefore, according to lemma \ref{eq_end} we obtain that $\dim{Q}\le s + 4 - \frac{2s}{n + 1}.$ Tending to the limit on $s$ we obtain the required statement of the theorem.
$\Box$

{\bf Acknowledgements}. The authors are grateful to an anonymous referee whose remarks fostered to improve the exposition.
The work of the second author was supported in 2021 by a grant in the form of a subsidy from the
federal budget of RF for support of the creation and development of
international mathematical centers, agreenent No. 075-15-2019-1620 of
08.11.2019 between Ministry of Science and Highier Education of RF and
PDMI RAS.

\end{document}